\pdfoutput=1
\documentclass[12pt]{amsart}

\usepackage{amsmath,amsthm,amsfonts,amssymb,anysize,bbm,euscript,stmaryrd,xcolor,mathscinet,mathtools,combelow}

\usepackage[latin1,utf8]{inputenc} 
\usepackage[english]{babel}
\usepackage{tensor}

\SetSymbolFont{stmry}{bold}{U}{stmry}{m}{n}
\usepackage[hmargin=3cm,vmargin=3cm]{geometry}

\usepackage{hyperref}
 \hypersetup{
   colorlinks, linkcolor={blue},
   citecolor={medium-blue}, urlcolor={medium-blue}
 }

\usepackage{palatino, euscript}

\usepackage{graphicx}
\usepackage{microtype}
\usepackage{ifthen}
\usepackage[all]{xy}

\definecolor{dark-red}{rgb}{0.7,0.25,0.25}
\definecolor{dark-blue}{rgb}{0.15,0.15,0.55}
\definecolor{medium-blue}{rgb}{0,0,0.65}
\definecolor{DarkGreen}{RGB}{0,150,0}

\newcommand{\arxiv}[1]{\href{https://arxiv.org/abs/#1}{\small  arXiv:#1}}
\newcommand{\doi}[1]{\href{https://dx.doi.org/#1}{{\small DOI:#1}}}
\newcommand{\euclid}[1]{\href{https://projecteuclid.org/getRecord?id=#1}{{\small  #1}}}
\newcommand{\mathscinet}[1]{\href{https://www.ams.org/mathscinet-getitem?mr=#1}{\small  #1}}
\newcommand{\googlebooks}[1]{(preview at \href{https://books.google.com/books?id=#1}{google books})}

\newcommand{\numdam}[1]{}

\theoremstyle{plain}
\newtheorem{theorem}{Theorem}[section]
\newtheorem{thm-nono}{Theorem}
\newtheorem{proposition}[theorem]{Proposition}
\newtheorem{corollary}[theorem]{Corollary}
\newtheorem{lemma}[theorem]{Lemma}
\newtheorem{conjecture}[theorem]{Conjecture}

\newtheorem{notation}[theorem]{Notation}

\theoremstyle{definition}
\newtheorem{definition}[theorem]{Definition}
\newtheorem{remark}[theorem]{Remark}
\newtheorem{example}[theorem]{Example}

\newcommand{\comment}[1]{}

\newcommand{\C}{\mathbb{C}}
\newcommand{\B}{\mathbb{B}}
\newcommand{\Sch}{\mathbb{S}}

\newcommand{\CC}{\mathcal{C}}

\newcommand{\CK}{\mathcal{K}}
\newcommand{\CP}{\mathcal{P}}

\newcommand{\SBim}{\mathrm{SBim}}
\newcommand{\Hilb}{\mathrm{Hilb}}
\newcommand{\CCC}{\mathbf{C}}
\newcommand{\HH}{\mathrm{HH}}
\newcommand{\subb}{\,\begin{picture}(0,0)(-.75,-1.5)\circle*{2}\end{picture}\ }
\newcommand{\supb}{\,\begin{picture}(0,0)(-.75,-3)\circle*{2}\end{picture}\ }
\newcommand{\End}{\mathrm{End}}
\newcommand{\one}{\mathbbm{1}}

\newcommand{\adds}{\Sigma}

\newcommand{\adddp}{\Sigma^{\Pi}}
\newcommand{\Ai}{A_\infty}
\newcommand{\AS}{\EuScript A}
\newcommand{\BS}{\EuScript B}
\newcommand{\CS}{\EuScript C}
\newcommand{\DS}{\EuScript D}
\newcommand{\IS}{\EuScript I}
\newcommand{\JS}{\EuScript J}
\newcommand{\MS}{\EuScript M}
\newcommand{\NS}{\EuScript N}

\newcommand{\XS}{\EuScript X}
\newcommand{\ZC}{\mathcal{Z}}

\newcommand{\BB}{\mathbf{B}}
\renewcommand{\k}{\mathbbm{k}}
\newcommand{\Z}{\mathbb{Z}}
\newcommand{\Hom}{\operatorname{Hom}}
\newcommand{\Ext}{\operatorname{Ext}}
\newcommand{\im}{\operatorname{im}}
\newcommand{\Sym}{\operatorname{Sym}}
\newcommand{\Br}{\operatorname{Br}}

\renewcommand{\d}{\delta}
\newcommand{\e}{\varepsilon}

\newcommand{\Cone}{\operatorname{Cone}}
\newcommand{\inv}{^{-1}}
\newcommand{\op}{\mathrm{op}}
\renewcommand{\mod}{\text{-mod}}
\newcommand{\rmod}{\text{mod-}}
\newcommand{\dgmod}{\text{-dgmod}}
\newcommand{\dgrmod}{\text{dgmod-}}

\newcommand{\Id}{\operatorname{Id}}
\renewcommand{\a}{\alpha}
\renewcommand{\b}{\beta}

\newcommand{\Tw}{\operatorname{Tw}}
\newcommand{\tw}{\operatorname{tw}}
\newcommand{\Ch}{\operatorname{Ch}}
\newcommand{\Sk}{\operatorname{Sk}}
\newcommand{\smMatrix}[1]{\left[\begin{smallmatrix}#1\end{smallmatrix}\right]}
\newcommand{\sqmatrix}[1]{\left[\begin{matrix} #1\end{matrix}\right]}
\newcommand{\PB}{\mathbf{P}}
\newcommand{\AB}{\mathbf{A}}

\newcommand{\hTr}{\operatorname{Tr}}
\newcommand{\whTr}{\widetilde{\operatorname{Tr}}}
\newcommand{\Tr}{\operatorname{Tr}}
\newcommand{\Trz}{\operatorname{Tr_0}}
\newcommand{\Trzz}{\operatorname{Tr_0}}

\newcommand{\pretr}[1]{\operatorname{Pretr}(#1)}
\newcommand{\pretrf}{\operatorname{Pretr}}
\newcommand{\Kar}{\operatorname{Kar}}
\newcommand{\KhR}{\operatorname{KhR}}
\newcommand{\AKhR}{\operatorname{AKhR}}
\newcommand{\Perf}{\operatorname{Perf}}
\newcommand{\uh}{\underline{h}}
\newcommand{\Vect}{\operatorname{Vect}}
\newcommand{\bbar}{\boldsymbol{|}\hskip-2.5pt\boldsymbol{|}}

\newcommand{\ev}{\mathrm{ev}}
\newcommand{\coev}{\mathrm{coev}}
\newcommand{\susp}[1]{\Sigma^{#1}}
\newcommand{\dsc}[1]{\mathcal{D}_{#1}}
\newcommand{\lon}{L}

\newcommand{\XB}{\mathbf{X}}
\newcommand{\ep}{z}
\newcommand{\ue}{\underline{e}}
\newcommand{\Rid}{\EuScript{R}_{\text{idem}}}

\newcommand{\FT}{\operatorname{FT}}
\newcommand{\Hqe}{\mathbf{Hqe}}
\newcommand{\ee}{\mathbf{e}}

\usepackage{tikz}
\usetikzlibrary{cd}
\usetikzlibrary{calc}
\tikzset{anchorbase/.style={baseline={([yshift=-0.5ex]current bounding box.center)}}}

\usetikzlibrary{matrix, arrows, calc}
\tikzset{frontline/.style={preaction={draw=white,-,line width=6pt}},}  

\newcommand{\Sdot}[3]{
\begin{tikzpicture}[anchorbase,xscale=#2,yscale=#3]
\draw[#1, line width=.5mm] (0,0) to (0,-1);
\fill[#1] (0,0) circle (2.5mm);
\end{tikzpicture}
}

\newcommand{\trival}[3]{
\begin{tikzpicture}[anchorbase,xscale=#2,yscale=#3]
\draw[#1, line width=.5mm] (-1,-1) to (0,0) to (0,1) (0,0) to (1,-1);
\end{tikzpicture}
}

\newcommand{\ann}[1]{
\draw[dashed] (-#1,-2) to (-#1,2) (#1,-2) to (#1,2);
}
\newcommand{\annh}[1]{
\draw[dashed] (-#1,-1) to (-#1,1) (#1,-1) to (#1,1);
}

\newcommand{\sidedot}[3]{
 \begin{tikzpicture}[anchorbase,xscale=#2,yscale=#3]
    \draw[#1, line width=.5mm] 
    (0,-1) to (0,2)
    (1,0) to (0,1);
    \fill[#1] (1,0) circle (2.5mm);
      \end{tikzpicture}
}

\newcommand{\ru}{to [out=0,in=270]}
\newcommand{\rd}{to [out=0,in=90]}
\newcommand{\ur}{to [out=90,in=180]}
\newcommand{\ul}{to [out=90,in=0]}
\newcommand{\lu}{to [out=180,in=270]}
\newcommand{\ld}{to [out=180,in=90]}
\newcommand{\dr}{to [out=270,in=180]}
\newcommand{\dl}{to [out=270,in=0]}
\newcommand{\pu}{to [out=90,in=270]}
\newcommand{\pr}{to [out=0,in=180]}
\newcommand{\pd}{to [out=270,in=90]}
\newcommand{\pl}{to [out=180,in=0]}

\newcommand{\cyltop}[4]{
  \begin{scope}[shift={(#1,#2)}]
\draw (0,0) to (#3,0) \ru (#3+.5,.25) \ul (#3,.5) to (0,.5) \ld (-.5,.25) \dr (0,0);
\draw  (-.5,.25) to (-.5,-#4+.25);
\draw  (#3+.5,.25) to (#3+.5,-#4+.25);
  \end{scope}
}
\newcommand{\cylbot}[3]{
  \begin{scope}[shift={(#1,#2)}]
\draw  (-.5,.25) \dr (0,0) to (#3,0) \ru (#3+.5,.25);
\draw[dashed] (#3+.5,.25) \ul (#3,.5) to (0,.5) \ld (-.5,.25);
  \end{scope}
}

\newcommand{\boxx}[5]{
  \begin{scope}[shift={(#1,#2)}]
\draw[thick, fill=white] (0,0) to (#3,0) to (#3,#4) to (0,#4) to (0,0);
\node at (#3/2,#4/2) {#5};
  \end{scope}
}

\newcommand{\htrfigure}[1]{
  \begin{tikzpicture}[anchorbase,scale=#1]
    \cylbot{0}{0}{3}
    \draw[thick] (2,0) \pu (1.75,1.5);
    \draw[thick] (2.25,2.5) \pu (2,4);
    \draw[thick] (1.75,2.5) \ul (1.5,2.75) to (0,2.75) \lu (-.5,3);
    \draw[thick, dashed] (-.5,3) \ur (0,3.25) \pr (2.25,3.25) \pr (3,1.75) \rd (3.5,1.5);
    \draw[thick] (3.5,1.5) \dl (3,1.25) \pl (2.75,1.25) \lu  (2.25,1.5);
    \boxx{1.5}{1.5}{1}{1}{\small $f$}
    \node at (1.75,.25) {\tiny $X$};
    \node at (1.75,3.75) {\tiny $X'$};
    \node at (3,1) {\tiny $Y_0$};
    \node at (1.5,3) {\tiny $Y_{r+2}$};
    \node at (-.25,3) {\tiny $Y_{0}$};
    \node at (.5,2.5) {\small $\underline{c}$};
    \draw[purple, fill=purple, fill opacity=.2] (0,0) to (0,4) to (1,4) to (1,0) to (0,0);
    \cyltop{0}{4}{3}{4}
  \end{tikzpicture}
\quad
:=
\quad
\begin{tikzpicture}[anchorbase,scale=#1]
  \cylbot{-1}{0}{5}
  \draw[thick] (3,0) \pu (2.75,1.5);
  \draw[thick] (3.25,2.5) \pu (3,4);
  \draw[thick] (2.75,2.5) \ul (2.5,2.75) to (-1,2.75) \lu (-1.5,3);
  \draw[thick, dashed] (-1.5,3) \ur (-1,3.25) \pr (3.25,3.25) \pr (4,1.75) \rd (4.5,1.5);
  \draw[thick] (4.5,1.5) \dl (4,1.25) \pl (3.75,1.25) \lu  (3.25,1.5);
  \boxx{2.5}{1.5}{1}{1}{\small $f$}
  \node at (2.75,.25) {\tiny $X$};
  \node at (2.75,3.75) {\tiny $X'$};
  \node at (4,1) {\tiny $Y_0$};
  \node at (2.5,3) {\tiny $Y_{r+2}$};
  \node at (-1.25,3) {\tiny $Y_{0}$};
  \node at (-1,2.5) {\tiny $c_{0}$};
  \node at (-1,2.75) {$\bullet$};
  \draw[purple] (-.8,0) to (-.8,4);
  \node at (-.6,2.5) {\tiny $c_{1}$};
  \node at (-.6,2.75) {$\bullet$};
  \draw[purple] (-0.4,0) to (-.4,4);
  \node at (-.2,2.5) {\tiny $c_{2}$};
  \node at (-.2,2.75) {$\bullet$};
  \draw[purple] (0,0) to (0,4);
  \node at (0.2,2.75) {$\bullet$};
  \node at (.2,2.5) {\tiny $c_{3}$};
  \node at (0.7,2) {$\cdots$};
  \draw[purple] (1.2,0) to (1.2,4);
  \node at (1.4,2.5) {\tiny $c_{r}$};
  \node at (1.4,2.75) {$\bullet$};
  \draw[purple] (1.6,0) to (1.6,4);
  \node at (2,2.5) {\tiny $c_{r+1}$};
  \node at (1.8,2.75) {$\bullet$};
  
  \cyltop{-1}{4}{5}{4}
\end{tikzpicture}
}

\newcommand{\htrcompfigure}[1]{
  \begin{tikzpicture}[anchorbase,scale=#1]
    \cylbot{0}{1}{3}
    \draw[thick] (2,1) \pu (1.75,1.5);
    \draw[thick] (2.25,2.5) \pu (1.75,3.5);
    \draw[thick] (2.25,4.5) \pu (2,6);
    \draw[thick] (1.75,2.5) \ul (1.5,2.75) to (0,2.75) \lu (-.5,3);
    \draw[thick, dashed] (-.5,3) \ur (0,3.25) \pr (2.25,3.25) \pr (3,1.75) \rd (3.5,1.5);
    \draw[thick] (3.5,1.5) \dl (3,1.25) \pl (2.75,1.25) \lu  (2.25,1.5);
    \boxx{1.5}{1.5}{1}{1}{\small $g$}
    \boxx{1.5}{3.5}{1}{1}{\small $f$}
    \node at (1.5,1.25) {\tiny $X$};
    \node at (.5,2.5) {\small $\underline{d}$};
    \begin{scope}[shift={(0,2)}]
      \draw[thick] (1.75,2.5) \ul (1.5,2.75) to (0,2.75) \lu (-.5,3);
      \draw[thick, dashed] (-.5,3) \ur (0,3.25) \pr (2.25,3.25) \pr (3,1.75) \rd (3.5,1.5);
      \draw[thick] (3.5,1.5) \dl (3,1.25) \pl (2.75,1.25) \lu  (2.25,1.5);
      \node at (.5,2.5) {\small $\underline{c}$};
    \end{scope}
    \node at (1.65,5.75) {\tiny $X''$};
    \draw[purple, fill=purple, fill opacity=.2] (0,1) to (0,6) to (1,6) to (1,1) to (0,1);
    \cyltop{0}{6}{3}{5}
  \end{tikzpicture}
  \quad
  :=
  \quad
  \begin{tikzpicture}[anchorbase,scale=#1]
    \cylbot{0}{0}{3}
    \draw[thick] (2,0) \pu (1.75,1.5);
    \draw[thick] (2.25,2.5) \pu (2,4);
    \draw[thick] (1.75,2.5) \ul (1.5,2.75) to (0,2.75) \lu (-.5,3);
    \draw[thick, dashed] (-.5,3) \ur (0,3.25) \pr (2.25,3.25) \pr (3,1.75) \rd (3.5,1.5);
    \draw[thick] (3.5,1.5) \dl (3,1.25) \pl (2.75,1.25) \lu  (2.25,1.5);
    \boxx{1.125}{1.5}{1.75}{1}{\small $\mu_{\nearrow}(f,g)$}
    \node at (1.75,.25) {\tiny $X$};
    \node at (1.65,3.75) {\tiny $X'$};
    \node at (.5,2.5) {\small $\underline{d}*\underline{c}$};
    \draw[purple, fill=purple, fill opacity=.2] (0,0) to (0,4) to (1,4) to (1,0) to (0,0);
    \cyltop{0}{4}{3}{4}
  \end{tikzpicture}
}

\newcommand{\centonhtrfigure}[1]{
  \begin{tikzpicture}[anchorbase,scale=#1]
    \cylbot{0}{0}{3}
    \draw[thick] (2,0) \pu (1.75,1.5);
    \draw[thick] (2.25,2.5) \pu (2,4);
    \draw[thick] (1.75,2.5) \ul (1.5,2.75) to (0,2.75) \lu (-.5,3);
    \draw[thick, dashed] (-.5,3) \ur (0,3.25) \pr (2.25,3.25) \pr (3,1.75) \rd (3.5,1.5);
    \draw[thick] (3.5,1.5) \dl (3,1.25) \pl (2.75,1.25) \lu  (2.25,1.5);
    \boxx{1.5}{1.5}{1}{1}{\small $e$}
    \node at (1.75,.25) {\tiny $X$};
    \node at (1.65,3.75) {\tiny $X'$};
    \node at (.5,2.5) {\small $\underline{c}$};
    \draw[purple, fill=purple, fill opacity=.2] (0,0) to (0,4) to (1,4) to (1,0) to (0,0);
    \cyltop{0}{4}{3}{4}
  \end{tikzpicture}
  \;\mapsto\;
  \begin{tikzpicture}[anchorbase,scale=#1]
    \cylbot{0}{0}{4}
    \draw[thick] (3,0) \pu (2.75,1.5);
    \draw[thick] (3.25,2.5) \pu (3,4);
    \draw[thick] (2.75,2.5) \ul (2.5,2.75) to (0,2.75) \lu (-.5,3);
    \draw[thick, dashed] (-.5,3) \ur (0,3.25) \pr (3.25,3.25) \pr (4,1.75) \rd (4.5,1.5);
    \draw[thick] (4.5,1.5) \dl (4,1.25) \pl (3.75,1.25) \lu  (3.25,1.5);
    \boxx{2.5}{1.5}{1}{1}{\small $e$}
    \node at (2.75,.25) {\tiny $X$};
    \node at (2.65,3.75) {\tiny $X'$};
    \node at (.5,2.4) {\small $\underline{c}_1$};
    \node at (1.75,2.4) {\small $\underline{c}_2$};
   
    \draw[purple, fill=purple, fill opacity=.2] (0,0) to (0,4) to (1,4) to (1,0) to (0,0);
    \draw[purple, fill=purple, fill opacity=.2] (1.25,0) to (1.25,4) to (2.25,4) to (2.25,0) to (1.25,0);
    \cyltop{0}{4}{4}{4}
  \end{tikzpicture}
  \;\mapsto\;
  \begin{tikzpicture}[anchorbase,scale=#1]
    \cylbot{0}{0}{4}
    \draw[thick] (3,0) \pu (2.75,1.5);
    \draw[thick] (1.5,0) \pu (2.25,2.25) \ul (2,2.5) to (2,3.25) \lu (1.5,4);
    \draw[thick] (3.25,2.5) \pu (3,4);
    \draw[thick] (2.75,2.5) \ul (2.5,2.75) to (0,2.75) \lu (-.5,3);
    \draw[thick, dashed] (-.5,3) \ur (0,3.25) \pr (3.25,3.25) \pr (4,1.75) \rd (4.5,1.5);
    \draw[thick] (4.5,1.5) \dl (4,1.25) \pl (3.75,1.25) \lu  (3.25,1.5);
    \boxx{2.5}{1.5}{1}{1}{\small $e$}
    \node at (2.75,.25) {\tiny $X$};
    \node at (1.25,.25) {\tiny $Z$};
    \node at (2.65,3.75) {\tiny $X'$};
    \node at (1.25,3.75) {\tiny $Z$};
    \node at (.5,2.4) {\small $\underline{c}_1$};
    \begin{scope}[shift={(2,2)}]
      \begin{scope}[rotate=60]
        \boxx{0}{0}{1.25}{1}{\small $\tau(\underline{c}_2)$}
      \end{scope}
    \end{scope}
    \draw[purple, fill=purple, fill opacity=.2] (0,0) to (0,4) to (1,4) to (1,0) to (0,0);
    \cyltop{0}{4}{4}{4}
  \end{tikzpicture}
}

\newcommand{\traciator}[1]{
  w_{X,Y} \;=\; 
  \begin{tikzpicture}[anchorbase,scale=#1]
    \cylbot{0}{0}{3}
    \draw[thick] (1.25,0) \pu (1,1.5) to (1,2.5);
    \draw[thick] (1.75,0) \pu (1.5,1.5) to (1.5,2.5);
    \draw[thick] (2,1.5) to (2,2.5) \pu (1.75,4);
    \draw[thick] (1.5,2.5) \pu (1.25,4);
    \draw[thick] (1,2.5) \ul (.75,2.75) to (0,2.75) \lu (-.5,3);
    \draw[thick, dashed] (-.5,3) \ur (0,3.25) \pr (2.25,3.25) \pr (3,1.75) \rd (3.5,1.5);
    \draw[thick] (3.5,1.5) \dl (3,1.25) \pl (2.25,1.25) \lu  (2,1.5);
    \node at (1,.25) {\tiny $X$};
    \node at (2,.25) {\tiny $Y$};
    \node at (1,3.75) {\tiny $Y$};
    \node at (2,3.75) {\tiny $X$};
    \cyltop{0}{4}{3}{4}
  \end{tikzpicture}
  \quad,\quad
  w^{-1}_{X,Y} \;=\; 
  \begin{tikzpicture}[anchorbase,scale=#1]
    \cylbot{0}{0}{3}
    \draw[thick] (1.75,0) \pu (2,2.5);
    \draw[thick] (1.25,0) \pu (1.75,4);
    \draw[thick] (1.25,4) \pd (1,1.25) \dl (.75,1) \lu (.5,1.25) \ul (0,1.75) \lu (-.5,2);
    \draw[thick, dashed] (-.5,2) \ur (0,2.25) \pr (2.5,2.25) \pr (3,1.75) \rd (3.5,1.5);
    \draw[thick] (3.5,1.5) \dl (3,1.25) \pl (2.75,1.25) \lu  (2.5,1.5) to (2.5,2.5) \ul (2.25,2.75) \ld (2,2.5);
    \node at (1,.25) {\tiny $Y$};
    \node at (2,.25) {\tiny $X$};
    \node at (0,1.4) {\tiny ${}^*X$};
    \node at (1,3.75) {\tiny $X$};
    \node at (2,3.75) {\tiny $Y$};
    \cyltop{0}{4}{3}{4}
  \end{tikzpicture}
  \quad,\quad
  w_{X} \;=\; 
  \begin{tikzpicture}[anchorbase,scale=#1]
    \cylbot{0}{0}{2}
    \draw[thick] (1,0) \pu (.75,1.5) to (.75,2.5) \ul (.5,2.75) to (0,2.75) \lu (-.5,3);
    \draw[thick, dashed] (-.5,3) \ur (0,3.25) \pr (1.25,3.25) \pr (2,1.75) \rd (2.5,1.5);
    \draw[thick] (2.5,1.5) \dl (2,1.25) \pl (1.5,1.25) \lu  (1.25,1.5) to (1.25,2.5) \pu (1,4);
    \node at (.75,.25) {\tiny $X$};
    \node at (.75,3.75) {\tiny $X$};
    \cyltop{0}{4}{2}{4}
  \end{tikzpicture}
}

\title{Derived traces of Soergel categories}
\author{Eugene Gorsky}
\author{Matthew Hogancamp}
\author{Paul Wedrich}

\address{E.G.: Department of Mathematics, 
University of California at Davis, 
One Shields Avenue, Davis CA 95616}
\address{E. G.: Moscow State University, 
Faculty of Mathematics and Mechanics, 
GSP-1, Moscow 119991, Russia}
\email{egorskiy@math.ucdavis.edu}

\address{M.H.: Department of Mathematics,
Northeastern University,
360 Huntington Ave, Boston, MA 02115}
\email{m.hogancamp@northeastern.edu}

\address{P.W.: Max Planck Institute for Mathematics,
Vivatsgasse 7, 53111 Bonn, Germany} 
\address{P.W.: Mathematical Institute, University of Bonn,
Endenicher Allee 60, 53115 Bonn, Germany}
\address{P.W.: Mathematical Sciences Research Institute,
17 Gauss Way, Berkeley, CA 94720. 
\href{http://paul.wedrich.at}{paul.wedrich.at}}
\email{p.wedrich@gmail.com}

\begin{document}

\begin{abstract} We study two kinds of categorical traces of (monoidal) dg
categories, with particular interest in categories of Soergel bimodules.  First,
we explicitly compute the usual Hochschild homology, or derived vertical trace,
of the category of Soergel bimodules in arbitrary types. Secondly, we introduce
the notion of derived horizontal trace of a monoidal dg category and compute the
derived horizontal trace of Soergel bimodules in type $A$. As an application we
obtain a derived annular Khovanov--Rozansky link invariant with an action of
full twist insertion, and thus a categorification of the HOMFLY-PT skein module
of the solid torus.
\end{abstract}

\maketitle
\setcounter{tocdepth}{1}
\tableofcontents

\section{Introduction} 
\subsection{Traces}
Traces are ubiquitous in mathematics. If $A$ is an algebra over a field $\k$,
its trace (or cocenter) is defined as the quotient $\HH_0(A)=A/[A,A]$. Given any
finite-dimensional $A$-module $M$, the trace $\mathrm{tr}_{M}\colon A\to \k$ of
the $A$-action on $M$ satisfies $\mathrm{tr}_{M}(xy)=\mathrm{tr}_{M}(yx)$ and
hence factors through $A/[A,A]$. The projection $A\to A/[A,A]$ can thus be
considered as a universal trace on $A$. 

We are interested in traces for categories. If $\CS$ is a $\k$-linear category,
then its trace is a vector space over $\k$. If $\CS$ is a monoidal category (or
a 2-category) then its trace is a 1-category. Abstractly, the trace of an
$n$-category satisfying certain assumptions is an $(n-1)$-category related to
the factorization homology of the circle (see e.g. \cite{2019arXiv190310961A}
for an introduction), but we will not employ this point of view in the present
paper. Categorical traces have been studied in various settings, see for example
Ocneanu~\cite{MR1317353}, Evans--Kawahigashi~\cite{MR1316301},
Walker~\cite{kw:tqft}, Ben-Zvi--Nadler \cite{2009arXiv0904.1247B},
Ponto--Shulman~\cite{MR3095324}, Beliakova--Lauda--Habiro--\v{Z}ivkovi\'{c}~
\cite{MR3606445}, Hoyois--Scher\-otzke--Sibilla~\cite{MR3607274},
Beliakova--Putyra--Wehrli~\cite{MR3910068}.

In this paper we construct and study
the derived traces of monoidal dg categories, with a view towards applications
in higher representation theory and link homology.

Given a $\k$-linear dg category $\CS$, one can define its Hochschild
homology $\HH_{\subb}(\CS)$ which is a vector space over $\k$. It is defined as
homology of the explicitly defined {\em cyclic bar complex} which we review in
Section \ref{ss:vert trace}. Keller proved \cite{Keller} that $\HH_{\subb}(\CS)$ is
a derived invariant of $\CS$, 
for example, if $\CS$ is the category of perfect
complexes over an algebra $A$ then $\HH_{\subb}(\CS)$ is isomorphic to the usual
Hochschild homology of $A$. Note that $\HH_0(A) = A/[A,A]$ suggesting the
interpretation of Hochschild homology as a \emph{derived (vertical) trace}.

If $\CS$ is a monoidal category (or a bicategory), then there is a richer notion
of {\em horizontal trace} $\Trzz(\CS)$ which is well-studied in various levels
of generality. This is a category equipped with the ``trace functor''
$\Trzz\colon \CS\to \Trzz(\CS)$. If $\CS$ has left duals, then $\Trzz$ is
initial among all trace-like functors\footnote{Variations of trace-like functors
are known under the names
\emph{shadows} \cite{MR3095324}, \emph{commutator functors}
\cite{MR2506324}, \emph{categorical traces} \cite{MR3578212}, \emph{trace functors} \cite{MR3837875}.}
$F\colon \CS \to \DS$, i.e. functors equipped with a natural transformation
$F(X\otimes Y)\to F(Y \otimes X)$ respecting the tensor product in $\CS$ (if
$\CS$ also has right duals, then the components of these natural transformations are necessarily
isomorphisms). The horizontal trace is indeed a richer notion than the vertical
trace, since the endomorphism algebra of $\Trzz(\one)$ in the horizontal trace
naturally agrees with the vertical trace $\HH_0(\CS)$.

\begin{example} The horizontal trace of the bicategory of tangles (where objects
are finite sets of points in $I^2$, $1$-morphisms are tangles in $I^3$, and
$2$-morphisms are tangle cobordisms in $I^4$ up to isotopy rel boundary) is
the category of links in the thickened annulus (objects are links in $I^2 \times
S^1$, morphisms are link cobordisms in $I^3 \times S^1$ up to isotopy rel boundary).
\end{example}

In this paper we define a derived
version of the horizontal trace and prove the following:

\begin{theorem}
\label{thm:trace functor}
There is a natural dg functor $\CS\to \hTr(\CS)$, which is \emph{homotopy
trace-like} i.e. it is equipped with transformations $\hTr(X\otimes Y)\to \hTr(Y\otimes X)$ that are
natural in $Y$ and natural up to coherent homotopy in $X$. The endomorphism algebra of $\hTr(\one)$ is
naturally isomorphic to $\HH_{\subb}(\CS)$.
\end{theorem}
If $\CS$ has left duals we expect that $\hTr$ is initial among all homotopy
trace-like dg functors out of $\CS$. 

We also define the notion of dg Drinfeld
center of $\CS$ and prove that it acts on $\hTr(\CS)$; see section
\ref{sec:center-on-trace}. 

It is desirable to consider the closure of $\hTr(\CS)$ with respect to mapping
cones and homotopy direct summands, which we denote by:
\begin{equation}
\label{eq:completed trace}
\whTr(\CS):=\pretr{\Kar^{dg}(\hTr(\CS))},
\end{equation}
where $\Kar^{dg}$ and $\pretr{-}$ respectively denote the homotopy idempotent
completion and pretriangulated hull; see sections \ref{sec: pretriangulated dg}
and \ref{sec:Karoubi} for details. 

\subsection{Traces of Soergel bimodules}

Next, we apply all the machinery of derived traces to categories of Soergel
bimodules, starting with a computation of the derived vertical trace. Let $W$ be
a Coxeter group with simple reflections $S\subset W$ and a realization $V$ over
$\C$, and $\SBim(W)$ the associated monoidal category of Soergel bimodules
\cite{MR1173115}, which is a categorification of the Hecke algebra associated to the Coxeter
system $(W,S)$. 
We set $R:=\Sym^\bullet(V^\ast)$, graded by placing $V^\ast$ in bidegree $(2,0)$ and
$\Lambda:=\Lambda^\bullet(V^\ast)$, graded by placing $V^\ast$ in bidegree $(2,-1)$.

\begin{theorem}
    \label{thm:main}
    We have an isomorphism of associative bigraded algebras
    $$
    \HH_{\subb}(\SBim(W))\ \cong \  \HH_{\subb}(R)\rtimes  \C[W].
     $$
\end{theorem}

\begin{remark}\label{rmk:HHR}
The Hochschild homology of the polynomial ring $R$ is canonically isomorphic to
$R\otimes \Lambda$.  After choosing a basis of $V$ we can identify
$\HH_{\subb}(\SBim(W))$ with the algebra
$\C[x_1,\ldots,x_r,\theta_1,\ldots,\theta_r]\rtimes  \C[W]$ in which $x_i$ are
even variables of degree $(2,0)$ and the $\theta_i$ are odd variables of degree
$(2,-1)$, and $r=\dim(V)$.
\end{remark}

 \begin{remark}
In Theorem \ref{thm:main} the generators of the wreath product algebra on the
right hand side are identified with the Hochschild cycles on the left hand side
as follows. The generators of $R$ correspond to cycles $x\in \Hom(\one,\one)$,
the generators of  $\Lambda$ correspond to the cycles $x\bbar \Id-\Id\bbar
x\in \Hom(\one,\one)\otimes \Hom(\one,\one)$ (see section~\ref{sec:bardg} for notation conventions), and the elements $w\in W$
correspond to $(-1)^{\ell(w)}\Id_{\Delta_w}\in \Hom(\Delta_w,\Delta_w)$ where
$\Delta_w$ is the Rouquier complex corresponding to positive braid lift of $w$,
and $\ell(w)$ is the length of $w$. See Theorem \ref{thm:HHSBim} for details. 
\end{remark} 
 
\begin{remark}
In \cite{MR3448186} the isomorphism $\HH_0(\SBim(W))\cong R\rtimes \C[W]$ was proved by
completely different methods, using cellularity of $\SBim(W)$.
\end{remark}

\begin{conjecture}
\label{conj: HH formal}
 $\CCC(\SBim(W))$ is formal as dg algebra, so higher $\Ai$-operations on
 Hochschild homology vanish.
\end{conjecture}

To support this conjecture, we prove a closely related Theorem \ref{thm: K formal} stating that $\End_{\hTr}(\Tr(K))$ is formal as a dg algebra, where $K$ is a certain ``cube complex'' built out of several copies of $\one$. Note that the conjecture describes formality of $\End_{\hTr}(\Tr(\one))$.

After proving Theorem~\ref{thm:trace functor} we specialise to Soergel bimodules
$\SBim_n$ for the symmetric group $S_n$, which feature in triply-graded
Khovanov--Rozansky link homology~\cite{MR2339573}. In \cite{GW} the first and
third authors studied the category of annular webs and foams, which can be
regarded as Karoubi completion of the horizontal trace $\Trzz(\SBim_n)$. In
particular, they proved that this Karoubi completion is generated by the direct
summands of $\Trzz(\one_n)$. Here we generalise this result to the derived
horizontal trace using slightly different methods. 

\begin{theorem}
  \label{mainthm:hTrSBimn}
The dg functor $\Hom_{\Tr(\SBim_n)}(\Tr(\one_n),-)$ induces a  quasi-equivalence
relating 
$\whTr(\SBim_n)$ to the category of perfect right $\Ai$-modules
over $\End(\hTr(\one_n))$.  In other words, we have a quasi-equivalence
\begin{eqnarray*}
\whTr(\SBim_n) & \cong &  \Perf(\HH_{\subb}(R)\rtimes \C[S_n])\\
&\cong & \Perf(\C[x_1,\ldots,x_n,\theta_1,\ldots,\theta_n]\rtimes \C[S_n]),
\end{eqnarray*}
where $\deg(x_i)=(2,0)$ and $\deg(\theta_i)=(2,-1)$.
\end{theorem}

Note that projective modules over the algebra $\HH_{\subb}(R)\rtimes \C[S_n]$ are
naturally indexed up to isomorphism by partitions of $n$, and under the
equivalence of categories above every object of $\whTr(\SBim_n)$ can be
expressed as a twisted complex whose terms are (direct sums of shifts of) these
projective modules. In other words, if $\ee_{\lambda}$ is a projector in
$\C[S_n]$ onto an irreducible representation $V_{\lambda}$, then we define the
\emph{Schur object}
\begin{equation}
\label{eq: def schur}
\Sch^{\lambda}:=\ee_{\lambda}\hTr(\one_n)
\end{equation}
in $\whTr(\SBim_n)$. A perfect $\Ai$-module over $\HH_{\subb}(R)\rtimes \C[S_n]$ 
is then a twisted complex built out of $\Sch^{\lambda}$.

\begin{remark} The type A Soergel bimodule categories taken together form the
  monoidal bicategory $\bigoplus_{n\geq 0} \SBim_n$, with the (new) tensor
  product $\boxtimes$ provided by induction $\SBim_n \times \SBim_m \to
  \SBim_{m+n}$. Its trace inherits the monoidal structure and is expected to
  admit a braiding that induces the $S_n$ action on $\Tr(\one_1)^{\boxtimes n} =
  \Tr(\one_n)$. The Schur object $\Sch^{\lambda}$ is designed to be the
  evaluation of the $\lambda$-Schur functor on $\Tr(\one_1)$.
  
\end{remark}
 
\newcommand{\annu}{\mathbb{A}}
\subsection{Derived annular link invariants and categorification of the skein module of the solid torus}
The main motivation for this paper is to develop a framework for the categorification
of the HOMFLY-PT skein module of the solid torus that is
compatible with expectations from topological field theory, while at
the same time allowing for explicit computations of the associated link invariant. 

To describe this skein module, recall that the type $A_{n-1}$ Hecke algebra $H_n$ can
be described as the linear span of braids on $n$ strands modulo \emph{skein
relations} and isotopies. The multiplication in $H_n$ is inherited from stacking
braids, and the unit is represented by identity braid. 

Similarly, the (positive half of the) skein module of the annulus $\Sk^{+}_n(\annu)$ is
defined as the linear span of annular braid closures modulo skein relations and
isotopies. It is easy to see from the definition that $\Sk^{+}_n(\annu)$ is
isomorphic to the cocenter of the Hecke algebra:
$$
\Sk^{+}_n(\annu)\simeq \frac{H_n}{[H_n,H_n]}.
$$ 
Any {\em trace function} $f$ on $H_n$, i.e. a linear function satisfying
$f(xy)=f(yx)$, naturally factors through the cocenter, and hence can be viewed as
a function on $\Sk^{+}_n(\annu)$.

Let $\Lambda_q$ denote the ring of symmetric functions in infinitely many
variables over $\C(q)$, and let $\Lambda_q^{(n)}$ denote
the subspace of degree $n$ symmetric functions.  The skein module of the annulus
enjoys the following properties:
\begin{itemize}
\item[(a)]$\Sk^{+}_n(\annu)$ is isomorphic to $\Lambda_q^{(n)}$.
 It has a basis of Schur functions $s_{\lambda}$ labeled by
partitions $\lambda$ with $n$ boxes.
\item[(b)] The HOMFLY-PT invariant of links yields a trace on $H_n\rightarrow \C(q)[a,a\inv]$
(called the Jones-Ocneanu trace), and can be computed by projecting $H_n$ to
$\Sk^{+}_n(\annu)\cong \Lambda_q^{(n)}$ and applying a certain algebra map $\Lambda_q\rightarrow \C(q)[a,a\inv]$.
\item[(c)] The center of $H_n$ naturally acts on its cocenter. On the level of
annular link diagrams this corresponds to cutting open the annular link diagram and inserting a
central element before closing it again. In particular, the
full twist is central in the braid group and hence acts on $\Sk^{+}_n(\annu)$. 
\end{itemize}

Links in $S^1 \times D^2$ may be studied by means of their diagrams in $\annu$,
after choosing a homeomorphism $S^1\times D^2 \cong \annu\times I$.  Such a
homeomorphism will be referred to as an \emph{I-bundle structure} on $S^1\times
D^2$.  A \emph{framing} is a choice of $I$-bundle structure up to isotopy.  A
choice of $I$-bundle structure gives us a well-defined link diagram associated
to generic links $L$, whereas a choice of framing determines a diagram only up
to Reidemeister moves.  Two different framings are related by some number of
twists, which on the level of link diagrams corresponds to the insertion of some
power of the full-twist braid, as in (c) above.

For this reason, if one is interested in (say, the positive half) of the skein
module of a 3-manifold $Y$ which is homeomorphic to $S^1\times D^2$ (but with
no preferred homeomorphism) then it is necessary to understand not just the
skein module $\Sk^{+}_n(\annu)$, but also the automorphism of full twist insertion.

The categorification of the skein module proceeds in several steps. First, the Hecke
algebra $H_n$ is categorified by the monoidal category $\SBim_n$ of Soergel
bimodules in type $A_{n-1}$, or by a closely related monoidal category of webs
and foams defined by Queffelec--Rose~\cite[Remark 3.24]{MR3545951}, see \cite{MR2671770},
\cite[Remark 3.3]{MR3982970}, and \cite[Section 4.4]{MR3590355} for the
connection.  In the second step one must categorify the cocenter of $H_n$.
Traditionally (see Queffelec--Rose~\cite{MR3729501},
Beliakova--Putyra--Wehrli~\cite{MR3910068}, and
Queffelec--Rose--Sartori~\cite{QRS}) this is done using the
underived horizontal trace $\hTr_0$.   This underived horizontal trace is
satisfactory for many purposes.  For instance Queffelec--Rose--Sartori proved in
\cite{QRS} that the triply-graded Khovanov-Rozansky homology
$\KhR$ \cite{MR2421131,MR2339573} factors through the underived horizontal
trace, which gives a categorification of (b) above.  Additionally, in \cite{GW}
the first and third authors connected the annular Khovanov--Rozansky invariant
of Queffelec--Rose~\cite{MR3729501} to the underived horizontal trace of type A
Soergel bimodules $\hTr_0(\SBim_n)$ and showed that a categorification of (a)
holds upon Karoubi completion. 

Note that $\hTr_0(\SBim_n)$ indeed categorifies $\Lambda_q^{(n)}$, as
it is generated by objects $\Sch^{\lambda}$ as in \eqref{eq: def schur}, which correspond to Schur
functions $s_{\lambda}$. So the Grothendieck group of $\hTr_0(\SBim_n)$ is naturally isomorphic to
$\Lambda_q^{(n)}$.

For a categorification of the skein module of the solid torus, we also need
automorphisms of the target category which realise changes in $I$-bundle
structure, as in property (c) above. The following example shows that the
ordinary annular Khovanov--Rozansky link invariant, which is constructed using the
underived horizontal trace, does not enjoy this property. 

\begin{example} The annular Khovanov--Rozansky invariant of a 2-component unlink
decomposes into two non-trivial direct summands $\AKhR(\one_2)\cong S^2\oplus
\wedge^2$. Here $S^2$ and $\wedge^2$ denote the Schur objects $\Sch^{(2)}$ and
$\Sch^{(1,1)}$ as in \eqref{eq: def schur}. Twisting the $I$-bundle structure
turns the unlink into an annular Hopf link---the braid closure of the full twist
on two strands---whose invariant is a chain complex 
\[
    \begin{tikzpicture}[baseline=0em]
    \tikzstyle{every node}=[font=\scriptsize]
    \node (a) at (0,0) {$\AKhR(\FT_2)$};
    \node at (1.25,0) {$\simeq $};
    \node (cb) at (3,0) {$\underline{\wedge^2(-2)\oplus\wedge^2} \;\; \oplus \;\;\Big(\underline{0}$};
    \node (cc) at (7,0) {$\wedge^2$};
    \node (cd) at (11,0) {$S^2(2)\Big)$};
    \path[->,>=stealth,shorten >=1pt,auto,node distance=1.8cm,
      thick]
    (cc) edge node {$x_1-x_2$} (cd);
    \path[->,>=stealth,shorten >=1pt,auto,node distance=1.8cm]
    (cb) edge (cc);
    \end{tikzpicture}.
    \] 
which decomposes into three
  non-trivial direct summands. Changes in $I$-bundle structure do not induce
  isomorphisms on annular Khovanov--Rozansky invariant. Here we have used the
  version of $\AKhR$ defined in \cite{GW}, but the same
  argument applies to all other constructions employing the horizontal trace. 
\end{example}

To remedy this issue, we use the derived horizontal trace.

\begin{definition} We define the \emph{derived annular Khovanov--Rozansky
  invariant} of a braid word $\underline{\b}$ on $n$ strands, denoted
  $\AKhR_{\mathrm{dg}}(\underline{\b})$, to be the derived horizontal trace
  class of the Rouquier complex of $\underline{\b}$ in $\Tr(\SBim_n)$.
\end{definition}

By Theorem~\ref{mainthm:hTrSBimn}, this invariant can be considered as taking
values in perfect $\Ai$ modules over the $\Ai$ algebra $\C[x_1,\ldots,x_n,
\theta_1,\ldots, \theta_n] \rtimes \C[S_n], $ where $S_n$ is supported in
cohomological degree zero and the variables $x_i$ and $\theta_i$ have
cohomological degree $0$ and $-1$ respectively. Unlike for $\AKhR$, changes in
$I$-bundle structure induce automorphisms on $\AKhR_{\mathrm{dg}}$. These arise
naturally through the action of the derived central Rouquier complex of the full
twist braid on the derived horizontal trace. 

\begin{example} The derived annular Khovanov--Rozansky invariant of the full
twist on two strands is a twisted complex
    \[
\begin{tikzpicture}[baseline=0em]
\tikzstyle{every node}=[font=\scriptsize]
\node (a) at (0,0) {$\AKhR_{\mathrm{dg}}(\FT_2)$};
\node at (1.4,0) {$\simeq $};
\node (cb) at (3,0) {$ \underline{\wedge^2(-2)} \;\; \oplus\;\; \Big(\underline{\wedge^2}$};
\node (cc) at (7,0) {$\wedge^2$};
\node (cd) at (11,0) {$S^2(2)\Big)$};
\draw[frontline,->,>=stealth,shorten >=1pt,auto,below,node distance=1.8cm,thick]
(cb) to[bend right=10] node {$\theta_1-\theta_2$} (cd);
\path[->,>=stealth,shorten >=1pt,auto,node distance=1.8cm,
  thick]
(cc) edge node {$x_1-x_2$} (cd);
\end{tikzpicture}.
\]
 Like $\AKhR_{\mathrm{dg}}(\one_2)\cong S^2\oplus \wedge^2$, this has two
indecomposable direct summands. The action of the derived central full twist on
$\AKhR_{\mathrm{dg}}$ sends $\wedge^2$ to $\wedge^2(-2)$ and $S^2$ to the twisted
complex shown as the second direct summand above.
\end{example}

\begin{remark}
The indecomposable summands of $\AKhR_{\mathrm{dg}}(\FT_n)$ are nothing but the
images of indecomposable summands of $\AKhR_{\mathrm{dg}}(\one_n)$ (that is,
$\Sch^{\lambda}$ for partitions $\lambda$ of $n$) under the action of the full
twist. Following the conjectures of the first author, Negu\cb{t} and Rasmussen
\cite{GNR}, we expect the action of the full twist to be closely related to the
action of Bergeron-Garsia operator $\nabla$ originating in the theory of Macdonald
polynomials \cite{nabla}. 
\end{remark}

In future work, we will use $\AKhR_{\mathrm{dg}}$ to study \emph{cabling
operations} for Khovanov--Rozansky link homologies. We also anticipate that the
technology of derived traces will be useful in the program to categorify skein
algebras, see e.g. \cite{1806.03416}, and for explicit computations of the vector space-valued
4-manifold invariants derived from Khovanov--Rozansky link homology
\cite{2019arXiv190712194M}.

\smallskip

\subsection{Comparison with character sheaves}
\label{subsec:character sheaves}

In this section we briefly compare our results to the theory of character
sheaves. Let $G$ be a semi-simple split algebraic group with a Borel subgroup
$B$ and Weyl group $W$. Recall that a geometric categorification of the Hecke
algebra for $W$ is given by the category of $B$-equivariant constructible
sheaves on the flag variety $G/B$, or, equivalently, $B$-biequivariant sheaves
on $G$. Other, very similar versions of the Hecke category include
$B$-equivariant $D$-modules on $G/B$, or Harish-Chandra bimodules. By the work
of Soergel the geometric Hecke category is closely related to the category of
Soergel bimodules. The polynomial ring R corresponds to the B-equivariant
cohomology of a point. For more details, see \cite{MR1029692}.

In \cite{MR2669705,2009arXiv0904.1247B} both the trace and the center of the
geometric Hecke category were identified with the category of Lusztig's
character sheaves \cite{MR732546}. In particular, the object $\Tr(\one)$ which
plays a prominent role in the paper corresponds to the so-called Springer
sheaf. The derived endomorphism algebra of the Springer sheaf is known (in
particular, it is isomorphic to $\C[W]$ in degree zero), and the formality
result similar to Theorem \ref{thm:main} was proved by Rider \cite{MR3010057}.
It is important to mention that the results of \cite{MR3010057} hold in the
category of {\em mixed perverse sheaves} which is equipped with an additional
grading, which is analogous to our $q$-grading. See also \cite[eq. 0.0.4]{MR3985610} and \cite{Li} for related results.
On the other hand, in Theorem \ref{mainthm:hTrSBimn} we get wo sets of variables $\theta_i,x_i$ while the other references 
have only one set of variables.

Finally, it is known \cite{MR732546,MR3644235} that in type A the summands of
the Springer sheaf generate the category of character sheaves, while this is not
the case in other types. 

Note that we do not claim any results about the Drinfeld center of $\SBim_n$ or
the corresponding category of complexes, but plan to compute it in the future
work.

\subsection{Comparison with Hilbert scheme of points}

In \cite{GNR} the first author, Negu\cb{t} and Rasmussen proposed a set of
conjectures relating the category of Soergel bimodules to the Hilbert scheme of
points on the plane $\Hilb^n(\C^2)$. In particular, they conjectured that both
the trace and the center of $\CK^b(\SBim_n)$ are closely related to the derived
category of coherent sheaves on $\Hilb^n(\C^2)$. We plan to work out the precise
connection between this work and \cite{GNR} in the future, and only comment on
one remarkable formal similarity. 

Haiman constructed in \cite{MR1839919} a rank $n!$ vector bundle $\CP$ on
$\Hilb^n(\C^2)$ called \emph{the Procesi bundle}. Its endomorphism algebra has
the form
\begin{equation}
\label{eq: procesi}
\Hom(\CP,\CP)=\C[x_1,\ldots,x_n,y_1,\ldots,y_n]\rtimes \C[S_n],\ \Ext^i(\CP,\CP)=0,\ i>0.
\end{equation}
It is known that the direct summands of $\CP$ generate the derived category of
$\Hilb^n(\C^2)$, and hence the functor
$$
\mathrm{RHom}(\CP,-)\colon D^b(\Hilb^n(\C^2))\to D^b(\C[x_1,\ldots,x_n,y_1,\ldots,y_n]\rtimes \C[S_n]\mod)
$$
is an equivalence \cite{MR1824990}. The equation \eqref{eq: procesi} is very
similar to the endomorphism algebra of $\Tr(\one_n)$ appearing in Theorem
\ref{thm:main}, but the odd variables $\theta_i$ of degree $(2,-1)$ are replaced
by the even variables $y_i$ of degree $(-2,2)$. It is likely that the dg
enhancement of $D^b(\Hilb^n(\C^2))$ is related to the horizontal trace of
$\SBim_n$ by some kind of Koszul duality.

\subsection{Organization of the paper}

In section \ref{sec: facts} we set up notation and conventions for differential
graded (dg) categories and functors between them. Throughout the paper we chose
to avoid any discussion of $\Ai$-categories and $\Ai$-functors, so we use the
formalism of quasi-functors instead (see subsection \ref{subsec:dg fun}). We
also review the notion of formality for dg algebras and its relation to Massey
products, see subsection \ref{subsec:Ainf}. 

In section \ref{sec: pretriangulated dg} and \ref{sec:Karoubi} we discuss various notions of
completion of dg categories with respect to direct sums, cones and homotopy
idempotents. In particular, we define the pretriangulated hull and the dg Karoubi
completion for an abstract dg category. This material is quite standard, and can
be found, for example, in Seidel's book \cite{MR2441780}, which we more or less
follow; however, we decided to present it as concretely as possible for the
readers' convenience. In particular, we avoid Yoneda embeddings altogether and
explicitly construct $\Ai$ lifts of homotopy idempotents (Proposition
\ref{prop:hi-Aii}) which allows us to give a dg model for the Karoubi
completion. 

In section \ref{sec:bardg} we define and study the 2-sided bar complex of a dg category,
its cyclic version and Hochschild homology of a dg category. In subsection
\ref{ss:vert trace} we compare the full Hochschild homology of the dg category
with its vertical trace. The main result of this section is Theorem
\ref{thm:cyclicbarcx-retract} where we prove that if a dg category admits a
semiorthogonal decomposition then its cyclic bar complex retracts onto the
direct sum of cyclic bar complexes for summands. This is a dg version of a
result of Kuznetsov \cite{Kuz} on additivity of Hochschild homology in
semiorthogonal decompositions.

The next section \ref{s:center and trace} is the technical core of the paper. We
define and study the derived Drinfeld center and derived horizontal trace for
monoidal dg categories. We prove Theorem \ref{thm:trace functor} on the
universal trace functor $\hTr\colon \CS\to \hTr(\CS)$ and its properties, and
also define an action of the derived center on the derived trace. Note that
$\hTr(\CS)$ is usually not pretriangulated or idempotent complete, but the
results of section \ref{sec:bardg} allow us to consider the corresponding
completions.

In section \ref{sec:HH-Scat} we apply this machinery to the monoidal category of
Soergel bimodules and prove Theorem \ref{thm:main}. 

In section \ref{sec:type A} we prove Theorem \ref{mainthm:hTrSBimn} and describe
an explicit ``annular simplification'' algorithm which allows us to identify the
trace of any type A Soergel bimodule with a homotopy summand in the direct sum
of several copies of $\Tr(\one)$.  We also discuss the connections of the
derived trace with annular Khovanov-Rozansky invariants and the work of the first
and third author \cite{GW}. In particular, we construct a ``forgetful functor''
from the derived to the ``underived'' horizontal trace, and show that
Khovanov-Rozansky homology of a braid closure factors through it, see Theorem
\ref{thm: hh from hooks complex}.

\section*{Acknowledgements}

The authors would like to thank Anna Beliakova, Roman Bezrukavnikov, Tudor
Dimofte, Ben Elias, Mikhail Gorsky, Tina Kanstrup, Bernhard Keller, Oscar
Kivinen, Andrei Negu\cb{t},  Alexei Oblomkov, Jacob Rasmussen, David Reutter,
Rapha\"{e}l Rouquier, Lev Rozansky, Sarah Scherotzke, Catharina Stroppel, Kostya
Tolmachov and Geordie Williamson for many useful discussions. 

Special thanks to Kostya Tolmachov who found a gap in a proof of Conjecture
\ref{conj: HH formal} that was claimed in the first version of this paper.

This work also profited from discussions during the AIM workshop ``Categorified
Hecke algebras, link homology, and Hilbert schemes'', the ESI Workshop
``Categorification in quantum topology and beyond'', the workshop ``Hidden
Algebraic Structures in Topology'' at Caltech, the workshop ``Hilbert schemes,
categorification and combinatorics'' at UC Davis, and the MSRI programs
``Quantum Symmetries'' and ``Higher Categories and Categorification''. We thank
the host institutions and sponsors of these activities for their support.

\section*{Funding}
Section \ref{sec:type A} was elaborated by E.G. with the support by the Grant
16-11-10018 of the Russian Science Foundation. The work of E.G. in the other
Sections was partially supported by NSF grants DMS-1700814 and DMS-1760329.

M. H. was partially supported by the NSF grant DMS-1702274.

P.W. was partially supported by the Australian Research Council grants `Braid
groups and higher representation theory' DP140103821 and `Low dimensional
categories' DP160103479 while at the Australian National University during early
stages of this project.

\section{Facts from homological algebra}
\label{sec: facts}
This section serves to recall important notions from homological algebra and to fix
notation and conventions. 

\subsection{Complexes}
\label{sec:DGdef}
Let $\k$ be a commutative ring. The category of complexes of $\k$-modules will
be denoted $\Ch(\k)$. Its objects are complexes of $\k$-modules, also called
\emph{dg $\k$-modules}, and we will use the cohomological convention for
differentials:
\[\cdots \xrightarrow{d} X^k \xrightarrow{d} X^{k+1} \xrightarrow{d} \cdots \]
 
In particular, the differentials are considered to be of cohomological degree
$1$. The morphism spaces between objects $X$ and $Y$ in $\Ch(\k)$ are the
complexes with
\[
  \Hom_{\Ch(\k)}^k(X,Y)=\prod_{i\in \Z} \Hom_{\k}(X^i,Y^{i+k})
  , \qquad 
  d_{\Hom_{\Ch(\k)}(X,Y)}(f) := d_Y \circ f - (-1)^{|f|} f \circ d_X
  \]
where $|f|$ denotes the cohomological degree of $f$. The full subcategory of bounded complexes will be denoted $\Ch^b(\k)$.
The categories $\Ch(\k)$ and $\Ch^b(\k)$ are
symmetric monoidal, with the tensor product defined on objects as

\[(X\otimes_{\k} Y)^k = \bigoplus_{i+j=k} X^i \otimes_\k Y^j, \qquad
d_{X\otimes_{\k} Y} = d_X\otimes \Id_Y + \Id_X \otimes d_Y \] 
and on morphisms
$f$, $g$ by $(f\otimes g)(x\otimes y) = (-1)^{|x||g|}f(x)\otimes g(y)$, with
braiding 
\[\tau_{X,Y}\colon X\otimes_{\k} Y \to  Y \otimes_{\k} X, \qquad \tau_{X,Y}(x\otimes y)= (-1)^{|x||y|} y \otimes x.\]

Finally, for a complex $X$ and $l\in \Z$, we denote by $\susp{l}X$ the complex with 
\[
  (\susp{l}X)^k=X^{k+l},\qquad d_{\susp{l}X}=(-1)^l d_X.
\]
In particular, for $l>0$, the translation $\susp{l}$ shifts the complex $X$ to the
left by $l$ steps.

\subsection{Differential graded categories}
A differential $\Z$-graded $\k$-linear category $\CS$, or short a \emph{dg
category} is a category enriched in $\Ch(\k)$. This means every morphism space
$\Hom_{\CS}(X,Y)$ is an object in $\Ch(\k)$ and composition of morphisms forms chain maps
\[
 \Hom_{\CS}(Y,Z) \otimes_{\k} \Hom_{\CS}(X,Y) \to \Hom_{\CS}(X,Z) 
\]
which means that the differentials satisfy the Leibniz rule with respect to composition:
$d_\CS(f \circ g) = d_\CS(f) \circ g + (-1)^{|f|} f \circ d_\CS(g)$.
Other abelian groups besides $\Z$ can be used for gradings, and later we will
consider $\Z\times \Z$-graded complexes with differentials of degree $(0,1)$. In
any case, the cohomological degree of a homogenous morphism $f$ will be denoted
by $|f|$.

\begin{example}
Any $\k$-linear category can and will be regarded as a dg category with morphism
complexes concentrated in cohomological degree zero, thus necessarily with
zero differential.
\end{example}

\begin{example}
For any $\k$-linear category $\AS$, the category $\Ch(\AS)$ (resp. $\Ch^b(\AS)$) of (bounded) complexes in
$\AS$ (whose definition mimics the one of $\Ch(\k)$ (resp. $\Ch^b(\k)$)) is a dg category.
\end{example}

\begin{example}
Any (dg) $\k$-algebra $A$ can and will be regarded as a dg category with one
object $\ast$ and $\End_{A}(\ast)=A$.
\end{example}

A morphism $f\in \Hom_{\CS}(X,Y)$ is said to be \emph{closed} if $d_{\CS}(f)=0$
and \emph{exact} or \emph{null-homotopic} if $f=d_{\CS}(h)$ for some $h\in
\Hom_{\CS}(X,Y)$, which in this case is called a \emph{null-homotopy} for $f$.
For $f,g\in \Hom^k_{\CS}(X,Y)$, we write $f\simeq g$ and say $f$ and $g$ are
\emph{homotopic} if $f-g$ is null-homotopic. 

The \emph{cohomology category} of $\CS$, denoted $H^0(\CS)$ is defined to be the
additive category with the same objects as $\CS$, and with 
\[\Hom_{H^0(\CS)}(X,Y):= \frac{\{f \in \Hom^0_{\CS}(X,Y)| d_{\CS}(f)=0\}}{d_{\CS}(\Hom^{-1}_{\CS}(X,Y)) }\]

By \emph{isomorphism} in $\CS$ we mean degree zero closed invertible morphisms.
If there exists an isomorphism in $\Hom_{\CS}(X,Y)$, we write $X\cong Y$. A
degree zero closed morphism $f \in \Hom_{\CS}(X,Y)$ is said to be a
\emph{homotopy equivalence} if it induces an isomorphism in $H^0(\CS)$, and in
this case we write $X \simeq Y$ and say $X$ and $Y$ are homotopy equivalent. If
$X\simeq 0$, then we say $X$ is \emph{contractible}.

For a dg category $\CS$, we denote by $\CS^{\op}$ the dg category with the same
objects as $\CS$ and with $\Hom_{\CS^{\op}}(X,Y):= \Hom_{\CS}(Y,X)$, where the composite of
$f\in \Hom_{\CS^{\op}}(X,Y)$ with $g\in \Hom_{\CS^{\op}}(Y,Z)$, denoted $g\circ_{\CS^{\op}} f$, is given by $(-1)^{|f||g|} f
\circ_\CS g$.  

For two dg categories $\CS,\DS$ we denote by $\CS
\otimes_\k \DS$ the category with objects given by pairs of objects $(X,Y)$ for
$X\in \CS$ and $Y\in \DS$ and morphisms given by complexes
\[
\Hom_{\CS\otimes_\k \DS}\Big((X,Y),(X',Y')\Big):= \Hom_{\CS}(X,X')\otimes_\k \Hom_{\DS}(Y,Y')
\]
with composition
\[
(f\otimes g)\circ (f'\otimes g') := (-1)^{|g||f'|}(f\circ f')\otimes (g\circ g').
\]

Let $\CS$ be a dg category.  A \emph{subcategory} $\IS\subset \CS$ is a
collection of objects $\operatorname{Obj}(\IS)\subset \operatorname{Obj}(\CS)$
with hom spaces $\Hom_{\IS}(X, X')$ for $X,X'\in \operatorname{Obj}(\IS)$ being
subcomplexes $\Hom_{\CS}(X, X')$, which are closed under composition.  The subcategory
$\IS\subset \CS$ is \emph{full} if $\Hom_{\IS}(X,X') = \Hom_{\CS}(X,X')$ for all $X,X'\in
\operatorname{Obj}(\IS)$.  The subcategory $\IS\subset \CS$ is \emph{unital} if
$\Id_X\in \IS$ whenever $X\in \operatorname{Obj}(\IS)$. Henceforth all
subcategories are unital.

\begin{example}\label{exa:subcat-idmaps}
We will denote by $\IS\subset \CS$ the subcategory spanned by the identity maps in
$\CS$.  More generally, if $\BS\subset \CS$ is a full subcategory, then we have
$\IS_\BS\subset \CS$, the (unital, but not full) subcategory spanned by the
identity morphisms in $\CS$ which are contained in $\BS$.
\end{example}

\subsection{DG functors}
\label{subsec:dg fun}
If $\BS,\CS$ are dg categories, a dg functor $F\colon \BS\rightarrow \CS$ is a
functor whose action on hom complexes $\Hom_\BS(X,Y)\rightarrow
\Hom_\CS(F(X),F(Y))$ is a degree zero chain map.  The collection of dg functors
$\BS\rightarrow \CS$ itself forms a dg category. Objects of this functor
category are functors, and morphisms are natural transformations, as defined
next. If $F,G$ are functors $\BS\rightarrow \CS$, a natural transformation
$\a\colon F\rightarrow G$ of degree $k$ is an assignment $X\mapsto \a_X \in
\Hom^k_\CS(F(X),G(X))$ such that
\[
G(f) \circ \a_X = (-1)^{k |f|} \a_Y \circ F(f)
\]
for all morphisms $f\in \Hom_\BS(X,Y)$. The differential of $\a$ by definition
sends $X\mapsto d_{\CS}(\a_X)$ (the naturality of $d(\a_X)$ so defined follows
from the Leibniz rule).

By an \emph{isomorphism of dg functors} we mean a degree zero closed invertible
natural transformation of dg functors. If $F,G\colon \BS \rightarrow \CS$ are
isomorphic, we write $F\cong G$.

Any dg functor $F\colon \BS\rightarrow \CS$ naturally induces a functor between
the corresponding homotopy categories $H^0(F)\colon H^0(\BS)\rightarrow
H^0(\CS)$.

  A dg functor $F\colon \CS\rightarrow \DS$ is an \emph{equivalence} of dg categories if
  there is a dg functor $G\colon \DS\rightarrow \CS$ such that $F\circ G\cong \Id_{\DS}$
  and $G\circ F\cong \Id_\CS$. 
  
  A dg functor $F\colon \CS\rightarrow \DS$ is \emph{quasi-fully faithful}, if
  restricts to quasi-isomorphisms on hom complexes, i.e. for every pair of
  objects $X$, $Y$ in $\CS$, the induced map $H^{\supb}(F)\colon
  H^{\supb}(\Hom_{\CS}(X,Y))\to H^{\supb}(\Hom_{\DS}(F(X),F(Y)))$ is an
  isomorphism; it is \emph{quasi-essentially surjective} if the functor $H^0(F)$
  between the respective homotopy categories is essentially surjective. If $F$
  is quasi-fully faithful and quasi-essentially surjective, then it is a called
  a \emph{quasi-equivalence}.

For many applications in this paper we will need a weaker notion of a functor
between dg categories. If $\mathrm{dgcat}$ is the category where the objects are
all (small) dg categories and the morphisms are all dg functors, one can define
\cite{MR2112034,MR2276263} the category $\Hqe$ as a localization of $\mathrm{dgcat}$ with respect to
quasi-equivalences. More abstractly, $\mathrm{dgcat}$ has a model category
structure whose weak equivalences are the quasi-equivalences, and $\Hqe$ is the
corresponding localization. 

Two dg categories $\CS$ and $\DS$ are called \emph{quasi-equivalent} if 
there exist dg categories $\BS_1,\ldots,\BS_n$ and a chain of quasi-equivalences
\[
\CS\leftarrow \BS_1\rightarrow \ldots \leftarrow \BS_n\rightarrow \DS.
\]
Then $\CS$ and $\DS$ are quasi-equivalent if and only if 
they are isomorphic in $\Hqe$.

A \emph{quasi-functor} between two dg categories is a morphism in $\Hqe$. For
example, if we have a dg functor $F\colon \CS\to \BS$ and a quasi-equivalence
$G\colon \DS\to \BS$ then $F$ induces a quasi-functor from $\CS$ to
$\DS$. Since quasi-equivalences induces equivalences of homotopy
categories, a quasi-functor between dg categories $\CS$ and $\DS$ induces an
honest functor between the homotopy categories $H^0(\CS)\to H^0(\DS)$.

\begin{remark}
Instead of working with $\Hqe$ and quasi-functors, one could choose to work with
$\Ai$-functors between dg categories. Over a field $\k$, this is essentially an
equivalent viewpoint since every quasi-equivalence admits an inverse, which is
in general not a dg functor but an $\Ai$-functor \cite{MR2441780}. However, we decided to stay
away from $\Ai$-functors in this paper. 
\end{remark}

\subsection{Bimodules}
\label{ss:bimodules}
If $\CS$ and $\DS$ are dg categories, a $\DS,\CS$-bimodule $\MS$ is the data of
\begin{itemize}
\item for each pair of objects $Y\in \DS$, $X\in \CS$, a dg $\k$-module $Y\MS X$.
\item for each quadruple of objects $Y,Y'\in\DS$ and $X,X'\in \CS$, action maps
\[
\Hom_{\DS}(Y,Y') \otimes_\k Y\MS X\otimes_\k \Hom_{\CS}(X',X) \rightarrow Y'\MS X'
\]
satisfying the usual associativity constraints.
\end{itemize}
The action maps are required to be chain maps of degree zero.  This is
equivalent to $|f\cdot m \cdot g| = |f|+|m|+|g|$ and the Leibniz rule
\[
d_{\MS}(f\cdot m \cdot g) = 
d_\DS(f)\cdot m \cdot g +(-1)^{|f|} f\cdot d_{\MS}(m)\cdot g + (-1)^{|f|+|m|} f\cdot m\cdot d_{\CS}(g)
\]
for all $f\in \Hom_{\DS}(Y,Y')$, $m\in Y\MS X$, $g\in \Hom_{\CS}(X',X)$.

The notation ${}_\DS \MS_\CS$ will be used to indicate that $\MS$ is a
$\DS,\CS$-bimodule.

A left $\CS$-module is the same as a $\CS,\k$-bimodule, and a right $\CS$-module
is the same as a $\k,\CS$-bimodule, by definition.

\begin{example}
If $\CS$ is a dg category and $X,Y\in \CS$ are objects, then we denote
\[
Y\CS X := \Hom_\CS(X,Y),\qquad\qquad \CS X :=\bigoplus_Y Y\CS X,\qquad\qquad Y \CS := \bigoplus_X Y\CS X.
\]
The composition of morphisms in $\CS$ equips $\CS X$ with the structure of a
left $\CS$-module, $Y\CS$ with the structure of a right $\CS$-module (called the
\emph{Yoneda modules}), and $\CS = \bigoplus_{X,Y} Y\CS X$ with the structure of a
$\CS,\CS$-bimodule (the \emph{regular bimodule}). 
\end{example}

We will use the Yoneda modules and the regular bimodule to streamline notation in certain places.

\begin{remark}
For dg categories $\CS,\DS$, the following notions are equivalent:
\begin{enumerate}
\item $\DS,\CS$-bimodules,
\item left $\DS\otimes_\k \CS^{\op}$-modules,
\item right $\CS\otimes_\k \DS^{\op}$-modules,
\item functors $\DS\otimes_\k \CS^{\op}\rightarrow \k\text{-dgmod}$.
\end{enumerate}
However, such identifications necessarily involve choices and hidden signs; for
this reason, we will typically not use them.
\end{remark}

Given dg categories $\BS,\CS,\DS$ and bimodules ${}_\DS \MS_{\CS}$, ${}_\CS \NS_\BS$, 
their tensor product ${}_{\DS}(\MS\otimes_\CS \NS)_{\BS}$ is the bimodule with
\[
Z(\MS\otimes_\CS \NS)X :=\bigoplus_{Y\in \CS} \Big(Z\MS Y\otimes Y\NS X\Big)\Big/\sim,
\]
where $\sim$ is the equivalence relation $(m\cdot f)\otimes n \sim m\otimes (f\cdot n)$ 
for all $m\in Z\MS Y$, $f\in Y\CS Y'$, $n\in Y'\NS X$.

\begin{remark}
The category of $\CS,\CS$-bimodules is monoidal with tensor product as defined
above, and monoidal identity given by the regular bimodule $\CS$.
\end{remark}

\begin{remark}
We often regard $\CS$ and $\DS$ as the (very big) non-unital algebras
\[
\CS = \bigoplus_{X',X\in \CS} X'\CS X,\qquad \qquad \DS = \bigoplus_{Y',Y\in \DS} Y'\DS Y,
\]
and a $\DS,\CS$-bimodule $\MS$ as the (very big) dg bimodule
\[
\MS = \bigoplus_{Y\in \DS,X\in \CS} Y\MS X.
\]
In this language an object $X\in \CS$ corresponds to the distinguished idempotent $\Id_X$
in the big algebra $\bigoplus_{X',X} X'\CS X$. In this way, essentially all of one's
intuition from the usual world of algebras and bimodules carries over into the
world of dg categories and their bimodules.
\end{remark}

\subsection{\texorpdfstring{$\Ai$}{A-infinity} algebras and deformation retracts}
\label{subsec:Ainf}

Recall, an $\Ai$-algebra is a graded $\k$-module $A$ equipped with maps
$\mu_n:A^{\otimes n}\rightarrow A$ degree $2-n$, $n\geq 1$, satisfying the
following family of identities for $M\geq 1$
\[
\sum_{M=r+s+t} (-1)^{r+s t} \mu_{r+1+t}(\Id^{\otimes r}\otimes \mu_s \otimes \Id^{\otimes t}) =0.
\]
In particular a dg algebra is an $\Ai$-algebra in which $\mu_n$ vanish for
$n\neq 1,2$. In this case $\mu_1:A\rightarrow A$ is the differential and
$\mu_2:A\otimes A\rightarrow A$ is honestly associative and satisfies the
Leibniz rule
\[
\mu_1\circ \mu_2 - \mu_2\circ (\mu_1\otimes \Id_A) - \mu_2\circ (\Id_A\otimes \mu_1)=0
\]

Suppose $X,Y\in \CS$ are objects in a dg category.  A deformation retract
$X\rightarrow Y$ is the data of closed degree zero morphisms $\pi:X\rightarrow
Y$, $\sigma:Y\rightarrow X$ and a degree $-1$ homotopy $h\in \End^{-1}(X)$ such
that
\[
\pi\circ \sigma = \Id_Y,\qquad d(h) = \Id_X - \sigma\circ \pi,\qquad h\circ \sigma = 0 = \pi\circ h.
\]

The following is well known \cite{KellerAinfty,Kadeishvili,Merkulov}.
\begin{theorem}
    \label{thm:Ainfty-defretracts}
If $A$ is a dg algebra and $V$ is a dg $\k$-module, then any deformation
 retract $A\rightarrow V$ gives $V$ the structure of an $\Ai$-algebra quasi-isomorphic to $A$.
\end{theorem}

\begin{proof}
 The construction follows \cite{Merkulov}. Let $m:A\otimes A\to A$ be the
 multiplication in $A$. We define recursively a sequence of maps
 $\lambda_k:V^{\otimes k}\to A$ by $\lambda_2=m(\sigma\otimes \sigma)$ and
 $$
 \lambda_n=-m(\sigma\otimes h\lambda_{n-1})+\sum_{s=2}^{n-2}(-1)^{s+1}m(h\lambda_{s}\otimes  
 h\lambda_{n-s})+(-1)^{n+1}m(h\lambda_{n-1}\otimes \sigma),\ n\ge 3
 $$
 Then $\mu_n=\pi\circ \lambda_n$ defines the structure of a strictly unital $\Ai$-algebra on $V$.
\end{proof}
\begin{remark}
Suppose that $A$ and $V$ have an additional grading which is preserved by the
differential and the maps $\sigma,\pi$ and $h$. Then the $\Ai$-structure maps on
$V$ can be chosen to preserve this grading as well.
\end{remark}
\begin{remark}
  \label{rem:Ainfty-structure-linear}
 If $R\subset A$ is a commutative dg subalgebra then the multiplication on $A$ is $R$-bilinear, that is, 
descends to the map of $R$-bimodules
\[
m:{}_{R}{A}\otimes_R A_{R}\to {}_{R}{A}_{R}.
\]
Similarly, if $V$ admits the structure of
 a dg $R$-bimodule, and the data of the deformation retract (that is $\sigma,\pi$
 and $h$) can be chosen to be $R$-bilinear, then the $\Ai$-structure maps on $V$
 can be chosen to be $R$-linear, in the sense that they descend to the quotient
 \[
 \mu_n:{}_{R}{V}\otimes_R V\otimes_R \cdots \otimes_R V_{R}\rightarrow {}_{R}{V}_{R}.
 \]
 This elementary fact is often very useful.
 \end{remark}

If $\k$ is a field, then $A$ deformation retracts onto its homology $H(A)$
(regarded as a dg $\k$-module with zero differential), and so $H(A)$ inherits
the structure of an $\Ai$-algebra with $\mu_1=0$.
A differential graded algebra $A$ is called {\em formal} if it is
 quasi-isomorphic to its homology $H(A)$. The above discussion shows that a dg
 algebra $A$ over a field $\k$ 
 is formal if and only if the  $\Ai$
 structure on $H(A)$ is trivial, that is, $\mu_k=0$ for $k>2$.

\section{Standard dg categorical constructions}
\label{sec: pretriangulated dg}
In this and the following section we describe the processes of adjoining finite
direct sums, suspension, and twists to a dg
category.

A functor $F\colon \CS\rightarrow \k\dgmod$ is said to be \emph{representable}
if there is an object $X\in \CS$ such that $F$ is isomorphic to $\Hom_{\CS}(X,-)$ (or
$\Hom_{\CS}(-,X)$ if $F$ is contravariant). 

The dg category $\CS$ is \emph{additive} if for each finite collection of
objects $X_i\in \CS$ the functor $Y\mapsto \bigoplus_i Y\CS X_i$ is
representable. This means $\CS$ has finite coproducts, for which we use the
symbol $\oplus$. 

The dg category $\CS$ is \emph{suspended} (or \emph{has suspension}) if for each
$X\in \CS$ and $l \in \Z$ the functor $Y\mapsto \susp{l}(Y\CS X)$ is representable.
This is equivalent to the existence of an object $\susp{l}X$, an
\emph{$l$-translate}, for every object $X$ of $\CS$, together with a given
closed degree $l$ invertible morphism $\susp{l}X\to X$.

Let $\a\in \End^1_\CS(X)$ be an endomorphism in $\CS$ satisfying the
\emph{Maurer--Cartan equation} $d_\CS(\a)+ \a\circ \a = 0$. Then we have a
functor $\phi_{\a}\colon\CS \to \k\dgmod$ sending an object $Y$ to the complex
$(Y\CS X, d_\a)$ with \emph{twisted differential} $d_\a(f) := d(f) - (-1)^{|f|}
f\circ \a$.

The category dg $\CS$ is said to \emph{have twists} if for each Maurer--Cartan
element $\a$ in $\CS$ the functors $\phi_{\a}$ is representable.

\subsection{Additive suspended envelope}

If $\CS$ is a dg category, the \emph{additive suspended envelope} $\adds\CS$ of
$\CS$ is the dg category whose objects are collections $\{\susp{a_i}X^i\}_{i\in I}$
where $I\subset \Z$ is a finite set, $X^i\in \CS$ and $a_i\in \Z$. Morphism
complexes in $\adds\CS$ are by definition
\[
\Hom^l_{\adds\CS}\left(\{\susp{a_i}X^i\}_{i\in I},\{\susp{b_j}Y^j\}_{j\in J}\right) 
= \prod_{i\in I}\bigoplus_{j\in J} \Hom^{l+b_j-a_i}_{\CS}(X^i,Y^j)
\]
with differential 
\begin{equation}\
  \label{eq:diff-susp-env}
d_{\adds\CS}((f_{ji})_{(j,i)\in J\times I}) = ((-1)^{b_j}d_{\CS}(f_{ji}))_{(j,i)\in J\times I}.
\end{equation}

An element of this hom space can be thought of as a $J\times I$ matrix
$(f_{ji})$ of morphisms $f_{ji}\in \Hom_\CS(X^i,Y^J)$.  Composition of morphisms
is given by usual matrix multiplication and composition in $\CS$.

There is a canonical fully faithful dg functor $\CS\rightarrow \adds\CS$ defined
object-wise by $X\mapsto \{X\}$ (with indexing set $I$ a singleton), and we may
identify $\CS$ with its image in $\adds\CS$. It is straightforward to verify
that $\adds\CS$ is additive and suspended (with $\susp{l}$ indicating
$l$-translates), and we will henceforth abuse notation by writing
\[
\bigoplus_{i\in I} \susp{a_i}X^i := \{\susp{a_i}X^i\}_{i\in I} \in \adds\CS.
\]
We also write $0$ for the empty direct sum, corresponding to the case
$I=\emptyset$. 

It is not hard to check that $\adds\adds\CS \cong \adds \CS$ and this idempotent
property of the assignment $\CS \mapsto \adds \CS$ justifies the name additive
suspended envelope.

\begin{remark} We also define an additive suspended envelope $\adddp\CS$ with
countable direct products, where the finiteness assumption on the indexing sets
$I\subset \Z$ is removed.  Homogeneous morphisms in $\adddp\CS$ are by definition
matrices of morphisms in $\CS$, each row of which has only finitely many
nonzero entries.
\end{remark}

\subsection{Twisted envelope and pretriangulated hull}

The \emph{twisted envelope} of $\CS$ can be constructed explicitly as follows.
Let $\Tw(\CS)$ be the category with objects $\tw_\a(X)$ where $X\in\CS$ and $\a\in
\End^1_\CS(X)$ satisfying $d_{\CS}(\a)+\a^2=0$ as above. The morphism complexes
in $\Tw(\CS)$ are by definition
\[
\Hom_{\Tw(\CS)}\Big(\tw_\a(X),\tw_\b(Y)\Big) := \Hom_\CS(X,Y)
\]
with differential
\[
d_{\Tw(\CS)}(f) = d_{\CS}(f) + \b\circ f - (-1)^{|f|} f\circ \a.
\]

We say that $\CS$ \emph{has twists} if the obvious fully faithful dg functor
$\CS \rightarrow \Tw(\CS)$ sending $X\mapsto \tw_0(X)$ is an equivalence. The
natural inclusion $\Tw(\CS)\rightarrow \Tw(\Tw(\CS))$ sending $\tw_\a(X)\mapsto
\tw_0(\tw_\a(X))$ is an equivalence, with an inverse equivalence
$\Tw(\Tw(\CS))\rightarrow \Tw(\CS)$ defined by $\tw_\b(\tw_\a(X))\mapsto
\tw_{\a+\b}(X)$ (compare with the procedure of taking the total complex of a
bicomplex).  Thus, $\Tw(\CS)$ has twists. This idempotent property of the
assignment $\CS\mapsto \Tw(\CS)$ justifies our referring to $\Tw(\CS)$ as the
twisted envelope.

\begin{example}
Any additive category $\AS$ can be thought of as a dg category with zero
differential and trivial grading (all morphisms are placed in degree zero). Then
$\Tw(\Sigma\AS)$ is equivalent to the usual category of bounded complexes
$\Ch^b(\AS)$.
\end{example}

If $f\in \Hom^1_\CS(X,Y)$ is a degree $1$ closed morphism, then the cone of $f$
is the object $\tw_\a(X\oplus Y)$ with $\a=\smMatrix{0& 0\\ f & 0}$ inside
$\Tw(\CS)$.  If instead $f\colon X\rightarrow Y$ is a degree zero closed
morphism and $\CS$ has suspension, then we first replace $f$ by a degree 1
closed morphism $\susp{1}X\rightarrow Y$ and apply the previous construction. 

\begin{definition}
We say a dg category $\CS$ is \emph{pretriangulated} if it is suspended and
closed under taking cones.

The \emph{pretriangulated hull} $\pretr{\CS}$ of a dg category $\CS$ is the full subcategory of
$\Tw(\adds\CS)$ generated by $\adds\CS$ under taking mapping cones.
\end{definition}

It follows from the discussion above that $\pretr{\CS}$ is pretriangulated, and
$\CS$ itself is pretriangulated if and only if the natural embedding $\CS\to
\pretr{\CS}$ is an equivalence.

Objects in the pretriangulated hull can be expressed as iterated mapping cones
of objects in $\adds\CS$, also known as \emph{one-sided twisted complexes}. 

\begin{example} For $X=\tw_\a(\bigoplus_{j\in J}\susp{a_j}X^j)$ we can collect terms
with equal shifts, i.e. set $Y^i:=\bigoplus_{j\colon a_j=-i} X^j$ for $i\in \Z$, and write
$\a_{j,i}\in \Hom^{1+i-j}(Y^i, Y^j)$ for the components of the twist. Then $X$ can be
illustrated as:
  \[
    \begin{tikzpicture}[baseline=0em]
    \tikzstyle{every node}=[font=\scriptsize]
    \node (a) at (0,0) {$\cdots$};
    \node (b) at (-2.5,0) {$\Sigma^{-1-i}Y^{i+1}$};
    \node (c) at (-5,0) {$\Sigma^{-i}Y^i$};
    \node (d) at (-7.5,0) {$\Sigma^{1-i}Y^{i-1}$};
    \node (e) at (-10,0) {$\cdots$};
    \draw[frontline,->,>=stealth,shorten >=1pt,auto,node distance=1.8cm,thick]
    (c) to[bend left=50] node {$\alpha_{i+2,i}$} (a);
    \draw[frontline,->,>=stealth,shorten >=1pt,auto,node distance=1.8cm,thick]
    (d) to[bend left=50] node {$\alpha_{i+1,i-1}$} (b);
    \draw[frontline,->,>=stealth,shorten >=1pt,auto,node distance=1.8cm,thick]
    (e) to[bend left=50] node {$\alpha_{i,i-2}$} (c);
    \draw[frontline,->,>=stealth,shorten >=1pt,auto,node distance=1.8cm,thick]
    (d) to[bend left=70] node {$\alpha_{i+2,i-1}$} (a);
    \draw[frontline,->,>=stealth,shorten >=1pt,auto,node distance=1.8cm,thick]
    (e) to[bend left=70] node {$\alpha_{i+1,i-2}$} (b);
    \path[->,>=stealth,shorten >=1pt,auto,node distance=1.8cm,
      thick]
    (b) edge node 	{$\alpha_{i+2,i+1}$} (a)
    (c) edge node  {$\alpha_{i+1,i}$} (b)
    (d) edge node {$\alpha_{i,i-1}$} (c)
    (e) edge node {$\alpha_{i-1,i-2}$} (d);
    \end{tikzpicture}
    \]
The Maurer--Cartan equation for $\a$ now is $(-1)^j d_\CS(\a_{j,i})  +
\sum_{j>k>i}\a_{j,k}\circ \a_{k,i} = 0$.
\end{example}

\begin{remark} Let $A$ be a dg algebra supported in non-positive degrees,
considered as a dg category with one object.  Then $\pretr{A}=\Tw(\adds A)$
since every twisted complex is one-sided, i.e. an iterated cone.
\end{remark}

\subsection{$\Ai$-categories}

In this section we briefly discuss $\Ai$-categories and refer the reader to
\cite{Keller,Lefevre,MR2441780} for more complete exposition.

An $\Ai$-category $\CS$ consists of a set of objects $\mathrm{Ob}(\CS)$, a
graded vector space $\Hom(X,Y)$ for each pair of objects $X$, $Y$ and degree $2-d$
composition maps
$$
\mu_d:\Hom(X_1,X_0)\otimes   \cdots \otimes \Hom(X_{d},X_{d-1})\to \Hom(X_d,X_0)
$$
satisfying $\Ai$-equations as in Section \ref{subsec:Ainf}.

An $\Ai$-category $\CS$ with one object $X$ is the same data as the $\Ai$-algebra
$\End_{\CS}(X)$. If $\mu_d=0$ for all $d\ge 3$ in some $\Ai$-category $\CS$,
then $\CS$ is just a dg category with differential $\mu_1:\Hom(X_1,X_0)\to
\Hom(X_1,X_0)$ and composition $\mu_2:\Hom(X_1,X_0)\otimes \Hom(X_2,X_1)\to
\Hom(X_0,X_2)$ which in this case is strictly associative.

Given an $\Ai$-category $\CS$, we can consider its homotopy category where the
objects are the same as in $\CS$, and the morphisms are given by the homology
with respect to $\mu_1$. The composition of morphisms is induced by $\mu_2$.

Given objects $X,X'$ of an $\Ai$ category $\CS$ and $k\in \Z$, we say that $X'$
is the $k$-fold suspension of $X$, written $X'\cong X[k]$ if there is a morphism
$\phi:X\rightarrow X'$ of degree $k$ which is closed ($\mu_1(\phi)=0$) and
invertible (there exists $\phi':X'\rightarrow X$ of degree $-k$ such that
$\mu_2(\phi,\phi')=\Id_{X'}$ and $\mu_2(\phi',\phi)=\Id_X$).

Given an $\Ai$ category $\CS$ we let $\Sigma \CS$ denote the closure of $\CS$
with respect to finite direct sums and suspensions.  Replacing $\CS$ by
$\Sigma\CS$ if necessary, below we will assume that $\CS$ is closed under finite
direct sum and suspension.

Given an object $X\in \CS$, a Maurer-Cartan endomorphism of $X$ is a degree 1
element $\a\in \End_\CS(X)$ such that
\[
\sum_{d\geq 1}\mu_d(\a,\ldots,\a) =0.
\]

We say that $\a$ is \emph{one-sided} if there exists a direct sum decomposition
$X\cong \bigoplus_{i\in I} X_i$ where $I$ is a finite poset, with respect to
which $\a$ is represented by a strictly lower triangular matrix.

The category $\pretr{\CS}$ is the category with
\begin{itemize}
\item objects of $\pretr{\CS}$ are formal expressions $\tw_\a(X)$ in which $\a$
is a one-sided twist, 
\item $\Hom_{\pretr{\CS}}(\tw_\a(X),\tw_\b(Y)) =\Hom_{\CS}(X,Y)$,
\item the higher composition $\mu_d'$ of a sequence of morphisms
\[
\tw_{\a_0}(X_0) \buildrel f_1\over \leftarrow\tw_{\a_1}(X_1) \buildrel f_2\over \leftarrow \cdots \buildrel f_d\over \leftarrow \tw_{\a_d}(X_d),
\]
is defined by
\[
\mu'_d(f_1,\ldots,f_d) \ := \ \sum_{r_0,\ldots,r_d\geq 0} \pm \mu_{d+r_0+\cdots+r_d}\left(\underbrace{\a_0,\ldots,\a_0}_{r_0},f_1,\underbrace{\a_1,\ldots,\a_1}_{r_1},\ldots, f_d, \underbrace{\a_d,\ldots,\a_d}_{r_d}\right)
\]
(note that the one-sidedness of the Maurer-Cartan elements $\a_i$ guarantees finiteness of the above sum).
\end{itemize}

\section{Homotopy idempotents and the Karoubi envelope}
\label{sec:Karoubi}

\subsection{Homotopy idempotents}
\label{ss:homotopy idempts}

Let $\CS$ be a dg category.  A \emph{homotopy idempotent} in $\CS$ is a closed
endomorphism $e\in \End_\CS^0(X)$ such that $e^2\simeq e$ (i.e.~ an idempotent
in $H^0(\CS)$.  We say that $Y$ is an \emph{image of $e$} if there exist closed
degree zero morphisms $\sigma\colon Y\rightarrow X$, $\pi\colon X\rightarrow Y$
such that
\[
\pi\circ \sigma\simeq \Id_Y,\qquad\qquad \sigma\circ \pi\simeq e.
\]

\begin{lemma}\label{lemma:properties of images}
Images of homotopy idempotents satisfy the following basic properties:
\begin{enumerate}
\item Suppose $e_1,e_2\in \End_\CS(X)$ are homotopy idempotents with $e_1\simeq
e_2$.  If $Y_i$ is an image of $e_i$ ($i=1,2$) then $Y_1\simeq Y_2$.  In
particular the image of a homotopy idempotent is unique up to homotopy
equivalence.
\item If $Y$ is an image of a homotopy idempotent $e\in \End_\CS(X)$ then any
homotopy idempotent $e'\in \End_\CS(Y)$ determines a homotopy idempotent $e''\in
\End_\CS(X)$ with the property that $Z$ is an image of $e'$ if and only if it is
an image of $e''$.
\item If $Y_0,Y_1$ are the images of homotopy idempotents $e_i\in \End_\CS(X_i)$
($i=0,1$) then
\[
\Hom_{H^0(\CS)}(Y_0,Y_1) \cong e_1 \Hom_{H^0(\CS)}(X_0,X_1) e_0.
\]
\end{enumerate}
\end{lemma}
\begin{proof}
Exercise.
\end{proof}

\begin{definition}\label{def:idemp compl}
We say that $\CS$ is \emph{homotopy idempotent complete} if $H^0(\CS)$ is
idempotent complete, i.e.~each homotopy idempotent in $\CS$ has an image.
\end{definition}

Our goal in this section is to construct the \emph{homotopy Karoubi envelope}
$\Kar^{dg}(\CS)$ and prove the following.

\begin{theorem}
  \label{thm:Kar}
Every dg category $\CS$ admits an embedding $\CS\hookrightarrow \Kar^{dg}(\CS)$
into a homotopy idempotent complete dg category characterized up to
quasi-equivalence by the following universal property: if $\DS$ is a homotopy
idempotent complete dg category equipped with a dg functor $\CS\rightarrow \DS$,
then there is unique morphism (quasi-functor) $\Kar^{dg}(\CS)\rightarrow \DS$ in
$\Hqe$ such that the following diagram commutes:
\[
\begin{tikzcd}
  \CS\arrow[r] \arrow[d,hook] & \DS\\
  \Kar^{dg}(\CS) \arrow[ur, dashed]   & 
  \end{tikzcd}
\]
Furthermore, $\CS$ is idempotent complete if and only if the canonical functor
$\CS\rightarrow \Kar^{dg}(\CS)$ is a quasi-equivalence.
\end{theorem}

\subsection{$\Ai$-idempotents}
\label{ss:a infty idempts}

Note that if $e_0\in \End^0_\CS(X)$ is a homotopy idempotent and $h_0\in
\End^{-1}_\CS(X)$ satisfies $d_\CS(h_0) = e_0\circ (\Id_X-e_0)$, then $e_0\circ
h_0 - h_0\circ e_0$ is automatically closed.  This morphism obstructs certain
constructions, and it is natural to require it to be null-homotopic (as we will
see below, one can choose a homotopy $h_0$ such that this holds, but not every
homotopy satisfies this condition), via an endomorphism we will denote $e_1$.
There is a higher family of obstructions which is natural to require to be
trivial, via homotopies $e_k$, $h_k$ for $k\in \Z_{\geq 1}$.  This results in the
notion of an $\Ai$-idempotent or idempotent up to coherent homotopy, which we
describe next.
\begin{definition}\label{def:Ainfty idempt}
An \emph{$\Ai$-idempotent} in $\CS$ is a triple $(X,\ue,\uh)$, consisting of an
object $X$, and a collection of endomorphisms $\ue = \{e_k\in
\End^{-2k}_\CS(X)\}_{k=0}^\infty$, $\uh = \{h_k\in
\End^{-1-2k}_\CS(X)\}_{k=0}^\infty$ satisfying
\begin{subequations}
\begin{eqnarray}
d(e_k) &=& \sum_{i+j=k-1}(e_ih_j - h_i e_j)\label{eq:d of e 2}\\
d(h_k) &=& e_k-\sum_{i+j=k} e_i e_j  - \sum_{i+j=k-1} h_i h_j\label{eq:d of h 2}
\end{eqnarray}
\end{subequations}
The \emph{complement} of the $\Ai$-idempotent $(X,\ue,\uh)$ is the $\Ai$
idempotent $(X,\ue,\uh)^\perp := (X, \ue^\perp, \uh^\perp)$ where
\[
e_0^\perp = \Id_X- e_0,\qquad e^\perp_k=-e_k\quad (k\geq 1), \qquad h^\perp_k=h_k \quad(k\geq 0).
\]
Verification that this defines an $\Ai$-idempotent is left to the reader.
\end{definition}

\begin{definition}\label{def:image of Ai idemp}
If $(X,\ue,\uh )$ is an $\Ai$-idempotent, then let $Z(X,\ue,\uh)\in
\Tw(\adddp\CS)$ denote the twisted complex of the form $\tw_\d(\prod_{k\geq
0}\Sigma^{-k}X)$ and differential $\d$ given in terms of components by 
\[
\d_{ji} =\begin{cases}
-e_k^\perp & \text{for $i$ even, $j = i+1+2k$}\\
h_k & \text{ for $i$ even, $j=i+2k+2$}\\
e_k & \text{ for $i$ odd, $j = i+1+2k$}\\
-h_k & \text{ for $i$ odd, $j=i+2k+2$}\\
 \end{cases}
\]
(recall that $-e_k^\perp = e_k$ for $k\geq 1$).
\end{definition}
This twisted complex can be visualized as
\[
\begin{tikzpicture}[baseline=0em]
\tikzstyle{every node}=[font=\scriptsize]
\node (a) at (0,0) {$\underline{X}$};
\node (b) at (2.5,0) {$\Sigma\inv X$};
\node (c) at (5,0) {$\Sigma^{-2} X$};
\node (d) at (7.5,0) {$\Sigma^{-3} X$};
\node (e) at (10,0) {$\cdots$};
\draw[frontline,->,>=stealth,shorten >=1pt,auto,node distance=1.8cm,thick]
(a) to[bend left=50] node {$h_0$} (c);
\draw[frontline,->,>=stealth,shorten >=1pt,auto,node distance=1.8cm,thick]
(b) to[bend left=50] node {$-h_0$} (d);
\draw[frontline,->,>=stealth,shorten >=1pt,auto,node distance=1.8cm,thick]
(c) to[bend left=50] node {$h_0$} (e);
\draw[frontline,->,>=stealth,shorten >=1pt,auto,node distance=1.8cm,thick]
(a) to[bend left=70] node {$e_1$} (d);
\draw[frontline,->,>=stealth,shorten >=1pt,auto,node distance=1.8cm,thick]
(b) to[bend left=70] node {$e_1$} (e);
\path[->,>=stealth,shorten >=1pt,auto,node distance=1.8cm,
  thick]
(a) edge node 	{$e_0-\Id_X$} (b)
(b) edge node  {$e_0$} (c)
(c) edge node {$e_0-\Id_X$} (d)
(d) edge node {$e_0$} (e);
\end{tikzpicture}.
\]
with length $>3$ components of the differential not pictured.  

\begin{example}\label{ex:Z(X,1,0)}
The identity of $X$ gives an $\Ai$-idempotent with $e_0=\Id_X$ and
$e_{k+1}=0=h_k$ for $k\geq 0$.  The resulting twisted complex
\[
Z(X,\Id_X,0) \ \ = \ \ \begin{tikzpicture}[baseline=0em]
\tikzstyle{every node}=[font=\small]
\node (a) at (0,0) {$\underline{X}$};
\node (b) at (2.5,0) {$\Sigma\inv X$};
\node (c) at (5,0) {$\Sigma^{-2} X$};
\node (d) at (7.5,0) {$\Sigma^{-3} X$};
\node (e) at (10,0) {$\cdots$};
\path[->,>=stealth,shorten >=1pt,auto,node distance=1.8cm,
  thick]
(a) edge node 	{$0$} (b)
(b) edge node  {$\Id_X$} (c)
(c) edge node {$0$} (d)
(d) edge node {$\Id_X$} (e);
\end{tikzpicture}
\]
is homotopy equivalent to $X$ after the cancellation of contractible summands.
\end{example}

We wish to prove that $Z(X,\ue,\uh)$ is a well-defined twisted complex, and is
an image of of $e_0$.  To prove this requires a considerable amount of
bookkeeping, for which it is useful to consider the generating functions $e(\ep)
= \sum_{k\geq 0} e_k \ep^k$ and $h(z) = \sum_{k\geq 0} h_k \ep^k$ where $\ep$ is
a formal indeterminate of degree 2.  Below we work an abstract dg category which
formalizes the relations satisfied by these generating functions.

\subsection{Abstract $\Ai$-idempotents and their images}
\begin{definition}\label{def:representing ai idemp}
Let $\ep$ be a formal indeterminate of cohomological degree 2.  Let $\Rid$ be
the $\k\llbracket\ep\rrbracket$-linear dg category with one object $\XB$ whose
endomorphism complex is freely generated by endomorphisms $e,h\in
\End_{\Rid}(\XB)$ satisfying
\[
\deg(e) = 0,\qquad\qquad \deg(h)=-1,
\]
\begin{subequations}
\begin{eqnarray}
d(e) &=& \ep(eh - he)\label{eq:d of e}\\ 
d(h) &=& e-e^2 - \ep h^2\label{eq:d of h},
\end{eqnarray}
\end{subequations}
extended to arbitrary morphisms by the Leibniz rule.
\end{definition}
To check that $\Rid$ is a dg category, one must check that $d^2=0$ on all
morphisms.  It suffices to check on the generating morphisms, which is
straightforward.  For instance to verify that $d^2(e)=0$ it suffices to verify
that $d(eh) = d(he)$, which follows from the computations
\[
d(eh) = d(e)h + e d(h) = \ep(eh-he)h + e(e-e^2 - \ep h^2) = e^2-e^3 - \ep heh
\]
and
\[
d(he) = d(h)e - h d(e) = (e-e^2-\ep h^2)e - \ep h (eh-he) =e^2-e^3  - \ep heh.
\]
The proof that $d^2(h)=0$ is equally straightforward. 

The notion of an $\Ai$-idempotent in $\CS$ can now be described as follows.
Consider the dg category $\CS\llbracket \ep\rrbracket$ with the same objects as
$\CS$, and morphism complexes
\[
\Hom_{\CS\llbracket \ep\rrbracket}(X,Y) \ := \ \Hom_\CS(X,Y)\otimes_\k \k\llbracket \ep\rrbracket.
\]
Then an $\Ai$-idempotent in $\CS$ is equivalent to a dg functor $\Rid\rightarrow
\CS\llbracket \ep\rrbracket$.  The image of $\XB$ in $\CS$ is an object $X\in
\CS$, and the images of $e,h$ are formal series of morphisms $h(\ep)$,
\[
e(\ep) = \sum_{k\geq 0} e_k \ep^k\qquad \quad (e_k \in \End^{-2k}_\CS(X)),
\]
\[
h(\ep) = \sum_{k\geq 0} h_k \ep^k \qquad \quad (h_k \in \End^{-1-2k}_\CS(X)),
\]
satisfying the identities \eqref{eq:d of e}, \eqref{eq:d of h}.  In terms of
components, this yields \eqref{eq:d of e 2},\eqref{eq:d of h 2}.  

The following is responsible for the notion of complementary $\Ai$-idempotents.

\begin{lemma}\label{lemma:symmetry of Rid}
There is an automorphism of $\Rid$ which sends $e\mapsto 1-e$ and fixes $h$.\qed
\end{lemma}

Now we have an analogue of Definition \ref{def:image of Ai idemp}.  

\begin{definition}\label{def:abstract image}
Let $Z(e,h)$ denote the twisted complex $\tw_\a(\XB\oplus \Sigma\inv \XB)$ where
\[
\a = \sqmatrix{\ep h &  \ep e  \\ e-1 & -\ep h}.
\]
\end{definition}

We also have an analogue of $X$, viewed through Example \ref{ex:Z(X,1,0)}.

\begin{definition}\label{def:abstract X}
Let $Z(1,0)\in \Tw(\Sigma \Rid)$ denote the twisted complex $\tw_\b(\XB\oplus
\Sigma\inv \XB)$ where
\[
\b = \sqmatrix{0&\ep\\ 0 & 0}.
\]
\end{definition}

Our goal is show that $Z(e,h)$ is the image of a homotopy idempotent acting on
$Z(1,0)$.

\begin{lemma}\label{lemma:abstract idempotent}
The object $Z(1,0)=\tw_\b(\XB\oplus \Sigma\inv \XB)$ is a well-defined twisted
complex in $\Tw(\adddp \Rid)$, and the following defines a homotopy idempotent
$E$ acting on $Z(1,0)$:
\[
 E\  :=\  \sqmatrix{e&0\\ he-eh & e}.
\]
\end{lemma}
\begin{proof}
The first statement is clear since $d(\b)=0=\b^2$.  For the second statement,
let
\[
H \ := \ \sqmatrix{h&0\\h^2&-h},
\]
regarded as an endomorphism of $\XB\oplus \Sigma\inv \XB$.  It is
straightforward to check that 
\[
d(H) + \b H + \b H = E - E^2.
\]
For this, one must keep in mind that a sign appears in the bottom row of the
matrix representing $d(H)$, due to the signs involved in differentiating
morphisms in the suspended envelope \eqref{eq:diff-susp-env}:
\[
d\left(\sqmatrix{h&0\\h^2&-h}\right) = \sqmatrix{d(h)&0\\ -d(h^2)&d(h)}.
\]
\end{proof}

\begin{lemma}\label{lemma:abstract image}
The object $Z(e,h)=\tw_\a(\XB\oplus \Sigma\inv \XB)$ from Definition
\ref{def:abstract image} is a well-defined twisted complex in $\Tw(\adddp
\Rid)$; moreover $Z(e,h)$ is an image of the homotopy idempotent $E$ acting on
$Z(1,0)$ (from Lemma \ref{lemma:abstract idempotent}).
\end{lemma}
\begin{proof}
To show that $Z(e,h)$ is a well-defined twisted complex we must check that
$\a\in \End_{\Sigma \Rid}(\XB\oplus \Sigma\inv \XB)$ satisfies the Maurer-Cartan
equation
\[
\sqmatrix{\ep d(h) &  \ep d(e)  \\ -d(e) & \ep d(h)}+\sqmatrix{\ep h &  \ep e  \\ e-1 & -\ep h}\sqmatrix{\ep h &  \ep e  \\ e-1 & -\ep h} = 0.
\]
This is easily verified. 

Now, define maps $\sigma:Z(e,h)\rightarrow Z(1,0)$ and $\pi:Z(1,0)\rightarrow
Z(e,h)$ by the matrices
\[
\sigma = \sqmatrix{1 & 0 \\ h & e},\qquad\qquad \pi = \sqmatrix{e & 0 \\ -h & 1}.
\]
Observe that $\sigma\circ \pi = E$.  To check $\sigma$ is closed is a
computation:
\[
\sqmatrix{0 & 0 \\ -d(h) & -d(e)} + \sqmatrix{0 & \ep \\ 0 & 0}\sqmatrix{1 & 0 \\ h & e} -  \sqmatrix{1 & 0 \\ h & e}\sqmatrix{\ep h &  \ep e \\ e-1 & -\ep h} = 0,
\]
and that $\pi$ is closed is the computation
\[
\sqmatrix{d(e) & 0 \\ d(h) & 0} + \sqmatrix{\ep h &  \ep e \\ e-1 & -\ep h} \sqmatrix{e & 0 \\ -h & 1} - \sqmatrix{e & 0 \\ -h & 1}\sqmatrix{0 & \ep \\ 0 & 0} =0,
\]
both of which are straightforward.

 It remains to show that $\pi\circ \sigma\simeq \Id_{Z(e,h)}$.  Observe that
\[
\pi\circ \sigma = \sqmatrix{e & 0 \\ -h & 1} \sqmatrix{1 & 0 \\ h & e} = \sqmatrix{e & 0 \\ 0 & e}.
\]
Now, let $K=\smMatrix{0&1\\ 0&0}\in \End_{\Sigma \Rid}(\XB\oplus \Sigma\inv
\XB)$.  The following computes $d(K)+\a K + K \a$:
\[
d\left(\sqmatrix{0&1\\ 0&0}\right) 
+ \sqmatrix{\ep h &  \ep e \\ e-1 & -\ep h}\sqmatrix{0&1\\ 0&0}
+\sqmatrix{0&1\\ 0&0}\sqmatrix{\ep h &  \ep e \\ e-1 & -\ep h} 
= \sqmatrix{e-1 &\ep h - \ep h\\ 0& e-1},
\]
which shows that $\pi\circ \sigma \simeq \Id_{Z(e,h)}$.  This completes the
proof that $Z(e,h)$ is an image of the idempotent $E$ acting on
$Z(1,0)$.\end{proof}

The twisted complexes $Z(1,0)$ and $Z(e,h)$ and the maps relating them can be
pictured diagrammatically as
\[
\begin{tikzpicture}[baseline=0em]
\tikzstyle{every node}=[font=\small]
\node (a) at (0,2.5) {$\XB$};
\node (b) at (3,2.5) {$\Sigma\inv \XB$};
\node (c) at (0,0) {$\XB$};
\node (d) at (3,0) {$\Sigma\inv \XB$};
\node (e) at (0,-2.5) {$\XB$};
\node (f) at (3,-2.5) {$\Sigma\inv \XB$};
\node (g) at (-4,2.5) {$Z(1,0)$};
\node (h) at (-4,0) {$Z(e,h)$};
\node (i) at (-4,-2.5) {$Z(1,0)$};
\node at (-2.5,2.5) {$=$};
\node at (-2.5,0) {$=$};
\node at (-2.5,-2.5) {$=$};
\path[->,>=stealth,shorten >=1pt,auto,node distance=1.8cm, thick]
(b) edge node {$\ep$} (a)
(c) edge [loop left, looseness = 7,in=218,out=142] node[left] {$\ep h$}	(c)
([yshift=3pt] c.east) edge node[above] {$e-1$}		([yshift=3pt] d.west)
([yshift=-2pt] d.west) edge node[below] {$\ep e$}		([yshift=-2pt] c.east)
(d) edge [loop left, looseness = 5,in=30,out=330] node[right] {$-\ep h$}	(d)
(f) edge node {$\ep$} (e)
(a) edge node {$e$} (c)
(a) edge node {$-h$} (d)
(b) edge node {$1$}(d)
(c) edge node {$1$} (e)
(c) edge node {$h$} (f)
(d) edge node {$e$}(f)
(g) edge node {$\pi$} (h)
(h) edge node {$\sigma$} (i)
;
\end{tikzpicture}
\]

Suppose $(X,\ue,\uh)$ is an $\Ai$-idempotent and $\Phi\colon \Rid\rightarrow
\CS\llbracket \ep \rrbracket$ the corresponding dg functor.  We may regard $\CS\llbracket \ep
\rrbracket$ as a (non-full) subcategory of $\Sigma^\Pi \CS$ via the functor
sending $X\mapsto X\llbracket \ep\rrbracket := \prod_{k\geq 0} \Sigma^{-2k} X$
with the formal endomorphism $\ep$ given by the rightward shift on $X\llbracket
\ep\rrbracket  = X \times \Sigma^{-2} X \times \Sigma^{-4}X\times\cdots$.

Thus the images of the abstract twisted complexes $Z(1,0)$ and $Z(e,h)$ under
$\Phi$ can be viewed as twisted complexes in $\Tw(\adddp \CS)$.  A moment's
thought confirms that these twisted complexes are precisely $Z(X,\Id_X,0)$ from
Example \ref{ex:Z(X,1,0)} and $Z(X,\ue,\uh)$ from Definition \ref{def:image of
Ai idemp}.  To see this, note that half terms in $Z(X,\ue,\uh)$ yield a copy of
$X\llbracket \ep\rrbracket = \prod_{k\geq 0}\Sigma^{-2k}X$, while the other half
yield a copy of $\Sigma\inv X\llbracket \ep\rrbracket = \prod_{k\geq
0}\Sigma^{-2k-1}X$.

Since $Z(X,\Id_X,0)\simeq X$, Lemma \ref{lemma:abstract image} shows that
$Z(X,\ue,\uh)$ is the image of some homotopy idempotent acting on $X$.  It is
not hard to see that this homotopy idempotent is $e_0$, thereby proving the
following.

\begin{proposition}
\label{prop:image-of-idempotent}
  Let $(X,\ue,\uh)$ be an $\Ai$-idempotent in $\CS$.  The object $Z(X,\ue,\uh)$ is
a well-defined twisted complex in $\Tw(\adddp \CS)$; moreover this twisted
complex is an image of the homotopy idempotent $e_0$ acting on $X$.\qed
\end{proposition}

\begin{remark}
\label{rem: truncated z}
The complex $Z(X,\ue,\uh)$ has a natural endomorphism $\ep$ representing the
2-periodicity in this construction. This endomorphism is null-homotopic (by an
explicit homotopy), so $\ep^k$ is null homotopic for all $k\ge 1$ as well. 

The cone of $\ep^k$ is homotopy equivalent to a finite twisted complex (a
truncated version of $Z(X,\ue,\uh)$), which by the above represents the image of
$e_0$ acting on $X\oplus \susp{2k-1}X$. 
\end{remark}

\subsection{From homotopy idempotents to $\Ai$-idempotents}
\label{ss:lifting idempotents}
Finally, we show that any homotopy idempotent $e$ in $\CS$ can be given the
structure of an $\Ai$-idempotent. This is well known to experts, but we will
give an explicit construction of the higher homotopies following ideas of Seidel
\cite[Lemma 4.2]{MR2441780}. See \cite[Propositions 3.2 and 3.4]{MR1214458} for
an alternative proof. We were not able to find explicit formulas for $\ue$ and
$\uh$ in the literature.

\begin{proposition}
  \label{prop:hi-Aii}
Suppose that $e,h\in \End_\CS(X)$ are such that $e$ is degree zero and closed,
and $h$ is degree $-1$ and satisfies $d_\CS(h)=e^2-e$. Then there exist
endomorphisms $h^{(k)}\in \End^{-k}(X)$ for $k\geq 1$ such that $h^{(0)}=1-e$,
$h^{(1)}=h$ and
\[
  d(h^{(k)}) =\sum_{i=0}^{k-1} (-1)^i h^{(i)}h^{(k-1-i)}+\begin{cases}
    -h^{(k-1)}& \text{if}\ k\ \text{is odd}\\
    0 &  \text{if}\ k\ \text{is even}.\\
    \end{cases}  
\]
\end{proposition}

\begin{corollary}
  Suppose that $e,h\in \End_\CS(X)$ are such that $e$ is degree zero and closed,
and $h$ is degree $-1$ and satisfies $d_\CS(h)=e-e^2$. Let $h^{(k)}$ denote the
morphisms obtained from Proposition~\ref{prop:hi-Aii} starting at $(e,-h)$. We define
$e_0:=e$ and $e_k=(-1)^{k+1}h^{(2k)}$ for $k \geq 1$, as well as
$h_k:=(-1)^{k+1} h^{(2k+1)}$ for $k\geq 0$. Then $(X,\ue,\uh)$ is an $\Ai$
idempotent in $\CS$.
\end{corollary}
\begin{proof} It is straightforward that these morphisms satisfy the equations
from Definition~\ref{def:Ainfty idempt}.
\end{proof}

\begin{proof}[Proof of Proposition \ref{prop:hi-Aii}] We will construct
 $h^{(n)}$ inductively. Recall that $h^{(0)}=1-e$. Suppose that we found
 $h^{(1)},\ldots,h^{(n-1)},h_{\mathrm{temp}}^{(n)}$ such that
\begin{align}
  \label{eqn:Ai-idem-induction}
d(h^{(k)})& =\sum_{i=0}^{k-1} (-1)^i h^{(i)}h^{(k-1-i)}+\begin{cases}
-h^{(k-1)}& \text{if}\ k\ \text{is odd}\\
0 &  \text{if}\ k\ \text{is even}.\\
\end{cases},
\quad 1\le k\le n-1
\\
\label{eqn:Ai-idem-inductiontwo}
d(h_{\mathrm{temp}}^{(n)})& =\sum_{i=0}^{n-1} (-1)^i h^{(i)}h^{(n-1-i)}+\begin{cases}
  -h^{(n-1)}& \text{if}\ n\ \text{is odd}\\
  0 &  \text{if}\ n\ \text{is even}.\\
  \end{cases}
\end{align}
For $n=1$, this follows from the assumptions of the proposition if we set $h_{\mathrm{temp}}^{(1)}=h$.
Our goal is to find an $h_{\mathrm{temp}}^{(n+1)}$ and a closed $y_n$, such that
\eqref{eqn:Ai-idem-inductiontwo} will be satisfied for $n\mapsto n+1$ if we set
$h^{(n)}=h_{\mathrm{temp}}^{(n)} +y_n$. Moreover, \eqref{eqn:Ai-idem-induction} will
then hold for $k=n$ since $y_n$ is closed. To this end, we define
\begin{align*}
x_n&= h^{(0)}h_{\mathrm{temp}}^{(n)} +
\sum_{i=1}^{n-1} (-1)^{i} h^{(i)}h^{(n-i)}+ (-1)^n h_{\mathrm{temp}}^{(n)}h^{(0)}+\begin{cases}
  -h_{\mathrm{temp}}^{(n)}& \text{if}\ n\ \text{is even}\\
  0& \text{if}\ n\ \text{is odd}.
\end{cases}
\\
q_n&= -h^{(1)}h_{\mathrm{temp}}^{(n)} +
\sum_{i=2}^{n-1} (-1)^{i} h^{(i)}h^{(n+1-i)}+ (-1)^n h_{\mathrm{temp}}^{(n)}h^{(1)}
\end{align*}

It is not hard to check that $d(x_n)=0$. For $x\in \End_\CS(X)$ we define 
$$A(x)=ex-xe,\ B(x)=ex-x(1-e).$$ 
We then set $y_n=B(x_n)$ if $n$ is even,
and $y_n=A(x_n)$ if $n$ is odd. In either case we again have $d(y_n)=0$. Further, one can
check that  $d(q_n)=A(x_n)$ if $n$ is even, and $d(q_n)=B(x_n)$ if $n$ is
odd, which helps to verify:
\begin{align*}
  B(y_n)=B(B(x_n)) = d(hx_n +3x_n h + 2 q_n e - q_n) + x_n &  \quad \text{ if}\ n\ \text{is even}\\
  A(y_n)=A(A(x_n)) = d(hx_n -3x_n h - 2 q_n e + q_n) + x_n &  \quad \text{ if}\ n\ \text{is odd}.
\end{align*}
Now we set $h^{(n)}=h_{\mathrm{temp}}^{(n)} +y_n$ and
$h_{\mathrm{temp}}^{(n+1)}= -(hx_n +3x_n h + 2 q_n e - q_n)$ if $n$ is even, and
$h_{\mathrm{temp}}^{(n+1)}= -(hx_n -3x_n h - 2 q_n e + q_n)$ if $n$ is odd. Let
us check that \eqref{eqn:Ai-idem-inductiontwo} is now satisfied for $n\mapsto
n+1$. We only consider the case of odd $n+1$, as the other one is analogous. 
\begin{align*}
  d(h_{\mathrm{temp}}^{(n+1)}) = x_n-B(y_n) &=  h^{(0)}h_{\mathrm{temp}}^{(n)} 
 &&+ \sum_{i=1}^{n-1} (-1)^{i} h^{(i)}h^{(n-i)} &&+ h_{\mathrm{temp}}^{(n)}h^{(0)}
    &&-h_{\mathrm{temp}}^{(n)}\\
    &\;\; + h^{(0)}y_n  && &&+y_n h^{(0)} &&-y_n
    \\
    &=  h^{(0)}h^{(n)} &&+
  \sum_{i=1}^{n-1} (-1)^{i} h^{(i)}h^{(n-i)} &&+ h^{(n)}h^{(0)}
    &&-h^{(n)}
\end{align*}
where we have used $-B(y_n) = -e y_n - y_n e + y_n = h^{(0)}y_n + y_n h^{(0)} -y_n$.
\end{proof}

\subsection{The Karoubi envelope}
\label{ss:karoubi and tw}

\begin{definition}\label{def:karoubi envelope}
For a dg category $\CS$ we define the dg Karoubi envelope
$\Kar^{dg}(\CS)$ as the full dg subcategory of $\Tw(\adddp\CS)$ with objects the
twisted complexes homotopy equivalent to $Z(X,\ue,\uh)$ for some
$\Ai$-idempotent $(X,\ue,\uh)$ in $\susp{} \CS$.
\end{definition}

Note that by Example~\ref{ex:Z(X,1,0)}, we have $\CS \hookrightarrow \Kar^{dg}(\CS)$.

\begin{lemma}\label{lemma:kar is karoubian}
The category $\Kar^{dg}(\CS)$ is homotopy idempotent complete.
\end{lemma}
\begin{proof}
If $e_0\in \End_\CS(X)$ is a homotopy idempotent, then $e_0$ admits a lift to an
$\Ai$-idempotent by Proposition \ref{prop:hi-Aii}, and Proposition
\ref{prop:image-of-idempotent} constructs an image of $e_0$.

On the other hand, every object of $\Kar^{dg}(\CS)$ is the image of some
idempotent in $H^0(\CS)$ by construction, and conversely every idempotent in
$\CS$ has an image in $\Kar^{dg}(\CS)$. So if $Y\in \Kar^{dg}(\CS)$ is the image
of a homotopy idempotent $e_0\in \End_\CS(X)$ then all images of all homotopy
idempotents in $\End_{\Kar^{dg}(\CS)}(Y)$ can be constructed as images of some
induced homotopy idempotents $e''_0\in \End_\CS(X)$ by part (2) of Lemma
\ref{lemma:properties of images}.
\end{proof}

\begin{lemma}
\label{lem: Kar homotopy}
There is an equivalence of additive categories $H^0(\Kar^{dg}(\CS))\simeq \Kar
H^0(\CS)$, where the latter denotes the usual idempotent completion of the
additive category $H^0(\CS)$.
\end{lemma}

\begin{proof}
Lemma \ref{lemma:kar is karoubian} implies  that $H^0(\Kar^{dg}(\CS))$ is
idempotent complete in the usual sense, for additive categories.   Now, the
canonical functor $H^0(\CS)\rightarrow H^0(\Kar^{dg}(\CS))$ induces a functor on
the Karoubi envelope $\Kar(H^0(\CS))\rightarrow H^0(\Kar^{dg}(\CS))$ because the
latter is Karoubian.  This functor is essentially surjective because every
object in $\Kar^{dg}(\CS)$ is the image of some homotopy idempotent acting on
some object of $\CS$, and fully faithful by part (3) of Lemma
\ref{lemma:properties of images}.
\end{proof}

\begin{lemma}
\label{lem: Kar qe}
A dg category $\DS$ is homotopy idempotent complete if and only if the natural functor
$\DS\to \Kar^{dg}(\DS)$ is a quasi-equivalence.
\end{lemma}

\begin{proof}
Assume that $H^0(\DS)$ is idempotent complete.  Then $H^0(\Sigma \DS)$ is also
idempotent complete.  By Lemma \ref{lem: Kar homotopy} the natural functor
$\Sigma\DS\to \Kar^{dg}(\Sigma\DS)$ induces an equivalence $H^0(\Sigma\DS)\to
H^0(\Kar^{dg}(\Sigma\DS))\simeq \Kar(H^0(\susp\DS))$, hence $\DS\rightarrow
\Kar^{dg}(\DS)$ is a quasi-equivalence.

For the converse, suppose that the canonical functor $\DS\rightarrow
\Kar^{dg}(\DS)$ is a quasi-equivalence.  Then we have an equivalence of
categories $H^0(\DS)\rightarrow H^0(\Kar^{dg}(\DS))$.  Since the latter category
is idempotent complete, so is the former, i.e.~$\DS$ is homotopy idempotent
complete.
\end{proof}

We are ready to check the universal property of $\Kar^{dg}(\CS)$.

\begin{proof}[Proof of Theorem \ref{thm:Kar}] Let $\DS$ be a homotopy idempotent
complete dg category, and let $F\colon \CS\to \DS$ be a dg functor.  We
extend $F$ to a dg functor $\Tw(\adddp\CS)\rightarrow
\Tw(\adddp)$; this restricts to a dg functor on the full subcategories
$\Kar^{dg}(\CS)\rightarrow \Kar^{dg}(\DS)$.  Since $\DS$ is homotopy idempotent
complete, the canonical dg functor $\DS\rightarrow \Kar^{dg}(\DS)$ is a
quasi-equivalence by Lemma \ref{lem: Kar qe}, and composing with the inverse (in
$\Hqe$) gives a quasi-functor $\widetilde{F}\colon \Kar^{dg}(\CS)\rightarrow \DS$
lifting $F$.

The uniqueness of this lift up to homotopy (again in $\Hqe$) follows because if
$e\in \End_\CS(X)$ is a homotopy idempotent then $F(\im e)$ is determined
uniquely up to homotopy by $F(X)$ and $F(e)$. 
\end{proof}

\begin{lemma}
  \label{lem:twist-of-idem}
  Let $I$ be a finite poset and suppose we have objects $X_i\in \CS$ equipped with
  homotopy idempotents $e_i\in \End^0_\CS(X)$.  Suppose $Y_i\simeq \im e_i$ for
  $i\in I$.  Then any one-sided twist $\tw_\a(\bigoplus_i Y_i)$ is the image of
  some homotopy idempotent $f$ acting on some one-sided twist $\tw_\b(\bigoplus_i
  X_i)$.  
  \end{lemma}
  
  \begin{proof}
  It is sufficient to prove the statement for two-term twisted complexes, i.e.
  cones. Suppose that $Y_0,Y_1$ are homotopy summands of $X_0,X_1$ corresponding
  to homotopy idempotents $e_0$, $e_1$ and $g\colon Y_0\rightarrow Y_1$ is a closed
  degree zero morphism. There are inclusions $\sigma_i\colon Y_i\to X_i$, projections
  $\pi_i\colon X_i\to Y_i$ and homotopies $h_i\colon Y_i\to Y_i$ such that
  $\pi_i\sigma_i=\Id_{Y_i}+d(h_i)$. Define $g':=\sigma_{1}g\pi_0$,  and consider
  the following chain maps between $\Cone(g)$ and $\Cone(g')$:
  $$
  p:=\begin{tikzcd}
  X_0 \arrow{d}{\pi_0}\arrow{r}{g'} \arrow{dr}{H} & X_1\arrow{d}{\pi_1}\\
  Y_0  \arrow{r}{g} & Y_1
  \end{tikzcd},
\qquad
  s:=\begin{tikzcd}
  X_0 \arrow{r}{g'} & X_1\\
  Y_0 \arrow{u}{\sigma_0}\arrow{r}{g} \arrow{ur}{H'} & Y_1 \arrow{u}{\sigma_1}
  \end{tikzcd}
  $$
  where $H=h_1g\pi_0$ and $H'=\sigma_1gh_0$. It is easy to see that $p\circ s$
  is homotopic to identity, so $\Cone(g)$ is a homotopy summand of $\Cone(g')$.
  \end{proof}

\begin{theorem}
If $\CS$ is pretriangulated, then so is $\Kar^{dg}(\CS)$.
\end{theorem}
\begin{proof} 
Every object of $\Kar^{dg}(\CS)$ is (isomorphic to) a homotopy summand of an object in
$\CS$.  A one-sided twisted complex constructed from homotopy summands of $Y^i\in \CS$
($i\in I$) is also a homotopy summand of a twisted complex constructed from
$Y^i$ by Lemma~\ref{lem:twist-of-idem}; such objects live in $\Kar^{dg}(\CS)$
since $\CS$ is pretriangulated and $\Kar^{dg}(\CS)$ is idempotent complete.
\end{proof}

\begin{remark}
\label{rem: Ai Karoubi}
All the above constructions naturally extend to $\Ai$-categories, see \cite{MR2441780} for details. In particular, the equation for an $\Ai$-idempotent in Proposition \ref{prop:hi-Aii} should be replaced with 
$$
\sum_d \sum_{i_1,\ldots,i_d} \pm \mu_d(h^{(i_1)},\ldots,h^{(i_d)})=\begin{cases}
h^{(k-1)} & \text{if}\ k\ \text{is odd},\\
0 & \text{if}\ k\ \text{is even},
\end{cases}
$$
where the sum is over all partitions $i_1+\ldots+i_d=k-d+1$.  

The analogue of Proposition \ref{prop:hi-Aii} is \cite[Lemma 4.2]{MR2441780}
which states that any homotopy idempotent can be lifted to an $\Ai$-idempotent.
This allows one to define the Karoubi completion $\Kar(\CS)$ of an $\Ai$-category
$\CS$, which is again an $\Ai$-category and satisfies the universal property as
above.
\end{remark}

\subsection{Perfect complexes are Karoubian}
\label{ss:perf definitions}
Let $A$ be a dg algebra. Throughout this section we will assume that $A$ is
supported in {\em non-positive} cohomological degrees. 

Let $\langle A\rangle$ to be the category of finitely generated free
$A$-modules, that is, direct sums of finitely many copies of $A$. We define
$\pretr{A}$ to be the category of (finite) twisted complexes of free $A$-modules
or, equivalently, the pretriangulated hull of $\langle A\rangle$. We define
$\Kar^{dg}\langle A\rangle$ to be the category of projective $A$-modules, that
is, homotopy direct summands of objects in $\langle A\rangle$. Let $\Perf(A)$ be
the category of {\em perfect twisted complexes} over A, that is, bounded twisted
complexes built out of objects in $\Kar^{dg}\langle A\rangle$. In other words,
$$
\Perf(A)=\pretr{\Kar^{dg}\langle A\rangle}.
$$
Finally, we define $\Ch^b(A)$ (resp. $\Ch^b(\Kar^{dg}\langle A\rangle)$)
to be the category of bounded complexes of free (resp.
projective) $A$-modules.   

By definition, an object of $\pretr{A}$ is a graded free dg module $M=\oplus_i
\susp{-i}M_i$  with $M_i$ in $\langle A\rangle$, and a differential given by a
matrix $d=(d_{ij})$ where $d_{ij}\colon M_j\to M_i$ and $d_{ii}$ agrees with the
internal differential on $M_i$. The differential $d$ is homogeneous of
cohomological degree $1$, so $d_{ij}$ has cohomological degree $1+j-i$. Since
$A$ is non-positively graded, we have $d_{ij}=0$ unless $i\ge j$,
 so the matrix $d=(d_{ij})$ is 
lower-triangular. Also, $d^2=0$ which means
\begin{equation}
\label{eq: twisted d square}
\sum_{k\colon i\ge k\ge j} d_{ik}d_{kj}=0\ \text{for all}\ i> j.
\end{equation}

Similarly, if $M=(\oplus_i \susp{-i}M_i,d)$ and $N=(\oplus_i \susp{-i}N_i,d)$ then a morphism 
from $M$ to $N$ is given by an lower-triangular matrix of morphisms $f=(f_{ij})$ where
$f_{ij}$ is a morphism from $M_j$ to $N_i$ of cohomological degree $j-i$.

This immediately implies the following:

\begin{proposition}
\label{prop: epsilon}
 Suppose that the differential on $A$ vanishes
 and let
$(M,d)$ be a twisted complex in $\pretr{A}$, and 
$d'=\sum_{i}d_{i+1,i}$. Then $(d')^2=0$, so $(M,d')$ is a well-defined
chain complex in $\Ch^b(A)$.  
\end{proposition}

\begin{remark}
A similar construction works if the differential on $A$ is non-trivial, but its
component from $A_{-1}$ to $A_0$ vanishes. 
\end{remark}

\begin{remark}
We can define a similar construction for an $\Ai$-algebra $A$ with $\mu_1=0$. Indeed, the $\Ai$-Maurer-Cartan equation for $d$
still implies $d_{i+2,i+1}d_{i+1,i}=0$.
\end{remark}

We will call $(M,d')$ the underlying
complex of the twisted complex $(M,d)$.  
This construction defines dg functors
$$
\epsilon\colon \pretr{A}\to \Ch^b(A),\quad \epsilon\colon \Perf(A)\to \Ch^b(\Kar^{dg}\langle A\rangle).
$$
We will denote both of these functors by $\epsilon$ since it will be clear from
the context which one is used.

 The following lemma is a standard application of perturbation theory. 
 \begin{lemma}
 \label{lem:epsilon homotopy equivalence}
 Suppose that $X$ is an object in $\pretr{A}$ or $\Perf(A)$.
 If $\epsilon(X)$ is contractible then $X$ is contractible. If $f\colon X\to Y$ is a morphism 
 such that $\epsilon(f)$ is a homotopy equivalence, then $f$ is a homotopy equivalence.
 \end{lemma}

 \begin{proof}
 For the first part see e.g \cite[Lemma 2.19]{GH}. For the second part, observe
 that $f$ is a homotopy equivalence if and only if the cone of $f$ is
 contractible. 
 \end{proof}

\begin{theorem}
\label{th: perf complete}
Assume that $A$ is non-positively graded and the differential on $A$ vanishes.
Then the category $\Perf(A)$ is homotopy idempotent complete, and 
$$
\Kar^{dg}(\pretr{A})=\Perf(A).
$$
\end{theorem}

\begin{remark}
The proof is similar to \cite[Lemma 1.5.6(iii)]{BeiVol} and
\cite[Appendix]{GW}, but we include it here for completeness
and add more details for reader's convenience.
\end{remark}

\begin{proof}
Let $X$ be a twisted complex in $\Perf(A)$, and let $e$ be a homotopy idempotent
endomorphism of $X$. Without loss of generality we can suppose that $X$ is a
(one-sided) twist of
$$
\bigoplus_{0<i<(2N-1)} \susp{-i}X_i
$$
for some sufficiently large $N$.

As in Section \ref{ss:lifting idempotents}, we can lift $e$ to an $A_{\infty}$
idempotent $(X,\ue,\uh)$ and consider a twisted complex $Z$ built out of several
copies of $X$ connected by the maps from $\ue$ and $\uh$ with appropriate signs.
As in Remark \ref{rem: truncated z} we build a bounded
$Z$ from $2N$ of copies of $X$, which is homotopy equivalent to the image of $e$
acting on $X\oplus \susp{2N-1}X$ via an argument similar to Lemma
\ref{lemma:abstract image}. By our assumptions, $X$ is a direct sum of copies of
$\susp{-i} A$ with $i>0$, and $\susp{2N-1}X$  is a direct sum of copies of
$\susp{-i} A$ with $i<0$. There are no nonzero $A$-module maps from $\susp{-j}
A\rightarrow \susp{-i} A$ unless $i>j$, so we get the following diagram:
\begin{center}
\begin{tikzcd}
  \susp{2N-1}X \arrow{d}\arrow{dr} \arrow[bend left]{drr} &     & X \arrow{d}\\
Z_{<0}\arrow{r} \arrow{d} \arrow{drr} \arrow[bend left]{rr}   & Z_0 \arrow{r} \arrow{dr} & Z_{>0} \arrow{d}\\
\susp{2N-1}X &     & X\\
\end{tikzcd} 
\end{center}
where the top half of the diagram represents the projection $\pi\colon X\oplus
\susp{2N-1}X\to Z$ and the bottom half represents the inclusion $i\colon Z\to
X\oplus \susp{2N-1}X$.  We decomposed $Z=\bigoplus_a \Sigma^{i_a} A$ into its
summands $Z_{>0}$, $Z_0$, $Z_{<0}$ corresponding to those indices $a$ for which
$i_a>0$, $i_a=0$, or $i_a<0$, respectively (recall that $\Sigma$ is the
suspension so in fact lowers cohomological degree).

Observe that $\pi\circ i$ is homotopic to identity on $Z$ and vanishes on $Z_0$. 
Since the differential on $A$ vanishes, by Proposition \ref{prop: epsilon} one
can define the forgetful dg functor $\epsilon\colon \Perf(A)\to \Ch(\Kar^{dg}(A))$. By
applying it to $Z$ we get a chain complex $\epsilon(Z)$ and a homotopy $h\in
\End(\epsilon(Z))$ such that $dh+hd|_{Z_0}=\Id_{Z_0}$. Now 
$$
(dh)^2=(dh+hd)dh=dh,\ (hd)^2=hd(dh+hd)=hd,
$$
so $dh$ and $hd$ are two orthogonal idempotents on $Z_0$. Since $\Kar^{dg}\langle
A\rangle$ is Karoubian, we can split $Z_0=Q\oplus Q'$.

For the final step, let $T$ be the natural twisted complex built out of $Q$ and
$Z_{>0}$ (abstractly speaking, this is the cone of the composition
$\susp{-1}Q\hookrightarrow \susp{-1}Z_0\to Z_{>0}$).  It is easy to see that $T$
is a subcomplex of $Z$ and we can restrict the maps $i$ and $\pi$ from the above
diagram to $i'\colon T\to X$ and $\pi'\colon X\to T$. Then $i'\circ \pi'=i\circ
\pi|_{X}\simeq e$ and $\epsilon(\pi'\circ i')\simeq \Id_{\epsilon(T)}$. By Lemma
\ref{lem:epsilon homotopy equivalence} we get that $\pi'\circ i'$ is homotopic
to $\Id_{T}$. Therefore $T$ represents the image of $e$, and we conclude that
$X$ is homotopy equivalent to a perfect $A$-module.

This shows that $\Perf(A)$ is homotopy Karoubian.
\end{proof}

Suppose that a dg algebra $A$ retracts to its dg module $V$. Recall that in section
\ref{subsec:Ainf} we defined $\Ai$-operations $\mu_k\colon V^{\otimes k}\to V$.  

\begin{lemma}
\label{lem:perf quasi equivalent}
Define $\Perf(V)=\pretr{\Kar\langle V\rangle}$ where $\langle V\rangle$ is the $\Ai$-category with a single object with endomorphism $\Ai$-algebra $V$.
Then the $\Ai$-categories $\Perf(A)$ and $\Perf(V)$ are quasi-equivalent.
\end{lemma}

\begin{proof}
We can consider $A$ and $V$ as
$\Ai$-categories $(*,A)$ and $(*,V)$ with one object with endomorphism algebra
$A$ or $V$, respectively. We claim that $(*,A)$ and $(*,V)$ are
quasi-equivalent.

Indeed, there exist $\Ai$-algebra homomorphisms between $A$ and $V$
defined by $\Ai$-maps $\lambda_k:V^{\otimes k}\to A$  from the proof of Theorem \ref{thm:Ainfty-defretracts}. By Theorem
\ref{thm:Ainfty-defretracts}), these $\Ai$-algebra homomorphisms induce
quasi-isomorphisms on homology, so $(*,A)$ and $(*,V)$ are quasi-equivalent as
$A_{\infty}$ categories. 

The corresponding categories of free dg
modules $\langle A\rangle$ and $\langle V\rangle$ are nothing but additive closures of 
$(*,A)$ and $(*,V)$, so these are quasi-equivalent.
Therefore $\Kar^{dg}\langle A\rangle$ and $\Kar^{dg}\langle V\rangle$ are
quasi-equivalent by Theorem \ref{thm:Kar} and their pretriangulated hulls are
quasi-equivalent as well. 
\end{proof}

\section{The bar complex for dg categories}
\label{sec:bardg}
One can motivate the introduction of the (two-sided) bar complex of a dg
category as the object which governs the notion of ``naturality up to coherent
homotopy''.  To illustrate, let $\CS$ and $\DS$ be dg categories and
$F,G:\CS\rightarrow \DS$ dg functors.  We have already introduced the notion of
a (say, degree zero) natural transformation, which is a choice of morphism
$\a_X\in \Hom_\DS(F(X),G(X))$ for all $X\in \CS$, natural with respect to
morphisms in $\CS$. This means for every morphism $f\colon X \to Y$ in $\CS$ we have
\[G(f) \circ \a_X - \a_Y \circ F(f) = 0.\] 
Said differently, the pair of dg functors $F,G$ determines a
$\CS,\CS$-bimodule $\BS(G,F)$, which for $X,Y\in \CS$ has
\[
Y \BS(G,F) X:= G(Y)\DS F(X),
\]
(in the notation of \S \ref{ss:bimodules}) and a natural transformation is
simply a map of $\CS,\CS$-bimodules
\[
\a\colon \CS \rightarrow \BS(G,F).
\]
The image of $\Id_X\in \CS$ is the chosen morphism $\a_X$, and naturality is
equivalent to $\a$ being a $\CS,\CS$-bimodule map.

Now suppose that we are in a situation
where naturality does not hold on the nose, but only up to homotopies $\a_f\in
\Hom_{\DS}^{-1}(F(X),G(Y))$:
\[d_{\DS}(\a_f) = G(f) \circ \a_X - \a_Y \circ F(f).\] Let us consider another
morphism $g \colon Y \to Z$ in $\CS$ and suppose that $f$ and $g$ are closed.
It is straightforward to check that the expression
\[\alpha_{g \circ f} - G(g)\circ \alpha_{f} - \alpha_{g} \circ F(f)\] is closed.
It is often desirable to suppose this expression is also exact, i.e.~ there
exists a homotopy $\a_{g,f}\in \Hom_{\DS}^{-2}(F(X),G(Z))$.  Roughly speaking,
this says that the assignment $f\mapsto \a_f$ satisfies a version of the Leibniz
rule, up to homotopy.

Now, assuming the existence of such ``higher homotopies'' $\a_{f,g}$ for all
composable morphisms $f,g$ one can then define a degree $-3$ closed morphism
associated to each triple of composable morphisms $f_1,f_2,f_3$, the exactness
of which would allow us to define a family of closed morphisms (obstructions)
associated to each 4-tuple of composable morphisms, and so on.

If all such obstructions are exact (and a family of homotopies realizing this
exactness is given) then then we say that the system
$(\a_X,\a_f,\a_{g,f},\cdots)$ is a \emph{homotopy coherent natural
transformation}.

The two-sided bar complex $\BB(\CS)$ of $\CS$, which we describe explicitly
below, can be considered as a free resolution of the trivial $\CS,\CS$-bimodule
$\CS$, and is spanned as a bimodule by sequences of composable morphisms of arbitrary finite
length $r\geq 0$. The data of a homotopy coherent natural transformation
$F\rightarrow G$ is then encoded as a map of $\CS,\CS$-modules 
\[
\BB(\CS)\rightarrow \BS(G,F).
\]
Various operations on natural transformations can be understood in terms of
various structures on the two-sided bar complex; for instance the composition of natural
transformations can be understood via a natural comultiplication on $\BB(\CS)$.

After this motivation, we now give an explicit description of the two-sided bar complex.

\subsection{The two-sided bar complex of a pair}
\label{ss:bar of pair}

Let $\CS$ be a dg category and $\IS\subset \CS$ a (unital) subcategory. We wish
to define the two-sided bar complex associated to $(\IS,\CS)$.  First, consider
$\CS,\CS$-bimodules of  the form $\CS\otimes_\IS \cdots \otimes_\IS \CS$.

\begin{example}
In case $\IS = \IS_\CS$ is the subcategory of identity maps from
Example~\ref{exa:subcat-idmaps}, $\CS\otimes_\IS \cdots \otimes_\IS \CS$ is
spanned by composable morphisms in $\CS$.  More generally, if $\BS\subset \CS$
is a full subcategory and $\IS=\IS_\BS\subset \CS$ is the subcategory of
identity morphisms in $\BS$, then
\[
\underbrace{\CS\otimes_\IS\otimes \cdots \otimes_\IS \CS}_{r+2 \text{ factors}}
\]
is spanned as a bimodule by sequences of composable morphisms, which we denote by $f_0\bbar
f_1\bbar \cdots \bbar f_{r+1}$, in which $f_1,\ldots,f_r$ are morphisms in the
full subcategory $\BS\subset \CS$.
\end{example}

The \emph{two-sided bar complex} associated to the pair $\IS\subset \CS$ is by definition
the total complex  $\BB(\CS;\IS)$ of the following bicomplex:
\begin{equation}\label{eq:bar cx}
\begin{tikzpicture}[baseline=0em]
\tikzstyle{every node}=[font=\small]
\node (a) at (0,0) {$\CS\otimes_\IS \CS$};
\node (b) at (-4,0) {$\CS\otimes_\IS \CS\otimes_\IS \CS$};
\node (c) at (-8,0) {$\CS\otimes_\IS\CS\otimes_\IS\CS\otimes_\IS \CS$};
\node (d) at (-12,0) {$\cdots$};
\path[->,>=stealth,shorten >=1pt,auto,node distance=1.8cm,
  thick]
(b) edge node 	{} (a)
(c) edge node  {} (b)
(d) edge node {} (c);
\end{tikzpicture}
\end{equation}
in which the horizontal arrows are
\[
f_0\bbar \cdots \bbar f_{r+1}\mapsto \sum_{i=0}^r (-1)^i \;f_0\bbar\cdots \bbar f_i\circ f_{i+1}\bbar\cdots \bbar f_{r+1}.
\]

More precisely,
\[
\BB(\CS;\IS)= \bigoplus_{r\geq 0} \bigoplus_{X_1,\ldots,X_{r+1}} \susp{r}\Big(\CS X_{1}\otimes  
X_1\CS X_2\otimes\cdots \otimes  X_{r+1}\CS\Big),
\]
where the direct sum is over objects $X_1,\ldots,X_{r+1}$ for which
$\Id_{X_i}\in \IS$ for $i=1,\ldots,r+1$.  The differential on the two-sided bar complex is
$d_{\oplus} + d_{\text{bar}}$, where $d_{\oplus}$ is the direct sum of
differentials on the terms above (inherited from $d_{\CS}$ by the usual tensor
product rule, and a sign $(-1)^r$ coming from the $r$-translation), and
$d_{\text{bar}}$ is defined by
\[
d_{\text{bar}}(f_{0}\bbar\cdots\bbar f_{r+1}) = 
\sum_{i=0}^r (-1)^{i} \;f_{0}\bbar \cdots\bbar f_i\circ f_{i+1}\bbar \cdots \bbar f_{r+1}.
\]
(The notation is suggestive: the bar differential is the alternating sum of erasing bars.)
The bimodule structure is given term-wise by
\[
g\cdot (f_{0}\bbar f_1\bbar \cdots \bbar f_r \bbar f_{r+1}) \cdot g' \ :=
 \ (-1)^{|g|r} (g\circ f_{0})\bbar f_1\bbar \cdots \bbar f_r  \bbar (f_{r+1}\circ  g').
\]
Note that the left-action of bimodules obtains a sign-twist under
translation.

By convention $\BB(\CS)$ denotes $\BB(\CS,\IS)$ where $\IS=\IS_\CS$ is the
subcategory of identity maps in $\CS$. 

\begin{remark} At the beginning of this section we have motivated the bar
complex by claiming that $\CS,\CS$-bimodule maps $\alpha\colon \BB(\CS)\to
\BS(G,F)$ encode homotopy coherent natural transformations between dg functors
$F,G \colon \CS \to \DS$. Indeed, the data of such a bimodule map is precisely
the choice of closed degree zero morphisms $\a_X:=\a(\Id_X \!\bbar\! \Id_X)$
from $F(X)$ to $G(X)$ in $\DS$, together with morphisms $\a(\Id_Y \!\bbar f
\bbar\! \Id_X)$ which realize naturality up to homotopy, and the higher
homotopies $\a(\Id_{X_0}\bbar \! f_1,\ldots,f_r \bbar\! \Id_{X_{r+1}})$ which
provide the required higher naturality data.
\end{remark}

\begin{proposition}
The two-sided bar complex $\BB(\CS,\IS)$ is a coalgebra object in the category of
$\CS,\CS$-bimodules:  we have maps of bimodules
$\BB(\CS,\IS)\rightarrow \CS$ and $\BB(\CS,\IS)\rightarrow \BB(\CS,\IS)\otimes_\CS
\BB(\CS,\IS)$ satisfying the usual counit and coassociativity relations.
\end{proposition}
\begin{proof}
The counit $\e\colon \BB(\CS,\IS)\rightarrow \CS$ is defined componentwise by
\[
f_0\bbar \cdots \bbar f_{r+1} \mapsto \begin{cases}
f_0\circ f_1 & \text{ if $r=0$}\\
0 & \text{ otherwise}
\end{cases}.
\]
The comultiplication $\Delta\colon \BB(\CS,\IS)\rightarrow \BB(\CS,\IS)\otimes_\CS
\BB(\CS,\IS)$ is defined componentwise by
\[
f_0\bbar \cdots \bbar f_{r+1} \mapsto \sum_{i=1}^{r+1} (-1)^{(r-i+1)(|f_0|+\cdots+|f_{i-1}|)} 
(f_0\bbar \cdots\bbar f_{i-1} \bbar \Id)\otimes (\Id\bbar f_i \bbar\cdots\bbar f_{r+1}).
\]
The counit and coassociativity axioms are easily checked.
Furthermore, it is clear that $\Delta$ and $\e$ are closed and commute with the
$\CS,\CS$-bimodule structure on $\BB(\CS,\IS)$.
\end{proof}

Let us give an alternate description of the two-sided bar complex.  The bimodule
$\CS\otimes_\IS \CS$ is a coalgebra object $C_\IS$ in $\CS,\CS$-bimodules, and
$\BB(\CS,\IS)$ is the counital idempotent $\PB_{C_{\IS}}$ associated to
$C_{\IS}$ \cite{hogancamp2017idempotents}.

If $\IS\subset \JS\subset \CS$, then the counit of $C_{\IS}$ factors through the
counit of $C_{\JS}$, so $C_\IS\leq C_{\JS}$ in the notation of
\cite{hogancamp2017idempotents}.  The following is an immediate consequence of
this.

\begin{proposition}
If $\IS\subset \JS\subset \CS$ are unital subcategories then
\[
\BB(\CS,\IS)\otimes_\CS \BB(\CS,\JS)\simeq \BB(\CS,\IS)\simeq \BB(\CS,\JS)\otimes_\CS \BB(\CS,\IS).
\]
\end{proposition}

\subsection{Shrinking the bar complex}

\begin{definition}\label{def:generates}
We say that a full subcategory $\BS\subset \CS$ \emph{generates $\CS$} if for
all dg functors $F:\CS\rightarrow \k\text{-dgmod}$,  $F(X)\simeq 0$ for all
$X\in \BS$ implies $F(X)\simeq 0$ for all $X\in \CS$.
\end{definition}

\begin{proposition}\label{prop:B generates if it generates}
Suppose every object of $\CS$ is homotopy equivalent to a finite one-sided
twisted complex constructed from objects of $\BS$.   Then $\BS$ generates $\CS$
in the sense of Definition \ref{def:generates}.
\end{proposition}
\begin{proof}
Let $F:\CS\rightarrow \k\text{-dgmod}$ be a dg functor such that $F(X)\simeq 0$
for all $X\in \BS$, and let $Z = \tw_\d(\bigoplus_i \susp{a_i}X_i)$ be a finite
one-sided twisted complex.  Then $F(Z)$ is isomorphic to a one-sided twisted
complex $\tw_{F(\d)}(\bigoplus_i \susp{a_i}F(X_i))$.  It follows that $F(Z)\simeq 0$
since finite one-sided twisted complexes built from contractible complexes are
contractible.
\end{proof}
\begin{remark}
The statement of Proposition \ref{prop:B generates if it generates} remains true
if we replace ``finite'' with ``bounded above''.
\end{remark}

\begin{proposition}
\label{prop:barcx-retract}
If $\BS\subset \CS$ generates, then the natural inclusion
$\BB(\CS,\IS_{\BS})\rightarrow \BB(\CS,\IS_{\CC})$ is the section of a
deformation retract.
\end{proposition}
\begin{proof}
Let $\IS=\IS_{\CS}\subset \CS$ be the subcategory of identity morphisms in
$\CS$, and let $\JS=\IS_\BS\subset \CS$ be the subcategory of identity morphisms
in $\BS$.  Let $C=\CS\otimes_\IS\CS$ and $D=\CS\otimes_\JS \CS$.

Note that $C = D\oplus E$, where $E=\bigoplus_{Y\not\in \BS} \CS Y\otimes Y\CS$.

Consider the $\CS,\CS$-bimodules $\PB_C,\PB_D$ and $\AB_D=\Cone(\PB_D\rightarrow
\CS)$.  Since $C = D\oplus E$, it follows that $\PB_C$ can be expressed as a
one-sided twisted complex constructed from one copy of $\PB_D$ and terms of the
form
\[
\AB_D\otimes_\CS E\otimes_\CS \AB_D\otimes_\CS \cdots \otimes_\CS E\otimes_\CS \AB_D.
\]
But one of the features of $\AB_D$ is $Y\AB_D\simeq 0\simeq \AB_D Y$
for all $Y\in \BS$.  Since $\BS$ generates, it follows that $Y\AB_D\simeq 0
\simeq \AB_DY$ for all $Y\in \CS$.  In particular $E\otimes_\CS \AB_D\simeq 0
\simeq \AB_D\otimes_\CS E$.  Contracting the contractible terms yields the
desired deformation retract
\[
\PB_C = \PB_{D\oplus E} \rightarrow \PB_D. \vspace{-.6cm}
\]
\end{proof}

\subsection{The bar complex of the pretriangulated hull}
\label{ss:bar of pretr}

Let $\CS$ be a dg category and $\pretr{\CS}$ its pretriangulated hull.  We have
a homotopy equivalence of bimodules
\[
\BB(\pretr{\CS},\IS_{\pretr{\CS}})\rightarrow \BB(\pretr{\CS},\IS_{\CS})
\]
implicitly constructed in the previous section.   It will be useful to
understand this map explicitly.

We first consider the relation between the two-sided bar complexes of $\CS$ and $\Sigma \CS$.
There is a map of $\Sigma\CS,\Sigma\CS$-bimodules

\[
\Phi\colon \BB(\Sigma\CS,\IS_{\Sigma\CS})\rightarrow \BB(\Sigma\CS,\IS_{\CS})
\]
defined as follows: given objects $X_0,\ldots,X_r\in \Sigma\CS$ of the form
\[
X_i = \bigoplus_j \susp{a_{i,j}}X_{ij}
\]
and a sequence of composable morphisms
\[
g_1\otimes \cdots\otimes g_r \ \in X_0(\Sigma\CS) X_1\otimes \cdots \otimes X_{r-1} (\Sigma\CS) X_{r},
\]
we define $\Phi(\Id\bbar g_1\bbar\cdots\bbar g_r\bbar \Id)$ to be the sum of
terms of the form
\[
\pm \sigma_{j_0}\bbar (g_1)_{j_0,j_1}\bbar  (g_2)_{j_1,j_2}\bbar\cdots \bbar  (g_r)_{j_{r-1},j_r} \bbar\pi_{j_r},
\]
where $(g_i)_{j_{i-1},j_i}$ indicates the component of $g_i$ living in
$X_{i-1,j_{i-1}} \CS X_{i,j_i}$ and $\sigma_{j_0}$ denotes the inclusion of
$X_{0,j_0}$ into $X_0$ and $\pi_{j_r}$ denotes the projection of $X_r$ onto
$X_{r,j_r}$.   It is an exercise to find the correct signs such that $\Phi$
defines a deformation retract.

Next we consider the relation between the two-sided bar complex of $\CS$ and $\Tw(\CS)$.
We define a map of $\CS,\CS$-bimodules
\[
\Psi\colon \BB(\Tw(\CS),\IS_{\Tw(\CS)})\rightarrow \BB(\Tw(\CS),\IS_{\CS})
\]
as follows.  Given objects $\tw_{\a_i}(X_i)\in \Tw(\CS)$ ($0\leq i\leq r$) and
an element
\[
\Id\!\bbar f_1\bbar \cdots \bbar f_r \bbar \! \Id \ \ \in  \ \ \BB(\Tw(\CS),\IS_{\Tw(\CS)}),
\]
we define $\Psi(\Id\!\bbar f_1\bbar \cdots \bbar f_r \bbar \! \Id)$ to be the
sum of terms of the form
\[
\pm \phi \bbar \underbrace{\a_0\bbar\cdots\bbar \a_0}_{i_0} \bbar f_1 
\bbar \underbrace{\a_1\bbar\cdots\bbar \a_1}_{i_1} \bbar 
\cdots 
\bbar \underbrace{\a_{r-1}\bbar\cdots\bbar \a_{r-1}}_{i_{r-1}} 
\bbar f_r \bbar \underbrace{\a_r\bbar\cdots\bbar \a_r}_{i_r} \bbar \phi\inv
\]
where $i_0,\ldots,i_r\geq 0$ and $\phi$ and $\psi$ denote $\Id_{X_0}$ and
$\Id_{X_r}$, regarded as degree zero (not necessarily closed) maps
$X_0\rightarrow  \tw_{\a_0}(X_0)$ and $X_r\rightarrow \tw_{\a_r}(X_r)$.  It is
an exercise to find the correct signs and check that this defines a deformation retract.

Combining these gives a bimodule map relating the two-sided bar complex of $\CS$
and its pretriangulated hull; this map is a deformation retract and, in
particular, a homotopy equivalence.

\subsection{The semi-orthogonal bar complex}
\label{ss:semiortho}
\begin{definition}\label{def:semiortho}
Let $\Gamma$ be a  poset, and let $\BS_\gamma\subset \CS$ be full dg
subcategories of $\CS$, indexed by $\gamma\in \Gamma$ such that
\begin{enumerate}
\item every object $Z\in \CS$ is homotopy equivalent to a one-sided twisted
complex $\tw_\d(\bigoplus_{\gamma\in \Gamma}Z_\gamma)$ where $Z_\gamma\in
\BS_\gamma$, only finitely many $Z_{\gamma}$ are nonzero and the twist $\d$ is
strictly lower triangular with respect to the partial order on $\Gamma$.
\item $\Hom_\CS(Y_{\gamma},Y_{\gamma'})$ is contractible unless $\gamma\leq
\gamma'$.
\end{enumerate}
Then we say that $\BS_{\gamma}$ defines a $\Gamma$-indexed semi-orthogonal
decomposition of $\CS$.
\end{definition}

Note that if $X\in \BS_{\gamma}\cap \BS_{\gamma'}$ for $\gamma\neq \gamma'$ then
$X$ is contractible.  Thus, we usually assume, without loss of generality, that
the subcategories $\BS_{\gamma}$ have no objects in common.

Let $\{\BS_\gamma\}_{\gamma\in \Gamma}$ be a semi-orthogonal decomposition of $\CS$.  
Define the full dg subcategory $\BS:=\bigcup_{\gamma \in \Gamma}\BS_\gamma \subset \CS$, 
which generates $\CS$ by Proposition~\ref{prop:B generates if it generates}.

For $X,Y\in \BS$ we write $X\leq Y$ if $X\in \BS_\gamma$ and $Y\in \BS_{\gamma'}$ with
$\gamma\leq \gamma'$.  This relation is transitive and reflexive, but not
anti-symmetric: $X\leq Y$ and $Y\leq X$ holds if and only if $X,Y\in \BS_\gamma$
for some $\gamma$.

The \emph{semi-orthogonal two-sided bar complex} of $\CS$ is defined to be the subcomplex 
$\BB^\Gamma(\CS,\IS_\BS)\subset \BB(\CS,\IS_\BS)$ spanned by elements of the form
\[
f_0\bbar \cdots \bbar f_{r+1} \ \ \textrm{in} \ \ 
\CS X_1 \otimes X_1\CS X_2\otimes\cdots\otimes X_{r+1}\CS
\] 
with $X_1\geq X_2\geq \cdots \geq X_{r+1}$.

\begin{proposition}\label{prop:semiortho bar}
Retain notation as above. The inclusion $\BB^\Gamma(\CS,\IS_\BS)\hookrightarrow
\BB(\CS,\IS_\BS)$ is the section of a deformation retract;  the data of this
deformation retract are $\CS,\CS$-bilinear.
\end{proposition}
\begin{proof}
The complex $\BB(\CS,\IS_\BS)$ is a one-sided twist of 
\[
\bigoplus_{r\geq 0} \bigoplus_{X_1,\ldots,X_{r+1}\in \BS} 
\susp{r} \Big( \CS X_1\otimes_\k X_1\CS X_2\otimes_\k \cdots\otimes_\k  X_{r+1}\CS\Big)
\]
But $X_i \CS X_{i+1} = \Hom_{\CS}(X_{i+1},X_i)$ is contractible unless $X_i\geq
X_{i+1}$.  Contracting all such terms gives the desired deformation retract.
Since the contractions are all of the form
\[
\CS X_1\otimes_\k \underbrace{X_1\BS X_2\otimes_\k \cdots\otimes_\k X_r \BS X_{r+1}}_{\text{contract some factor}}\otimes X_{r+1}\CS,
\]
(i.e.~we only ever contract an ``internal'' tensor factor) the data of the
deformation retract commute with the left and right $\CS$-actions. 
\end{proof}

Since $\BS$ generates $\CS$, Proposition~\ref{prop:barcx-retract} implies the
following.

\begin{theorem}
\label{thm:barcx-retract}
Retain notation as above. The full bar complex $\BB(\CS)$ deformation retracts 
onto the semi-orthogonal bar complex $\BB^\Gamma(\CS,\IS_\BS)$. The data of this
deformation retract are $\CS,\CS$-bilinear.
\end{theorem}

\begin{corollary}
The semi-orthogonal bar complex $\BB^{\Gamma}(\CS,\IS_{\BS})$ has the structure
of an $\Ai$-algebra object in $\CS,\CS$-bimodules, c.f.
Remark~\ref{rem:Ainfty-structure-linear}.\qed
\end{corollary}

\subsection{Hochschild homology (vertical trace)}
\label{ss:vert trace}

If $\AS$ is an additive $\k$-category then the (vertical) trace of $\AS$ is the
$\k$-module 
\[
\HH_0(\AS) := \bigoplus_{X\in \AS} \End(X) \ \Big / \  \mathrm{span}_{\k}\{ g_0\circ g_1-g_1 \circ g_0\}
\]
where $(g_0,g_1)$ runs over pairs of morphisms which are composable in either order.

\begin{notation}
The class of $f\in \End(X)$ in $\HH_0(\AS)$ will be denoted $[f]$.  We also write
$[X]:=[\Id_X]$ for the associated class in $\HH_0(\AS)$.
\end{notation}

\begin{remark}
If $\AS$ is monoidal, then $\HH_0(\AS)$ inherits the structure of an algebra via
$[f] \cdot [g]:= [f\otimes g]$.
\end{remark}

It is natural to consider a derived version of the vertical trace, in which the
relation $g_0\circ g_1 = g_1\circ g_0$ is not strictly enforced, but is achieved
by the formal adjunction of a homotopy $h(g_0,g_1)$.  In this setting it is also
natural to allow our additive category to be an arbitrary dg category $\CS$.
After adjoining such homotopies, symbols of the form $h(g_0g_1,
g_2)-h(g_0,g_1g_2) + h(g_2g_0,g_1)$ are automatically closed. In order avoid
creating such new closed elements, one is forced to adjoin higher homotopies
$h(g_0,g_1,g_2)$.  Continuing in this fashion, one obtains a complex
\begin{equation}
  \label{eqn:ccc}
\CCC(\CS) = \bigoplus_{r\geq 0} \bigoplus_{X_1,\ldots,X_{r+1}\in \CS} 
\susp{r} \Big(X_1\CS X_2\otimes\cdots \otimes X_{r+1} \CS X_1\Big)
\end{equation}
where the direct sum is over objects $X_1,\ldots,X_{r+1}\in \CS$. We write
$\bbar f_1\bbar f_2\bbar \cdots \bbar f_{r+1}$ with $f_i\in X_i\CS X_{i+1}$
(with indices taken modulo $r+1$) for elementary tensors in the degree $r$ part
of \eqref{eqn:ccc}.

The differential of such an element is
\begin{gather*}
d_{\CCC(\CS)}(\bbar f_1\bbar f_2\bbar \cdots \bbar
f_{r+1}) := (-1)^r \sum_{i=1}^{r+1} (-1)^{|f_1|+\cdots+|f_{i-1}|} 
\bbar f_1\bbar \cdots \bbar d_\CS(f_i) \bbar \cdots \bbar f_{r+1}\\
+ \sum_{i=1}^{r} (-1)^i \bbar f_1\bbar \cdots \bbar f_i f_{i+1} \bbar \cdots \bbar f_{r+1} 
+ (-1)^{|f_1|(r-1+|f_2|+\cdots+|f_{r+1}|)} \bbar f_2\bbar  \cdots \bbar f_{r+1}f_1.
\end{gather*} 

Note that the
differential $d_{\CCC}$ splits as $d_\oplus + \delta$ where $d_\oplus$ is the
direct sum of differentials on the terms in \eqref{eqn:ccc}, shown in the first
line, while $\delta$ is an additional contribution shown in the second line.

The complex just constructed is called the cyclic bar complex of $\CS$.

\begin{definition}\label{def:cyclic bar}
The \emph{cyclic bar complex} of $\CS$ (relative to a unital subcategory
$\IS\subset \CS$) is the dg $\k$-module $\CCC(\CS,\IS)$ obtained as the quotient
of the two-sided bar complex $\BB(\CS,\IS)$ modulo the $\k$-span $[\CS,
\BB(\CS)]$ of elements of the form $f \cdot m - (-1)^{|m||f|}m \cdot f$ for all
$m\in \BB(\CS,\IS)$ and all morphisms $f$ in $\CS$. 

The \emph{Hochschild homology} $\HH_{\subb}(\CS,\IS)$ is defined to be the homology of
$\CCC(\CS,\IS)$.  Hochschild homology is written with homological convention for gradings, and so we write $\HH_k(\CS,
\IS):=H^{-k}(\CCC(\CS,\IS))$ for $k \geq 0$.
\end{definition}

\begin{remark}  Instead of \eqref{eqn:ccc}, the definition describes the cyclic bar complex as:
  \begin{equation}
    \label{eqn:ccct}
  \CCC(\CS) = \bigoplus_{r\geq 0} \bigoplus_{X_0,X_1,\ldots,X_{r+1}\in \CS} 
  \susp{r} \Big(X_0\CS X_1\otimes X_1\CS X_2\otimes\cdots \otimes X_{r+1} \CS X_0\Big) / \sim
  \end{equation}
  where the linear relation $\sim$ is generated by 
  \[(-1)^{|f_0|(r+|f_1|+\cdots+|f_{r+1}|)} f_0 \bbar f_1\bbar f_2\bbar \cdots \bbar f_{r+1}
  \sim \Id \!\bbar f_1\bbar f_2\bbar \cdots \bbar f_{r+1} f_0 \] 
The previously introduced notation for elements of the cyclic bar complex simply drops the leading identity
term after such a rewrite---in the example above: $\bbar f_1\bbar f_2\bbar \cdots \bbar f_{r+1}
f_0 $.
\end{remark}

\begin{remark}
  \label{rem:ccc-full-subcat}
  When $\IS=\IS_{\CS}$ is the subcategory of identity morphisms, then we have
$\CCC(\CS,\IS)=\CCC(\CS)$ as in \eqref{eqn:ccc}. If $\BS\subset \CS$ is a full
subcategory and $\IS=\IS_{\BS}$ is the subcategory of identity morphisms in
$\BS$, then $\CCC(\CS,\IS_{\BS})= \CCC(\BS)$.
\end{remark}

\begin{remark}
If $\CS$ has trivial grading and trivial differential then $\HH_0(\CS)$
coincides with the vertical trace defined earlier.
\end{remark}

An inclusion of subcategories $\IS\subset \JS\subset \CS$ gives a
chain map $\CCC(\CS,\IS)\rightarrow \CCC(\CS,\JS)$.  In particular, if
$\BS\subset \CS$ is a full subcategory then we have a canonical inclusion
$\CCC(\BS)\hookrightarrow \CCC(\CS)$.

The following is our main tool for computing $\HH_{\subb}(\CS)$.

\begin{theorem}
\label{thm:cyclicbarcx-retract}
If $\BS\subset \CS$ generates $\CS$ (see Definition \ref{def:generates}) then
the natural inclusion $\CCC(\BS)\hookrightarrow \CCC(\CS)$ is the section of a
deformation retract. More generally, if $\Gamma$ is a finite poset and
$\{\BS_{\gamma}\}_{\gamma\in \Gamma} $ gives a semi-orthogonal decomposition of
$\CS$ then $\CCC(\CS)$ deformation retracts onto $\bigoplus_{\gamma\in
\Gamma}\CCC(\BS_\gamma)$.
\end{theorem}

\begin{proof}
For the first assertion, recall from Proposition~\ref{prop:barcx-retract}
that $\BB(\CS)$ deformation retracts onto $\BB(\CS,\IS_{\BS})$, 
which induces a deformation retract 
\[
  \CCC(\CS) = \BB(\CS) / [\CS, \BB(\CS)] 
\buildrel\simeq\over\twoheadrightarrow  \BB(\CS,\IS_{\BS})/[\CS, \BB(\CS,\IS_{\BS})]
= \CCC(\BS)
\]
by Remark~\ref{rem:ccc-full-subcat} since $\BS$ is a full subcategory of $\CS$.

To prove the second assertion, we set $\BS=\bigcup_\gamma \BS_\gamma$ and recall
from Theorem \ref{thm:barcx-retract} that $\BB(\CS)$ deformation
retracts onto the semi-orthogonal bar complex $\BB^{\Gamma}(\CS,\IS_{\BS})$ in
the category of $\CS,\CS$-bimodules. There is an induced deformation retract
\[
\CCC(\CS) = \BB(\CS) / [\CS, \BB(\CS)] 
\buildrel\simeq\over\twoheadrightarrow  \BB^{\Gamma}(\CS,\IS_{\BS})/[\CS, \BB^{\Gamma}(\CS,\IS_{\BS})].
\]

The complex $\BB^{\Gamma}(\CS,\IS_{\BS})$ is a one-sided twist of 
\[
\bigoplus_{r\geq 0} \bigoplus_{X_1,\ldots,X_{r+1}\in \BS} 
\susp{r}\Big(\CS X_1\otimes X_1\CS X_2\otimes \cdots\otimes  X_{r+1}\CS\Big)
\]
with $X_1\geq X_2\geq \cdots \geq X_{r+1}=:X_0$. 
By semi-orthogonality, the summands with $X_1>X_{0}$ 
become contractible in the quotient
$\BB^{\Gamma}(\CS,\IS_{\BS})/[\CS, \BB^{\Gamma}(\CS,\IS_{\BS})]$. 
Contracting these induces a deformation retract onto 
the subcomplex where $X_1,\ldots,X_{r},X_0\in \BS_\gamma$ for some $\gamma$.   
For each $\gamma$ the contribution of all such terms is the two-sided bar complex $\BB(\BS_\gamma)$.  
By inspection the differential preserves summands, and we obtain the direct sum
decomposition in the statement.
\end{proof}

The following corollary was proved in \cite{Kuz}.

\begin{corollary}
  Suppose $\{\BS_{\gamma}\subset \CS\}$ defines a $\Gamma$-indexed semi-orthogonal
  decomposition of $\CS$. Then we have
  $$
  \HH_{\subb}(\CC)\cong \bigoplus_{\gamma\in \Gamma} \HH_{\subb}(\BS_{\gamma}).
  $$
  \end{corollary}

\begin{corollary}
  \label{cor:HH-of-complexes}
  If $\CS$ is a dg category, then the
  cyclic bar complex of $\pretr{\CS}$ deformation retracts onto the cyclic bar
  complex of $\CS$, and we have a natural isomorphism $\HH_{\subb}(\CS)\cong \HH_{\subb}(\pretr{\CS})$ 
  induced by the inclusion $\CS\hookrightarrow \pretr{\CS}$.\qed
  \end{corollary}

\begin{example}
If $\AS$ is an additive category (regarded as a dg category trivially) then we
have the natural isomorphism $\HH_{\subb}(\AS)\cong \HH_{\subb}(\Ch^b(\AS))$ induced by the inclusion
$\AS\hookrightarrow \Ch^b(\AS)$.
\end{example}

If $X\in \CS$ is an object, then we write $[X]:=[\Id_X]$ for its class in 
$\HH_0(\CS)\subset \HH_{\subb}(\CS)$. If $X$ is a one-sided twisted complex constructed from objects
$X^i$, it is natural to ask for the relation between the classes $[X]$ and $[X^i]$. 

\begin{lemma}
For objects $X,Y\in \CS$ and $f\colon X \to Y$, we have:
\begin{subequations}
\begin{equation}
[X\oplus Y] = [X]+[Y]
\end{equation}
\begin{equation}
[\susp{1}X] = -[X]
\end{equation}
\begin{equation}
\label{eq:cone}
[\Cone(f:X\rightarrow Y)] = [Y]-[X]
\end{equation}
\end{subequations}
\end{lemma}

\begin{proof}
It is sufficient to prove \eqref{eq:cone}. Let $Z=\Cone(X\xrightarrow{f}Y)$, let
$i_X, i_Y$ denote the inclusions of $X$ and $Y$ into $Z$ and let $\pi_X,\pi_Y$
denote the projections of $Z$ onto $X$ and $Y$. Note that $i_X$ and $\pi_X$ have degrees $\pm 1$, and
\begin{gather*}
d_\CS(i_X)=i_Y f,\quad d_\CS(i_Y)=0,\quad d_\CS(\pi_X)=0,\quad d_\CS(\pi_Y)=-f \pi_X,\\
\pi_X i_X=\Id_X, \quad \pi_Y i_Y=\Id_Y,\quad i_X\pi_X+i_Y\pi_Y=\Id_Z.
\end{gather*}
So we compute
\begin{gather*}
  d_{\CCC}(\bbar i_X\bbar \pi_X ) = - \bbar \!\Id_X - \bbar i_X\pi_X  - \bbar i_Y f\bbar \pi_X,
  \\
  d_{\CCC}(\bbar i_Y\bbar \pi_Y) =  \bbar \! \Id_Y -\bbar i_Y\pi_Y  + \bbar i_Y \bbar f\pi_X,
  \\
  d_\CCC(\bbar i_Y\bbar f\bbar \pi_X) = 0-\bbar i_Y f\bbar \pi_X + \bbar i_Y\bbar f\pi_X,
\end{gather*}
since $\pi_X i_Y = 0$ (each of the above has cohomological degree $-1$ before applying $d_\CCC$).  Therefore
\[
d_\CCC(\bbar i_Y\bbar f\bbar \pi_X - \bbar i_X\bbar \pi_X  - \bbar i_Y\bbar \pi_Y) = \bbar\! \Id_Z - (\bbar \! \Id_Y - \bbar \! \Id_X),
\]
which shows that $[\Id_Z]\simeq [\Id_Y] - [\Id_X]$ in $\HH_0(\CS)$.\
\end{proof}

\begin{corollary}
  Let $\CS$ be a dg category and consider an object $X=\tw_\a(\bigoplus_i X^i)$
  in $\pretr{\CS}$, then under the identification $\HH_{\subb}(\pretr{\CS})\cong
  \HH_{\subb}(\CS)$ the class of $X$ is given by the Euler characteristic
  \[
  [X] = \sum_i (-1)^i [X^i].
  \]
  \end{corollary}

\subsection{Connes differential and HKR isomorphism}
\label{sec: Connes}
The cyclic bar complex of a dg category has a canonical differential $\B$ of cohomological degree $-1$
\cite{MR780077} that we now describe. 
Given 
$\bbar f_1\bbar \cdots\bbar f_{r+1}\in \CCC^{-r}(\CS)$ with $f_1 \in X_1 \CS X_2$,
we define
\begin{align*}
t(\bbar f_1\bbar \cdots\bbar f_{r+1})&=(-1)^{|f_{r+1}|(|f_1|+\cdots+|f_{r}|)} \bbar f_{r+1}\bbar f_1\bbar \cdots\bbar f_{r}\\
s(\bbar f_1\bbar \cdots\bbar f_{r+1})&=\bbar\! \Id_{X_1}\!\bbar f_1\bbar \cdots\bbar f_{r+1}.
\end{align*}
Now we define $N=1+\ldots+t^r$ and $\B=(1-t)sN$. The operator $\B$ is sometimes
called Connes differential.

\begin{example}
For $\bbar f_1\in X_1\CS X_1 \subset \CCC^0(\CS)$ we have 
$$
\B(\bbar f_1)=\bbar\!\Id_{X_1}\!\bbar f_1-\bbar f_1\bbar \!\Id_{X_1}\in\CCC^{-1}(\CS).
$$
For $\bbar f_1\bbar f_2\in \susp{1}\left( X_1 \CS X_2\otimes X_2\CS X_1 \right)\subset \CCC^{-1}(\CS)$ we have
\begin{gather*}
\B(\bbar f_1\bbar f_2)=(1-t)\left(\bbar\!\Id_{X_1}\!\bbar f_1\bbar f_2+ (-1)^{|f_2||f_2|}\bbar \!\Id_{X_2}\!\bbar f_2\bbar f_1\right)
\\
= \bbar\!\Id_{X_1}\!\bbar f_1\bbar f_2
+(-1)^{|f_2||f_1|}\bbar \!\Id_{X_2}\!\bbar f_2\bbar f_1- (-1)^{|f_2||f_1|}\bbar f_2\bbar \!\Id_{X_1}\!\bbar f_1
-\bbar f_1\bbar \!\Id_{X_2}\!\bbar f_2.
\end{gather*}
\end{example}
The following is well known.

\begin{lemma}[\cite{MR780077}]
We have $d_{\CCC}\B+\B d_{\CCC}=0$ and $\B^2=0$. 
\end{lemma}

\begin{theorem}[\cite{MR142598,MR780077}]
Let $A$ be the ring of functions on a smooth affine scheme $X$. Then 
there is an algebra isomorphism
$$
\Omega^{\supb}(X)\simeq \HH_{\subb}(A),
$$ 
which identifies de Rham differential $D$ on the algebra of differential forms
$\Omega^{\supb}(X)$ on the left with the (induced) Connes differential $\B$ on the
right hand side. 
\end{theorem}

\begin{example}
We have $\HH_1(A)=A\otimes A/(ab\otimes c-a\otimes bc+ac\otimes b)$. We can identify this with $\Omega^1(X)$
by sending $a\otimes b\to aD(b)$. 
Indeed, $aD(bc)=abD(c)+acD(b)$. Now $\B(a)=1\otimes a-a\otimes 1\in A\otimes A$ is identified with 
$D(a)-aD(1)=D(a)\in \Omega^1(X)$. 
\end{example}

\begin{example}
Let $R=\C[x_1,\ldots,x_n]$ be the algebra of functions on $\C^n$. Then 
$$\HH_k(R)\simeq \Omega^k(\C^n)\simeq \C[x_1,\ldots,x_n]\otimes \wedge^k[\theta_1,\ldots,\theta_n],$$
where $\theta_i=D(x_i)$. The de Rham differential can be written as 
\begin{equation}
\label{eq: de rham}
D=\sum_i \theta_i\frac{\partial}{\partial x_i}.
\end{equation}
\end{example}

\section{The dg monoidal center and trace}
\label{s:center and trace}

\subsection{Monoidal dg categories and the shuffle product}
\label{sec:monoidal}
Note that any (dg) algebra $A$ can be regarded as a dg category with one object.
By the same token, a commutative (dg) algebra $A$ can be viewed as a
\emph{monoidal} dg category with one object.

A monoidal dg category is a dg category $\CS$ equipped with an object $\one\in\CS$, 
a functor $\star: \CS\otimes_\k \CS\rightarrow \CS$, 
and closed degree zero natural isomorphisms (\emph{associator} 
and \emph{unitors})
\[
(X\star Y)\star Z\cong X\star (Y\star Z),\qquad \one\star X\cong X \cong X\star \one
\]
satisfying the usual coherence relations for monoidal categories (on the nose;
not up to homotopy). The associator and unitor isomorphisms will usually be
suppressed from the notation, and we will refer to $(\CS,\star,\one)$ as a dg
monoidal category. If the associators and unitors are in fact identity
morphisms, then $(\CS,\star,\one)$ is said to be \emph{strictly monoidal}.  The
usual Eckmann-Hilton argument shows that $\End_\CS(\one)$ is always a
commutative dg algebra.

We say that a dg monoidal category {\em has duals} if every object $X$ has a left dual
$X^*$ and a right dual ${}^*\!X$ with evaluation and coevaluation maps
\[
\ev_X\colon X^*\star X \to \one,\quad \coev_X\colon \one \to X\star X^*
,\quad
\ev'_X\colon X\star {}^*\!X \to \one,\quad \coev'_X\colon \one \to {}^*\!X\star X 
\]
satisfying the usual \emph{string-straightening} axioms. Note that the existence
of duals is a property of $\CS$.

For a dg monoidal category $\CS$ with duals, the operations $(-)^*$ and
${}^*\!(-)$ extend to (contravariant) monoidal dg functors. A \emph{pivotal
structure} on $\CS$ is a monoidal natural isomorphism between $(-)^*$ and
${}^*\!(-)$ (equivalently, a monoidal natural isomorphism between $\Id_\CS$ and $(-)^{**}$). 

We say a dg monoidal category $\CS$ is \emph{strictly pivotal} if it is strictly
monoidal and we have $X^*={}^*\!X$ for every object $X$, with the identity
natural transformation as pivotal structure.

A monoidal structure on $\CS$ endows the two-sided bar complex $\BB(\CS,\IS)$
with the additional structure of an algebra, via the so-called \emph{shuffle
product} (e.g. \cite{MR780077,EZ}) which we recall below.

Suppose $\CS$ is a dg monoidal category.  Let
\[
\underline{f} = f_0\bbar \cdots  \bbar f_{r+1} \ \ \in \ \  
\susp{r}\left( X_0\CS X_1\otimes \cdots\otimes X_{r+1}\CS X_{r+2}\right)
\subset \BB_r(\CS,\IS)
\]
and
\[
\underline{g} = g_0\bbar \cdots \bbar g_{s+1} \ \ \in \ \  
\susp{s} \left( Y_0\CS Y_1\otimes \cdots\otimes Y_{s+1}\CS Y_{s+2}\right)
\subset \BB_s(\CS,\IS)
\]
be two elements in $\BB(\CS,\IS)$. 
To define the product $\underline{f} \ast \underline{g}$ we set 
\[e_i=
\begin{cases}
f_i\star \Id & 1\leq i \leq r\\
\Id \star g_{i-r} & r < i \leq r+s
\end{cases}
  \]
where we keep an open mind about what objects the identity morphisms are associated with. 
Now we set
\[
\underline{f} \ast \underline{g} := 
(-1)^{|\underline{f}| s}\sum_{\pi\in  S_{(r,s)}} (-1)^{w\sigma(\pi,\underline{f},\underline{g})} (f_0\star g_0) \bbar e_{\pi(1)} 
\bbar \cdots \bbar e_{\pi(r+s)} \bbar (f_{r+1}\star g_{s+1})
\]
where
$S_{(r,s)}\subset S_{r+s}$ denotes the set of shuffle permutations, and
$w\sigma(\pi,\underline{f},\underline{g})$ denotes the \emph{weighted sign} of
the permutation $\pi$, to which a transposition of $f_i\star \Id$ and $\Id\star
g_j$ contributes $(-1)^{|f_i||g_j|}$. For the $e_{\pi(i)}$ in the summands of
this formula, we implicitly choose the identity morphism factor which makes the
sequence of morphisms in the summand composable. Note that different summands
require different choices, although this is suppressed in the notation. 

The shuffle product together with the coproduct $\Delta$ give the two-sided bar
complex $\BB(\CS,\IS)$ the structure of a bialgebra.

\subsection{The quadmodule associated to a dg monoidal category}
\label{sec:diag}
Let $\CS$ now be a dg monoidal category.  Fix objects $X_1,X_2,Y_1,Y_2\in \CS$, and
let $f\in \Hom_\CS(X_1\star X_2, Y_1\star Y_2)$ be given.  We also consider
objects $X_i',Y_i'\in \CS$ ($i=1,2$) and morphisms
\[
a_i\in \Hom_\CS(Y_i,Y_i'),\qquad\qquad b_i \in \Hom_\CS(X_i',X_i).
\]
Then we define the following operations:
\[
a_1\circ_1 f := (a_1\star \Id_{Y_2})\circ f,\qquad\qquad a_2\circ_2 f :=(\Id_{Y_1}\star a_2)\circ f,
\]
\[
f\circ_1 b_1 := f\circ (b_1\star \Id_{X_2}),\qquad\qquad f\circ_2 b_2:=f\circ (\Id_{X_1}\star b_2).
\]
These operations give a $(\CS\otimes_\k \CS), (\CS\otimes_\k \CS)$-bimodule structure on 
\[
\XS:= \bigoplus_{X_1,X_2,Y_1,Y_2} (Y_1\star Y_2)\CS (X_1\star X_2).
\]
We may regard $\XS$ just defined as a \emph{quadmodule} over $\CS$.   By
combining the monoidal structure and composition of morphisms in $\CS$, we have
morphisms
\[
\mu_{\nwarrow},\mu_{\nearrow}:\XS\otimes_\k \XS\rightarrow \XS   
\]
defined by
\[
\mu_{\nwarrow}(f,g):=(f\star \Id)\circ (\Id\star g)
\]
and
\[
\mu_{\nearrow}(f,g):=(\Id\star f)\circ (g\star \Id)
\]
whenever these compositions make sense.  These operations interact with the
quadmodule structure according to
\begin{align*}
\mu_{\nwarrow}(f\circ_2 a , g) &= \mu_{\nwarrow}(f, a \circ_1 g)\\
\mu_{\nwarrow}(a\circ_1 f , g\circ_2 b) &= a\circ_1 \mu_{\nwarrow}(f,  g)\circ_2 b\\
\mu_{\nwarrow}(a\circ_2 f , a'\circ_2 g) &= (a\star a')\circ_2\mu_{\nwarrow}(f,  g)\\
\mu_{\nwarrow}(f\circ_1 b ,  g\circ_1 b') &=\mu_{\nwarrow}(f,  g)\circ_1 (b\star b'),
\end{align*}
with similar identities involving $\mu_{\nearrow}$ (swapping the roles of
$\circ_1$ and $\circ_2$).

From $\XS$ we obtain a $\CS,\CS$-bimodule by forgetting the ``northeast'' and
``southwest'' actions of $\CS$.  Precisely, $\XS_{12}$ equals $\XS$, but with
$\CS,\CS$-bimodule structure 
\[
a\otimes f\otimes b \ \ \mapsto  \ \ a \circ_1 f\circ_2 b,\qquad\qquad a,b\in \CS,\qquad f\in \XS_{12}.
\]
One may define a bimodule $\XS_{21}$ in a similar fashion, but we will not need
it.
\begin{remark}
  \label{rem:quadmod-nw-algebra}
The map $\mu_{\nwarrow}$ makes $\XS_{12}$ into an algebra object in
$\CS,\CS$-bimodules, with unit $\CS\rightarrow \XS_{12}$ given by the bimodule
map sending $\Id_X\in \CS$ to the canonical isomorphism $\one\star X\rightarrow
X\star \one$, regarded as an element of $X\XS_{12} X$.
\end{remark}
\begin{remark}
The map $\mu_{\nearrow}$ also defines an associative multiplication on
$\XS_{12}$ which interacts with the bimodule structure in a nonstandard
way:
\begin{align*}
  \mu_{\nearrow}(f\circ_1 a , g) &= \mu_{\nearrow}(f, a \circ_2 g)\\
  \mu_{\nearrow}(a\circ_2 f , g\circ_1 b) &= a\circ_2 \mu_{\nearrow}(f,  g)\circ_1 b\\
  \mu_{\nearrow}(a\circ_1 f , a'\circ_1 g) &= (-1)^{|a'|(|a|+|f|)} (a'\star a)\circ_1\mu_{\nearrow}(f,  g)\\
  \mu_{\nearrow}(f\circ_2 b ,  g\circ_2 b') &= (-1)^{|b|(|b'|+|g|)}\mu_{\nearrow}(f,  g)\circ_2 (b'\star b),
  \end{align*}
\end{remark}

By fixing $Y_2=Z'$ and $X_1=Z$, we obtain sub-bimodules of $\XS_{12}$ of the
form
\[
\XS_{12}(Z',Z) \ := \ \bigoplus_{X,Y\in \CS} (Y\star Z')\CS (Z\star X) \ \subset \ \XS_{12}
\]

Note that $\mu_{\nwarrow}$ restricts to morphisms of $\CS,\CS$-bimodules
\[
\mu_{\nwarrow}  \ : \ \XS_{12}(Z',Z) \otimes_\CS \XS_{12}(U',U)\rightarrow \XS_{12}(Z'\star U',Z\star U),
\]
while $\mu_{\nearrow}$ restricts to morphisms
\[
\mu_{\nearrow}  \ : \ \XS_{12}(Z'',Z') \otimes_\k \XS_{12}(Z',Z)\rightarrow \XS_{12}(Z'', Z)
\]
In particular the bimodule $\XS_{12}(Z):=\XS_{12}(Z,Z)$ inherits an associative
multiplication (which respects the dg $\k$-module structure and is compatible
with the $\CS,\CS$-bimodule structure).

\subsection{The dg monoidal centralizer}

Let $\CS$ be a dg monoidal category, and fix an object $Z\in \CS$.  
We would like to discuss what it means for $Z$ to be central in $\CS$.  Just as
in the usual Drinfeld center, this is not a \emph{property} enjoyed by $Z$, but
rather additional \emph{structure} which must be provided.

Actually, we will consider the slightly broader problem of defining what it
means for $Z$ to centralize a full dg monoidal subcategory $\MS\subset \CS$. 
For $Z$ to
centralize $\MS$ (up to homotopy) requires the following data:
\begin{enumerate}
\item for each object $X\in \MS$, a degree zero closed morphism, called \emph{half-braiding}, $\tau_X\in (X
\star Z)\CS(Z \star X)$.
\item for each closed morphism $f\in Y\MS X$ a homotopy $h_f\in (X \star Z)\CS(Z
\star X)$ with $
d_\CS(h_f) = f\circ_1\circ \tau_X - \tau_Y\circ_2 f
$
\item  certain higher homotopies.
\end{enumerate}
To get a feeling for the sort of higher homotopies required, observe that for
each pair of closed morphisms $f_0 \in X_0\MS X_1$ and $f_1\in X_1\MS X_2$ we
have \emph{two ways} of commuting $f_0\circ f_1$ past $Z$.  First, we have the
homotopy $h_{f_0\circ f_1}$.  But we also have $h_{f_0}\circ_2 f_1 +
(-1)^{|f_0|} f_0\circ_1 \circ h_{f_1}$.  We should require the difference of
these two homotopies (which is a closed morphism of degree $-1$) to be
null-homotopic.   The various higher homotopies required are, in fact, already
organized for us in the form of the two-sided bar complex.

\begin{definition}\label{def:dg center}
Let $\CS$ be a dg monoidal category and $\MS\subset \CS$ a subcategory. The
\emph{dg monoidal centralizer of $\MS$ in $\CS$} is the dg category
$\ZC^{\text{dg}}_\CS(\MS)$ whose objects are pairs $(Z,\tau)$ where $Z\in \CS$
and $\tau\colon \BB(\MS)\rightarrow \XS_{12}(Z)$ is a map of $\CS,\CS$-bimodules
as well as a map of dg algebras.  The complex of morphisms in
$\ZC^{\text{dg}}_\CS(\MS)$ from $(Z,\tau)$ to $(Z',\tau')$ is the subcomplex of
$\Hom_\CS(Z,Z')$ consisting of those morphisms which commute strictly with the
\emph{structure maps} in the images of $\BB(\MS)$. I.e. those $z\in
\Hom_\CS(Z,Z')$ such that for $f \in Y \BB(\MS) X$ we have 
\[
z\circ_2 \tau(f) = \tau'(f)\circ_1 z
  \]
The \emph{dg Drinfeld center} is defined to be
$\ZC^{\text{dg}}(\CS):=\ZC^{\text{dg}}_{\CS}(\CS)$.
\end{definition}

\begin{remark}
  \label{rem:caveat}
It is perhaps better to weaken the condition on morphisms in
$\ZC^{\text{dg}}(\CS)$ to strictly commute with the structure maps from $\BB(\CS)$.
This would bring us into the world of $\Ai$-categories, which we choose to avoid
for the moment. In the rest of this paper, we only consider objects of the dg
Drinfeld center and not morphisms.  
\end{remark}

The dg Drinfeld center $\ZC^{\text{dg}}(\CS)$ has a tensor product defined by
\[
(Z,\tau)\star (Z',\tau') := (Z\star Z', \tau''),
\]
where $\tau''$ is the composition of maps
\[
\BB(\CS)\rightarrow \BB(\CS)\otimes_{\CS} \BB(\CS) \rightarrow \XS_{12}(Z)\otimes_{\CS} 
\XS_{12}(Z') \rightarrow \XS_{12}(Z\star Z').
\]
Here, the first map is the coproduct on the two-sided bar complex, and the second map sends
\[
f\otimes f' \mapsto (f\star \Id_{Z'})\star (\Id_Z\star f') \in X'' \XS_{12}(Z\star Z')X
\]
for all $f'\in X'\XS_{12}(Z')X $ and all $f\in X'' \XS_{12}(Z)X'$.

  The following is immediate.
\begin{proposition}
  There is a natural forgetful functor $\ZC^{\text{dg}}(\CS)\to \CS$. It is monoidal.
\end{proposition}

\begin{remark}
It is well known that the Drinfeld center of a monoidal category is braided. It
is natural to expect that dg Drinfeld center is braided in a dg sense, that is,
the braiding is natural up to homotopy. To define such a structure, it seems
likely that passage to the world of $\Ai$ (braided monoidal) categories and
functors is unavoidable.  We save such explorations for future work.
\end{remark}

\begin{remark} If $\CS$ is a monoidal category, considered as a dg category with
 trivial differential, then the Drinfeld center $\ZC^{0}(\CS)$ is a monoidal
 subcategory of $\ZC^{\text{dg}}(\CS)$.
\end{remark}

In the following remarks we spell out the meaning of some of the structure maps 
that are part of the data of a central object $(Z, \tau)$.

\begin{remark}
For each $X\in \CS$ the two-sided bar complex has a distinguished degree zero closed
element $\Id_X\!\bbar\! \Id_X$.  The image of this element in $\XS_{12}(Z)$
under $\tau$ will be denoted $\tau_X$ and is called the half-braiding of $X$ with $Z$.  
Then each degree zero element
$a \bbar b$ gets sent to $a\circ_1\tau_X \circ_2 b$ where $X$ is the codomain of
$b$ (same as the domain of $a$), since $\tau$ commutes with the
bimodule structures on $\BB(\CS)$ and $\XS_{12}(Z)$.
\end{remark}

\begin{remark}
In $\BB(\CS)$ we have $(\Id_X\!\bbar\! \Id_X)*(\Id_Y\!\bbar\! \Id_Y)=\Id_{X\star
Y}\!\bbar\! \Id_{X\star Y}$. By definition of the center, $\tau$ is a map of dg
algebras, which implies 
$$
\tau_{X\star Y}=(\Id_{X}\star \tau_Y)\circ (\tau_X\star \Id_{Y}) 
$$
and one recovers the familiar compatibility between half-braiding morphisms and
the monoidal structure in $\CS$.
\end{remark}

\begin{remark}
For each morphism $f\in Y\CS X$ of degree $l$ the two-sided bar complex has a degree $l-1$
element of the form $\Id_Y\!\bbar f\bbar\! \Id_X$, whose image under $\tau$ we
denote by $h_f$.  If $f$ is closed then we have
\[
d(h_f) = d(\tau(\Id_Y\!\bbar f\bbar\! \Id_X))  = \tau(f\bbar\! \Id_X - \Id_Y\!\bbar f) = f\circ_1 \tau_X - \tau_Y\circ_2 f,
\]
since $\tau$ commutes with the differentials in $\BB(\CS)$ and
$\XS_{12}(Z)$. This means $h_f$ is a homotopy for commuting $f$ through
the half-braiding with $Z$. 

Similarly for each sequence of composable (closed) morphisms $f,g$ we have an element
$h_{f,g} = \tau(\bbar f\bbar g\bbar)$ (dropping explicit occurrences of identity maps) satisfying
\begin{eqnarray*}
d(h_{f,g}) 
&=& d(\tau(\bbar f\bbar g\bbar)) \\ 
&=& \tau(d(\bbar f\bbar g\bbar)) \\ 
&=& \tau(f\bbar g\bbar - \bbar f\circ g \bbar + \bbar f \bbar g) \\ 
&=& (-1)^{|f|}f\circ_1 h_g - h_{f\circ g} + h_f\circ_2 g
\end{eqnarray*}
and so on. In summary, we see that $\tau$ gives the data of a half-braiding with $Z$ 
that is natural up to coherent homotopy. 
\end{remark}

\begin{remark}
  If $\CS$ is an arbitrary dg category, there is a functor from $Z^{dg}(\CS)$
  to the usual Drinfeld center of the homotopy category $Z^0(H^0(\CS))$. 
  It sends a central object $(Z,\tau)$ to $(Z,\tau_X)$ and forgets all higher homotopies.
  \end{remark}

\begin{lemma} 
  If $(Z,\tau)$ is an object of $\ZC^{\text{dg}}(\CS)$ and $X$ an object of $\CS$
  that has a right dual, then $\tau_X:X\star Z\rightarrow X\star Z$ is a homotopy
  equivalence.
  \end{lemma}
  \begin{proof} If $^*\!{X}$ is a right dual of $X$ with \[
    \ev' \colon X\star {}^*\!X \to \one,\quad \coev'\colon \one \to {}^*\!{X}\star X
    \] then we set
    \[
    \tau^{-1}_X:= (\ev'\star \Id_{Z\star X})
    \circ (\Id_X\star\;\tau_{^*\!{X}} \star \Id_{{}^*\!X})
    \circ (\Id_{X\star Z}\star \;\coev')
    \]
  and one can check that $\tau_X \tau^{-1}_X\simeq \Id_{X\star Z}$ and $\tau^{-1}_X \tau_X \simeq \Id_{Z\star X}$.
  \end{proof}

\begin{remark}
\label{rem:z commutes past complexes}
The dg monoidal center $\ZC^{\text{dg}}(\CS)$  embeds as a monoidal dg subcategory of
$\ZC^{\text{dg}}(\pretr{\CS})$. To see this, let $(Z,\tau)$ be an object of
$\ZC^{\text{dg}}(\CS)$. The half-braiding morphism $\tau$ extends to $\adds\CS$
in a trivial fashion, so $(Z,\tau)$ can be thought of as an object of
$\ZC^{\text{dg}}(\adds\CS)$. Next, given a morphism $f\colon X\to Y$, the
half-braiding of the cone $\Cone(f):=[X\xrightarrow{f} Y]$ past $Z$ is defined
by the following morphism:
\[
\tau_{\Cone(f)}:
\begin{tikzcd}
\big( Z\star X\arrow{r}{\Id_Z\star f} \arrow{d}{\tau_X} \arrow{dr}{h_f}& Z\star Y\big) \arrow{d}{\tau_Y}\\
\big( X\star Z\arrow{r}{f\star \Id_Z} & Y\star Z\big)\\
\end{tikzcd}
\]
Furthermore, the (higher) homotopies for half-braidings for cones are
analogously determined by $\tau$. Iterating this construction, we see that
$(Z,\tau)$ is derived central for one-sided twisted complexes, i.e. it
represents an object in $\ZC^{\text{dg}}(\pretr{\CS})$.
\end{remark}

\begin{remark}
Given an object $(Z,\tau)$ in $\ZC^{\text{dg}}(\CS)$, it is instructive to
compute the half-braiding of twisted complexes past $Z$
explicitly. Thus, suppose we have a twisted complex $\tw_\a(X)$ with $X\in \CS$
(or $X\in \Sigma \CS$).

Note that $Z\star \tw_\a(X) = \tw_{\Id_Z\star \a} (Z\star X)$ and
$\tw_\a(X)\star Z = \tw_{\a\star \Id_Z}(X\star Z)$.  The half-braiding morphism
$Z\star \tw_\a(X)\rightarrow \tw_\a(X)\star Z$ can then be constructed using
standard homological perturbation theory as:
\begin{align}\label{eq:tau of cx}
\tau_{\tw_\a(X)}&\colon Z\otimes \tw_\a(X) \rightarrow \tw_\a(X)\otimes Z,\\
\nonumber \tau_{\tw_\a(X)} \ &:= \ \sum_{r\geq 0} (-1)^{{r+1 \choose 2}}
\tau(\Id_X\bbar \underbrace{\a\bbar\ldots\bbar\a}_r\bbar\Id_X).
\end{align}
(This is well-defined if $X$ has the structure of a one-sided twisted complex.)  To see that this is a closed
degree zero morphism in $\Tw(\Sigma \CS)$, we must show that
\[
d(\tau_{(X,\a)}) 
+ (\a \otimes \Id_Z) \circ \tau_{(X,\a)} 
- \tau_{\tw_\a(X)}\circ (\Id_Z\otimes \a) = 0.
\]

We compute:
\begin{eqnarray*}
d(\tau_{\tw_\a(X)}) 
&=& \sum_{r\geq 0}(-1)^{{r+1 \choose 2}} d(\tau(\Id_X\!\bbar \underbrace{\a\bbar\ldots\bbar\a}_r\bbar\!\Id_X))\\
&=&\sum_{r\geq i\geq 1} (-1)^{{r+1 \choose 2}+r+i-1}\Id_X\!\bbar \underbrace{\a\bbar\ldots\bbar\a}_{i-1}\bbar d(\a)\bbar  
\underbrace{\a\bbar\ldots\bbar\a}_{r-i}\bbar\!\Id_X\\
&+& \sum_{r\geq 1} (-1)^{{r+1 \choose 2}+r-1} \a\circ \tau(\Id_X\!\bbar
\underbrace{\a\bbar\ldots\bbar\a}_{r-1}\bbar\!\Id_X)\\
&+&\sum_{r> i\geq 1} (-1)^{{r+1 \choose 2}+i}\Id_X\!\bbar \underbrace{\a\bbar\ldots\bbar\a}_{i-1}\bbar \a^2\bbar  
\underbrace{\a\bbar\ldots\bbar\a}_{r-i-1}\bbar\!\Id_X\\
&+& \sum_{r\geq 1}  (-1)^{{r+1 \choose 2}+r}\tau(\Id_X\!\bbar
\underbrace{\a\bbar\ldots\bbar\a}_{r-1}\bbar\!\Id_X)\circ \a\\
&=& - \a\circ \tau_{\tw_\a(X)} + (-1)^{|\tau_{\tw_\a(X)}|}  \ \tau_{\tw_\a(X)}\circ \a.
\end{eqnarray*}

Since ${r+1 \choose 2}+r-1 \equiv {r+2 \choose 2} \; \textrm{mod} \;2$ and
$d(\a)+\a^2=0$ the terms in the
first and third line cancel and we get 
$
d(\tau_{(X,\a)}) =
-(\a \otimes \Id_Z) \circ \tau_{(X,\a)} 
+ \tau_{(X,\a)} \circ (\a \otimes \Id_Z)
$
as desired. 
\end{remark}

\begin{remark} 
It an exercise to check that the braiding morphisms which commute $Z$ past one-sided twisted complexes satisfy the required compatibility with the monoidal structure.
\end{remark}
\subsection{The dg monoidal trace}

Let $\CS$ be a dg monoidal category (which we will soon assume to be strictly
monoidal). First, recall the $\CS,\CS$-bimodule $\XS_{12}$ with
\[
\XS_{12}\ := \ \bigoplus_{X,Y,Z,Z'\in \CS}\Hom(X\star Z, Z'\star Y) 
\ \ \ \ \ \  \Big(= (Z'\star Y)\CS(X\star Z)\Big),
\]
with bimodule structure defined by
\[
a\cdot f\cdot b := a \circ_1 f\circ_2 b\qquad\qquad \forall  
\ \ \ a,b \in \CS,\qquad f\in \XS_{12}.
\]

\begin{definition}
We define a dg category $\Tr(\CS)$, the \emph{(dg horizontal) trace of $\CS$},
as follows.    First, define the $\CS,\CS$ bimodule
\[
\BB(\CS)\otimes_\CS \XS_{12}.
\]
Elements of this bimodule are linear combinations of symbols of the form
$c_0\bbar c_1\bbar\cdots \bbar c_r \bbar (c_{r+1}\circ_1 f)$ where
$c_0,\ldots,c_{r+1}\in \CS$ is a sequence of composable morphisms and $f\in
\XS_{12}$.  Then we identify the left and right $\CS$-actions by forming the
quotient
\[
\CCC(\CS, \XS_{12}) \ := \ (\BB(\CS)\otimes_\CS \XS_{12})/\sim
\]
with respect to the relations of the form
\[
(-1)^{s}c_0\bbar c_1\bbar\cdots \bbar c_r \bbar  f 
\  \ \simeq \ \ \Id\bbar c_1\bbar\cdots\bbar c_r \bbar (f\circ_2 c_0),
\]
where the sign is determined by $s=|c_0|(|c_1|+\cdots+|c_r|+|f|)+|c_0|r$.  As usual we will typically
drop leftmost identity map from the notation, writing
\[
\bbar c_1\bbar\cdots\bbar c_r \bbar f \ = \ \Id\bbar c_1\bbar\cdots\bbar c_r \bbar f 
\]
Note that $\CCC(\CS, \XS_{12})$ is the Hochschild chain complex of $\CS$ with
coefficients in the bimodule $\XS_{12}$.  Recall that each pair of objects
$X,X'\in \CS$ determines a subbimodule
\[
\XS_{12}(X',X)\ := \ \bigoplus_{Y,Y'} \Hom_\CS(X\star Y, Y'\star X')  \ \subset \XS_{12}
\]
and so we have subcomplexes
\[
\CCC(\CS, \XS_{12}(X',X)) \subset \CCC(\CS, \XS_{12}).
\]
Now, define a dg category $\Tr(\CS)$ as follows.  Objects of $\Tr(\CS)$ are the
same as objects of $\CS$, though to avoid confusion we will write $\Tr(X)$ for
$X\in \CS$ regarded as an object of $\Tr(\CS)$.   The complex of morphisms is
given by
\[
\Hom_{\Tr(\CS)}(\Tr(X), \Tr(X')) \ := \ \CCC(\CS, \XS_{12}(X',X)),
\]
with composition induced induced by the shuffle product on the two-sided bar complex and
composition in $\CS$.
\end{definition}

Explicitly, this means that $\Hom_{\Tr(\CS)}(\Tr(X),\Tr(X'))$ has basis given by
formal symbols of the form
\[
(\underline{c}, f) := \bbar c_1\bbar \cdots \bbar c_{r} \bbar f,
\]
where $r\geq 0$, $c_i\in Y_i \CS Y_{i+1}$ and $f\in (Y_{r+1}\star X')\CS(X\star
Y_1)$ for objects $Y_1,\ldots,Y_{r+1}\in \CS$, and we abbreviate by writing
$\underline{c}=(c_1,\ldots,c_r)$.   When $r=0$, the sequence $\underline{c}$ is
empty, and we will use the notation $(\emptyset,f)$ and $\bbar f$
interchangeably.

 We also allow formal symbols of the form $c_0\bbar c_1\bbar \cdots \bbar c_{r}
 \bbar f$ but modulo the relations imposed on Hochschild chains, such an
 expression equals $\pm \bbar c_1\bbar \cdots \bbar c_{r} \bbar (f\circ_2 c_0)$.
 We picture these symbols as follows.

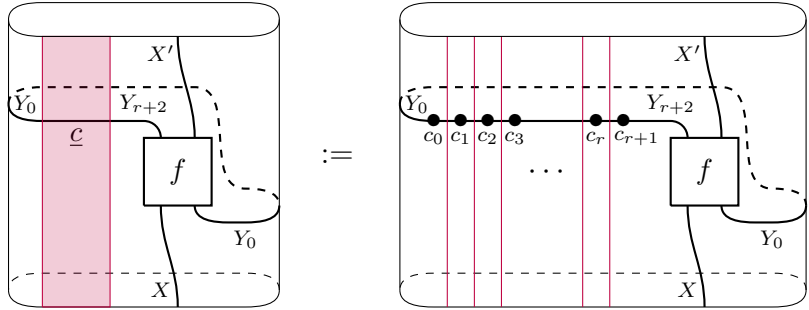
\[\htrfigure{.9}\]

The cohomological grading is
\[
\deg(\underline{c}, f) = -r + |f|
\]
We say that $(\underline{c},f)$ has \emph{bar degree} $r$.

The differential is the usual bar differential (an alternating sum of ways of
deleting bars) plus the terms involving the differentials of the individual
components: 

\begin{eqnarray*}
d_{\Tr(\CS)}(\bbar c_1\bbar \cdots \bbar c_r\bbar f) 
& = & c_1\bbar c_2\bbar \cdots \bbar c_r\bbar f \\
&+& \sum_{i=1}^{r-1}(-1)^i \bbar c_1\bbar \cdots \bbar c_i\circ c_{i+1}\bbar \cdots \bbar c_{r}\bbar f \\
& + &(-1)^r \bbar c_1\bbar \cdots \bbar c_{r-1}\bbar (c_r\circ_1 f)\\
&+&\sum_{i=1}^{r} (-1)^{s(i)} \bbar c_1\bbar \cdots\bbar  d(c_i)\bbar \cdots\bbar c_{r}\bbar f\\
 &+ &(-1)^{s(r+1)}  \bbar c_1\bbar \cdots\bbar c_{r}\bbar d(f).
\end{eqnarray*}
where $s(i) = |c_0|+\cdots+|c_{i-1}|$.  When $r=0$ the above generates to
$d(\bbar f) = \bbar d(f)$.

The composition of morphisms is defined by
\[
  (\underline{c}, f)\circ (\underline{d},g) := (\underline{d}\ast \underline{c}, \mu_{\nearrow}(f,g))
\]
where $\ast$ denotes the shuffle product from Section~\ref{sec:monoidal}, and
$\mu_{\nearrow}(f,g)$ is as defined in Section~\ref{sec:diag}. The composition
can be pictured as follows.
\[
  \htrcompfigure{.9}
\]

\begin{remark}
The identity endomorphism of $\Tr(X)$ in $\Tr(\CS)$ is given by
$(\emptyset,\phi)$ where $\phi$ is the canonical isomorphism $X\star
\one\rightarrow \one\star X$.  Then $(\emptyset,\phi)\circ (\emptyset, \phi)$ in
$\Tr(\CS)$ is by definition $(\emptyset, \psi)$ where $\psi$ is the canonical
isomorphism $X\star(\one\star\one)\rightarrow (\one\star \one)\star X$.  This
$\psi$ is homotopic, but not \emph{equal} to $\phi$.  Thus, strictly speaking
$(\emptyset,\phi)$ only acts as the identity of $\Tr(X)$ up to homotopy!  There
is a similar problem concerning the associativity of composition in $\Tr(\CS)$.
The essential issue is that the two-sided bar complex of $\CS$ is too large.  These
annoyances do not arise when $\CS$ is strict monoidal, and so this will be
assumed \emph{en force} in the sequel.  Without this assumption, $\Tr(\CS)$
should be regarded as an $\Ai$-category, not a dg category.
\end{remark}

We make some observations about $\Tr(\CS)$ below.  First, note that there is a dg functor 
$\hTr\colon \CS\rightarrow \Tr(\CS)$ sending $X\mapsto \Tr(X)$ and $f\mapsto \bbar f$.

\begin{remark} If $\AS$ is a an ordinary monoidal category (regarded as a dg
  category with trivial grading and differential) then $H^0(\Tr(\AS))$ is
  isomorphic to the usual \emph{horizontal trace}  $\Trz(\AS)$, as defined in
  \cite[Section 2.4]{MR3606445}. In this sense, $\Tr(\AS)$ is a derived version
  of $\Trz(\AS)$, much as Hochschild homology of an algebra $A$ is a derived
  version of its trace $A/[A,A]$. 
\end{remark}

\begin{lemma}
\label{lem: derived to underived trace}
If $\AS$ is a $\k$-linear monoidal category, considered as a dg category with
trivial grading and differential, then there is a functor $\hTr(\AS)\to
\Trz(\AS)$ defined on objects by $\Tr(X)\mapsto \Trz(X)$ and on morphisms by
$\bbar f\mapsto f$ and $(\underline{c}, f)\mapsto 0$ for all sequences of
$\underline{c}$ of length $r\geq 1$.
\end{lemma}

The well-known relationship between vertical and horizontal traces (turn head
by 90 degrees) transfers to the derived setting as follows. The following is
clear from the definitions.

\begin{proposition}
  \label{prop:verthor}
If $\CS$ is a dg monoidal category then
$$
\CCC_*(\CS)= \End_{\hTr(\CS)}(\hTr(\one)).
$$
as dg algebras.
\end{proposition}
Strictly speaking, $\End_{\hTr(\CS)}(\hTr(\one))$ is given by Hochschild chains
of $\CS$ with coefficients in the bimodule
\[
\bigoplus_{X,Y\in \CS} \Hom(\one\star X, Y\star \one),
\]
while $\CCC_*(\CS)$ is given by Hochschild cochains of $\CS$ with coefficients
in $\CS$.  If $\CS$ is strict monoidal, then these two bimodules are equal;
otherwise they are isomorphic via the unitor maps.

The following says that the natural functor $\CS\rightarrow \hTr(\CS)$ is
satisfies a categorical ``trace-like'' property, provided that $\CS$ has duals.

\begin{lemma}
\label{lem:trace relation}
Given two objects $X,Y$ in $\CS$, we define the \emph{traciator}
\[
w_{X,Y} \ \ : \ \ \Tr(X\star Y)\rightarrow \Tr(Y\star X)
\] 
to be the degree zero closed morphism associated to the identity map (or
associator in the non-strict monoidal case) $(X\star Y)\star X\rightarrow X\star
(Y\star X)$.  If $X$ has a right dual in $\CS$, then this map is a homotopy
equivalence. 
\end{lemma}

\begin{proof} 
Suppose ${}^\ast X$ is the right dual to
$X$, with structure maps  $\textrm{coev}_X \colon \one \to {}^\ast\!X\star X$
and $\textrm{ev}_X\colon X\star {}^\ast\!X \to \one$.  Define
\[w^{-1}_{X,Y}\colon \Tr(Y\star X)\rightarrow \Tr(X\star Y)\] 
to be the composite of 
\[
\bbar(\textrm{coev}_X  \star \Id_Y \star \Id_X)\colon \hTr(Y\star X) \to \hTr({}^\ast\! X\star X \star
Y\star X),
\]
followed by 
\[
w_{{}^\ast\! X,X\star Y\star X}\colon \hTr({}^\ast\! X\star X \star Y\star X) \to \hTr(X \star Y\star X \star {}^\ast\! X)
\] and finally 
\[
\bbar(\Id_X \star \Id_Y \star \textrm{ev}_X)\colon \hTr(X\star Y \star X \star
{}^\ast\! X) \to \hTr(X \star Y).
\]
One can now check that $\Id_{X\star Y}-w^{-1}_{X,Y}w_{X,Y}$ and $\Id_{Y\star X}-w_{X,Y}w^{-1}_{X,Y}$ are exact.
\end{proof}

These morphisms can be pictures as follows.
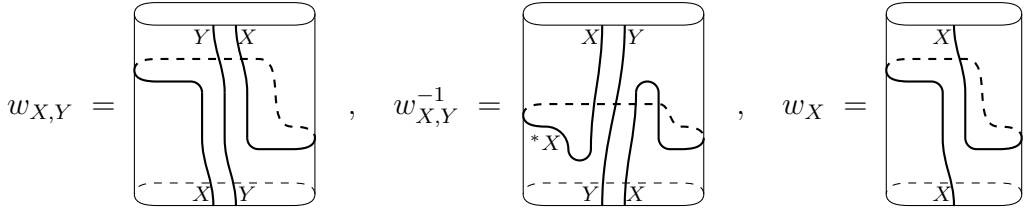
\[\traciator{.6}\]

\begin{remark}
\label{rem:higher traciators}
It is easy to see that the traciator $w_{X,Y}$ is natural in $Y$, but it is
natural in $X$ only up to coherent homotopy. For any closed map $f\colon X\to
X'$ the composition $w_{X',Y}\circ (f\star \Id_Y)\colon \Tr(X\star Y)\to
\Tr(Y\star X')$ is represented by the map $X\star Y\star X'\xrightarrow{f\star
\Id\star \Id} X'\star Y\star X'$ , while the composition $(\Id_Y\star f)\circ
w_{X,Y}$ is represented by the map $X\star Y\star X\xrightarrow{f\star \Id\star
\Id} X\star Y\star X'$. The difference of these two maps is given by the boundary
of the morphism
\begin{equation}
\label{eq:traciator homotopy}
w(f;Y):=\bbar f\bbar\Id_{X\star Y\star X'}
\end{equation}
More generally, for any sequence of composable morphisms $f_i\in X_{i}\CS X_{i+1}$ for $1\le
i\le r$ the derived trace contains \emph{higher traciators}
$$
w(f_{1}\bbar \ldots\bbar f_r;Y):=
\bbar f_{1} \bbar \cdots \bbar f_{r} \bbar \Id_{X_{r+1}\star Y\star X_{1}}
\in \Hom_{\hTr}(\Tr(X_{r+1}\star Y),\Tr(Y\star X_{1}))
$$
which provide the (higher) naturality data for the traciator.
\end{remark}

\begin{remark} 
\label{rem:rotator}
For an object $X$ in $\CS$ we define the \emph{rotator}
$w_X\colon \Tr(X)\rightarrow \Tr(X)$ to be the endomorphism given by $\bbar
\!\Id_{X\star X}$. Note that this agrees with the traciator $w_{X,\one}$ as
defined in the proof of Lemma~\ref{lem:trace relation} if $\CS$ is strict
monoidal. If $X$ has a right dual then $w_X$ is invertible up to homotopy.

In this way, the assignment $X\mapsto w_X$ is a degree zero endomorphism of the
canonical dg functor $\CS \to \hTr(\CS)$, natural up to coherent homotopy in the sense
of Remark \ref{rem:higher traciators}.
\end{remark}

\begin{remark}
  Similar to Remark \ref{rem:z commutes past complexes}, the higher traciators
  (from \eqref{eq:traciator homotopy}) for $\Tr(\CS)$ already carry enough
  information to determine traciators in $\Tr(\pretr{\CS})$. For example, given a
  morphism $f \colon X\to X'$ and its cone $\Cone(f):= [X\to X']$, the traciator
  $w_{\Cone(f),Y}\colon \Tr(\Cone(f)\star Y)\to \Tr(Y \star \Cone(f))$ is
  represented by the morphism
  \begin{equation}
  \label{eq:traciator two term}
  w_{\Cone(f),Y}=
  \begin{tikzcd}
  \big( \Tr(X\star Y) \arrow{r}{\Tr(f\star \Id_Y)} \arrow{d}{w_{X,Y}} \arrow{dr}{w(f;Y)}& \Tr(X'\star Y)\big) \arrow{d}{w_{X',Y}}\\
  \big( \Tr(Y\star X) \arrow{r}{\Tr(\Id_Y\star f)} & \Tr(Y\star X')\big)\\
  \end{tikzcd}
  \end{equation}
  \end{remark}

  It is very important to note that the dg category $\hTr(\CS)$ is additive but
not triangulated. Indeed, there are lots of morphisms in $\hTr(\CS)$  which do
not exist in $\CS$, so their cones do not exist as objects in $\CS$ or in
$\hTr(\CS)$. However, we can consider its pretriangulated hull as in section
\ref{sec: pretriangulated dg}, we will denote it by $\pretr{\hTr(\CS)}$.

\begin{lemma} If $\CS$ is a dg monoidal category, then $\pretr{\hTr(\CS)}\simeq
\pretr{\hTr(\pretr{\CS})}$.
\end{lemma}
\begin{proof} As it is a dg functor $\hTr$ sends twisted complexes to twisted
complexes. Thus we have a dg functor
\[\pretr{\hTr(\pretr{\CS})} \hookrightarrow \pretr{\pretr{\hTr(\CS)}} \simeq
\pretr{\hTr(\CS)}.\] Conversely we have a dg functor
$\pretr{\hTr(\CS)}\hookrightarrow \pretr{\hTr(\pretr{\CS})}$ induced by
$\CS\hookrightarrow \pretr{\CS}$ which is a quasi-inverse.
\end{proof}

  \subsection{Homotopy trace-like functors}
  \label{section:trace-like}

  In the introduction we highlighted that the (underived) trace functor
  $\Trzz\colon \CS \to \Trzz(\CS)$ for a monoidal category (or bicategory) $\CS$
  with left duals is initial among all trace-like functors from $\CS$ to another
  category $\DS$. Here we consider a dg analog of this situation.

  Let $\CS$ be a dg monoidal category and $\DS$ a dg category. We are interested
  in dg functors $\phi \colon \CS\to \DS$ that are \emph{homotopy trace-like} in
  the sense that for every pair of objects $X,Y$ in $\CS$ we get natural (up to
  coherent homotopy) maps $\phi(X\star Y) \to \phi(Y \star X)$, possibly even
  isomorphisms. We define a \emph{dg monoidal trace} on $\CS$ with values in
  $\DS$ to be a dg functor $\phi \colon \CS \to \DS$ that factors through
  $\hTr\colon \CS \to \hTr(\CS)$. I.e. there exists a dg functor $\phi'\colon
  \hTr(\CS)\to \DS$ such that $\phi=\phi'\circ \hTr$. 
  
  Note that a dg monoidal trace $\phi$ not only contains the data of specific
  morphisms $\phi'(w_{X,Y})\colon \phi(X\star Y) \to \phi(Y \star X)$, but also
  homotopies enforcing natural compatibility relations between these morphisms.
  In particular, by Remark \ref{rem:higher traciators} $\phi'(w_{X,Y})$ is
  natural in $Y$ but it is natural in $X$ only up to coherent higher homotopies
  $\phi'(w(f_{1}\bbar \ldots\bbar f_r;Y))$. 
  
\begin{remark}
For an underived trace-like functor $\phi=\phi'\circ \hTr_0$ one requires (see
e.g \cite{MR3837875}) that the traciators are compatible with associators in the sense
that diagrams of the following type commute:
\begin{equation}
 \begin{tikzcd}
  \label{eq: trace like associativity}
   \phi((X\star Y)\star Z) \arrow{d} \arrow{r}{\phi'(w_{X\star Y Z})}  & \phi( Z \star (X\star Y)) \arrow{r} & \phi((Z \star X)\star Y)\\
   \phi(X \star (Y\star Z)) \arrow{r}  & \phi((Y \star Z)\star X) \arrow{r} & \phi(Y \star (Z\star X)) \arrow{u}.
   \end{tikzcd}
   \end{equation}
In the derived case, similar condition are required for the images of higher traciators
$\phi'(w(f_{r}\bbar \ldots\bbar f_1;Y))$.
\end{remark}   

\begin{remark}
 It is likely that if $\CS$ is has left duals then knowing $\phi'(w_{X,Y})$ and the
 naturality data $\phi'(w(f_{1}\bbar \ldots\bbar f_r;Y))$ satisfying \eqref{eq:
 trace like associativity} and its analogues is enough to prove that reconstruct
 the functor $\phi'$ and thus verify that $\phi$ is trace-like, but we do not
 prove it here.
\end{remark}

Examples of dg monoidal traces are the following.

\begin{example}
  Any endofunctor of $\hTr(\CS)$ gives rise to a dg monoidal trace on $\CS$ by
pre-composition with the universal trace $\CS \to \hTr(\CS)$. A natural source
of endofunctors of $\hTr(\CS)$ is the dg monoidal center, see
Section~\ref{sec:center-on-trace}.
\end{example}

\begin{example}
  Let $X$ be an object of $\hTr(\CS)$, then another dg monoidal trace is given
  by the universal trace $\CS \to \hTr(\CS)$ composed with the representable functor
  $\Hom_{\hTr(\CS)}(X, -)$. The target is $\dgrmod\End_{\hTr(\CS)}(X)$. A
  particularly interesting case is $X=\hTr(\one)$, for which
  Proposition~\ref{prop:verthor}   identifies the target with $\dgrmod\CCC_*(\CS)$.
\end{example}

\begin{example}
If $\CS$ is a pivotal dg monoidal category, then $\Hom_{\CS}(\one,-)$ is a
trace-like functor with target $\dgrmod\End_\CS(\one)$. To
see that $\Hom_{\CS}(\one,-)$ is trace-like, we compute:
\[
  \Hom_{\CS}(\one, X\star Y) \cong\Hom_{\CS}(X^*, Y) 
  \cong \Hom_{\CS}(\one, Y\star X^{**})\cong \Hom_{\CS}(\one, Y\star X) 
\] where the last isomorphism uses pivotality.  
Since $\Hom_{\CS}(\one,-)$ is trace-like (strictly, not just up to homotopy), it
factors through the dg monoidal trace, providing a dg functor
\[\Tr(\CS)\to \dgrmod\End_\CS(\one)\]
\end{example}

\begin{example}
\label{ex: hom with central twist}
A similar computation for pivotal $\CS$ shows that if $(Z,\tau)$ is central,
then the dg functor $\Hom_{\CS}(\one,Z\star -)$ is trace-like. More generally,
if $(Z',\tau')$ is also central, then $\Hom_{\CS}(Z',Z\star -)$ is again
trace-like. An alternative way to construct these traces is to first factor
$\Hom_{\CS}(\one,-)$ through the trace $\CS \to \hTr(\CS)$ and then act by the
dg monoidal center.
\end{example}
 
\begin{example}
The full twist $\FT$ and its powers $\FT^k$ (with suitable half-braiding
data) are objects in dg Drinfeld center of the category of complexes of Soergel
bimodules \cite{EHcenter}. By applying Example \ref{ex: hom with central
twist} we see that 
$$\Hom_{\SBim}(\one,\FT^k\star -)$$ is a trace-like functor for all $k$. Such
functors play a prominent role in the work of the first and second author,
Negu\cb{t} and Rasmussen \cite{GNR,GH}. 
\end{example}

\begin{example} If $\CS$ is a monoidal category, then any (ordinary) trace-like
  functor $\phi \colon \CS \to \DS$ lifts to a dg monoidal trace. To see this, we
  factor $\phi=\phi'\circ \Trzz$ for $\Trzz\colon \CS \to \Trzz(\CS)$ and some
  functor $\phi' \colon \Trzz(\CS)\to \DS$. The induced dg monoidal trace is
  obtained by precomposing $\phi'$ with the functor $\hTr(\CS)\to \Trzz(\CS)$ from
  Lemma~\ref{lem: derived to underived trace} and the universal trace
  $\hTr\colon \CS \to \hTr(\CS)$.
  \end{example}

\subsection{Action of the dg monoidal center on the trace}
\label{sec:center-on-trace}
It is a basic observation that the center $Z(A)$ of an associative algebra $A$ acts on
the trace $A/[A,A]$. Similarly, if we think of the trace $\Trzz(\CS)$ of a
monoidal category as $\CS$ integrated over an annulus, the Drinfeld center
$\ZC^{0}(\CS)$ acts by ``cutting open'' the annulus and ``gluing in'' central
objects or morphisms, before ``resealing the annulus''. In this section we
describe the analogous action of the dg monoidal center $\ZC^{\text{dg}}(\CS)$
on the dg monoidal trace $\hTr(\CS)$.

Let $(Z, \tau)$ be an object of $\ZC^{\text{dg}}(\CS)$.  
The structure map $\tau$ is a map of $\CS,\CS$-bimodules $\BB(\CS)\rightarrow \XS_{12}(Z,Z)$.  
Using this map we have a map of $\CS,\CS$-bimodules given by composing:
\begin{eqnarray*}
\BB(\CS)\otimes_\CS \XS_{12}(X',X)
&\rightarrow & \BB(\CS)\otimes_\CS \BB(\CS)\otimes_\CS \XS_{12}(X',X)\\
& \rightarrow & \BB(\CS)\otimes_\CS \XS_{12}(Z,Z)\otimes_\CS \XS_{12}(X',X)\\
& \rightarrow & \BB(\CS)\otimes_\CS  \XS_{12}(Z\star X',Z\star X).
\end{eqnarray*}
The first of these maps is the comultiplication on the two-sided bar complex of $\CS$, the
second is an application of $\tau$, and the last is the algebra structure on the
bimodule $\XS_{12}$ (see Remark~\ref{rem:quadmod-nw-algebra}). 
Applying the functor which identifies
the left and right actions of $\CS$, we obtain a map of complexes
\[
\Hom_{\Tr(\CS)}(\Tr(X),\Tr(X')) \rightarrow \Hom_{\Tr(\CS)}(\Tr(Z\star X),\Tr(Z\star X')).
\]

In this way, $(Z,\tau)$ determines an endofunctor $\Xi_{(Z,\tau)}\colon
\Tr(\CS)\rightarrow \Tr(\CS)$ defined on objects by $\Tr(X)\mapsto \Tr(Z\star
X)$ and on morphisms by the above chain map. That this map respects composition
of morphisms follows from the assumption that $\tau$ is a morphism of dg
algebras and the compatibility of the shuffle product and coproduct on the
two-sided bar complex.

The action on morphism complexes can be pictured as follows.
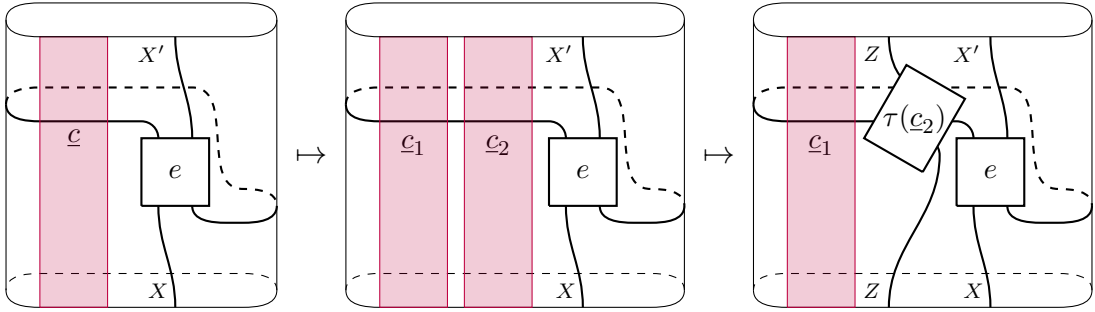
\[
\centonhtrfigure{.9}
\]

\begin{remark}
Morphisms $f\colon (Z,\tau) \rightarrow (Z',\tau')$ in $\ZC^{\text{dg}}(\CS)$
should give natural transformations $\Xi_f\colon \Xi_{Z, \tau} \to \Xi_{Z',
\tau'}$ that assemble into a monoidal functor
\[\Xi\colon \ZC^{\text{dg}}(\CS) \to \End(\hTr(\CS))\] from
$\ZC^{\text{dg}}(\CS)$ to the endofunctors of $\Tr(\CS)$.
We will not pursue this in detail, see also Remark~\ref{rem:caveat}.
\end{remark}

\section{Traces of the Soergel category}
\label{sec:HH-Scat}
\subsection{Soergel bimodules}
\label{subsec: SBim def}
Let $W$ be a Coxeter group with simple reflections $S\subset W$, length function
$\ell$, and Bruhat order $\leq$.  Let $(V, \{\alpha_s^{\vee}
\}, \{\alpha_s\})$ be a realization of $W$ in the sense of \cite{EW}, over $\C$, which we assume to be reflection
faithful and balanced unless stated otherwise. As usual, we consider the
polynomial ring $R=\C[V]:=\Sym^\bullet(V^\ast(-2))$, which is graded by declaring
elements in $V^*$ to be of degree $2$. In particular, we have $\alpha_s\in R$,
and these elements generate $R$ if they span $V^*$. 

\begin{remark}
In type $A$ we have $W=S_n$ with simple reflections (transpositions) indexed by
$i\in \{1,\dots, n-1\}$. We will consider the realization with $V=\C^n$, on
which $S_n$ acts by permuting standard basis vectors, and we identify
$R=\C[x_1,\ldots,x_n]$ and $\alpha_i=x_i-x_{i+1}$ for $1\leq i\leq n-1$.
\end{remark}

We now describe the category of Soergel bimodules $\SBim(W)$ associated to $W$ and its
chosen realization.  For each simple reflection $s\in S$ we let
$$B_s:=R\otimes_{R^s} R(1),$$ where $R^s = \{f\in R \:|\: s(f)=f\}$ and $(1)$ is
the ``downward'' grading shift. Then the category of Soergel bimodules
$\SBim(W)$ is the smallest full monoidal subcategory of graded $R,R$-bimodules
containing the bimodules $B_s$ for $s\in S$ whichis closed under grading shift, isomorphism, direct sums,
and direct summands. The monoidal structure on $\SBim(W)$ is denoted by $\star$
and tensor products of bimodules of the form $B_s$ are called
\emph{Bott-Samelson bimodules}. 

In this paper we will occasionally use the diagrammatic Hecke category $\dsc{W}$
of Elias--Williamson \cite{EW}, whose Karoubi completion is equivalent to the
category of Soergel bimodules under the assumptions taken. The
diagrammatic category $\dsc{W}$ has the advantages of being strictly monoidal
and manifestly strictly pivotal (by the balanced assumption) besides being
well-behaved for a larger class of realizations (which will not be relevant
here).

Let $\CS(W):=\Ch^b(\SBim(W))$ be the dg monoidal category of bounded chain
complexes of Soergel bimodules for $W$. If $\underline{\b}$ is a word in the
alphabet $\sigma_s, \sigma_s\inv$ with $i\in S$, then we have a finite complex
$F(\underline{\b})\in \CS(W)$ defined by
\begin{gather}
  \label{eqn:rouquiercx}
F(\sigma_s) := (\underline{B_s}\rightarrow R(1)),\qquad F(\sigma_s\inv) := 
(R(-1)\rightarrow \underline{B_s}),\\ 
\nonumber
F(\underline{\b}\cdot \underline{\b}') := 
F(\underline{\b})\star F(\underline{\b}'),
\end{gather}
where the maps $B_s\rightarrow R(1)$ and $R(-1)\rightarrow B_s$ are the canonical
bimodule maps, defined by
\[
1\otimes 1\mapsto 1,\qquad\qquad 1\mapsto \frac{1}{2}(\a_s\otimes 1 - 1\otimes \a_s),
\]
respectively. It is well known \cite{R,R2} that the complexes
$F(\underline{\b})$ satisfy the braid relations up to homotopy equivalence.  We
often abusively write $F(\b)$ where $\b\in \Br(W)$ is an element of the braid
group associated to $W$ (when in reality the complex $F(\b)$ depends on a choice
of braid word $\underline{\b}$ representing $\b$).

\begin{definition}
For each $w\in W$, let $\Delta_w$ and $\nabla_w$ denote the Rouquier complex of
the positive and negative braid lift of a chosen reduced expression of $w$.
\end{definition}

In particular, we have $\Delta_w \star \nabla_{w^{-1}} \simeq \one \simeq
\nabla_{w^{-1}} \star \Delta_w$ for any $w\in W$. We will sometimes write
$\Delta_w^{-1} := \nabla_{w^{-1}}$.

\subsection{Hochschild homology of the Soergel category --- linear structure}

The following is well known to experts. 

\begin{proposition}
  \label{prop:Soergel-semi-orth}
The complexes $\Delta_w$, $w\in W$, generate $\CS(W)$ with respect to
cones, shifts, sums, and homotopy equivalences (and similarly for $\nabla_w$).
These complexes satisfy
\[
\Hom_{\CS(W)}(\Delta_v, \Delta_{w})\simeq 0 \simeq \Hom_{\CS(W)}(\nabla_w,\nabla_v)\qquad\text{unless}\qquad v\leq w.
\]
In other words, $\{\Delta_v\}_{v \in W}$ and $\{\nabla_v\}_{v\in W}$ each
generate a semi-orthogonal decomposition of $\CS(W)$ (the latter one
with the opposite poset structure).
\end{proposition}

\begin{proof}
By the main result of \cite{LW}, one has
$$
\Hom_{\CS(W)}(\Delta_v,\nabla_w)
\simeq
\begin{cases}
R & \text{if}\ v=w,\\
0 & \text{otherwise}.\\
\end{cases}
$$
On the other hand, it is easy to see that $\Delta_w$ is filtered by $\nabla_u$
for $u\leq w$. Therefore if $\Hom_{\CS(W)}(\Delta_v,\Delta_w)\not\simeq 0$ then
$\Hom_{\CS(W)}(\Delta_v,\nabla_u)\not\simeq 0$ for some $u\leq w$, hence $v=u\leq w$. See also
\cite[Appendix]{GHMN} for more details.
\end{proof}

Thus the following is an immediate consequence of
Theorem~\ref{thm:barcx-retract}.
\begin{lemma}
There is a deformation retract from the cyclic bar complex $\CCC(\SBim(W))$ to
$\bigoplus_{w\in W} \CCC(\End(\Delta_w))$. \qed
\end{lemma}

Recall that we have chosen a realization $V$ of $W$.  Consider the following
complex of $\C$-vector spaces

\[
Z:= \left(V^\ast(-2) \rightarrow \underline{V^\ast(-2)\oplus V^\ast(-2)}\right),\qquad \qquad \phi\mapsto (\phi,-\phi).
\]

Complexes of $\C$-vector spaces form a symmetric monoidal (dg) category, and so
we have Schur functors.  The symmetric algebra of $Z$ is just
\[
\Sym^\bullet(V^\ast(-2)\oplus V^\ast(-2) \oplus \susp{1}V^\ast(-2)) \ \cong \ R\otimes R \otimes \Lambda,
\]
with its differential inherited from $Z$.  Here $\Lambda$ is the exterior algebra of $V^{\ast}$.
After choosing a basis of $V$ and
letting $x_1,\ldots,x_n\in V^\ast$ denote the dual basis, we can identify
\[
\Sym^\bullet(Z) \cong \C[x_1,\ldots,x_n,x_1',\ldots,x_n',\theta_1,\ldots,\theta_n],\qquad\qquad d(\theta_i) = x_i -x_i'.
\]
That is to say, $\Sym^\bullet(Z)$ is the Koszul resolution of $R$ as a bimodule over itself.
We can use this resolution to compute $\HH_{\subb}(R)$:
\[
\HH_{\subb}(R) \cong R\otimes \Lambda.
\]

We have a dg algebra map $R=\End_{\CS(W)}(\one)\rightarrow \End_{\CS(W)}(\Delta_w)$ sending
$f\mapsto \Id_{\Delta_w}\star f$.  The homotopy equivalence
\[
\Hom_{\CS(W)}(\Delta_w,\Delta_w) \simeq \Hom_{\CS(W)}(\Delta_w\inv\star \Delta_w, \one) \simeq \End_{\CS(W)}(\one) = R
\]
is clearly $R$-linear. 

Thus, combining everything up to this point gives:
\begin{corollary}
\label{cor:cbarcx-SBim}
There is a deformation retract from the cyclic bar complex
$\CCC_*(\CS(W))$ to its homology $R\otimes \Lambda \rtimes \C[W]$. Here the
generators $x_i$ of $R$ have degree $(2,0)$, the generators $\theta_i$ of
$\Lambda$ have degree $(2,-1)$ and $W$ has degree $(0,0)$.
\end{corollary}

\subsection{Hochschild homology of the Soergel category --- algebra structure}
 
 In this section we prove Theorem~\ref{thm:main}. Recall that by
 Theorem~\ref{thm:Ainfty-defretracts} any deformation retract $\CCC_*\twoheadrightarrow H$
 induces an $A_\infty$ structure on $H$. 
 
 Recall that the dg algebra $\CCC_*$ is called formal if all higher Massey
 products vanish (i.e.~the higher $\Ai$-maps on $H$ are trivial).
 
 \begin{theorem}
  \label{thm:HHSBim}
  We have an isomorphism of algebras
$$
\HH_{\subb}(\SBim(W))\cong R\otimes \Lambda\rtimes  W
 $$
with multiplication $\mu_2$. 
 \end{theorem}
 
 \begin{proof}
By Theorem \ref{thm:cyclicbarcx-retract} we have a canonical isomorphism
  \[\HH_{\subb}(\SBim(W))\cong \HH_{\subb}(\CS(W)),\] and we will compute the latter. All
  homs in this proof are taken in $\CS(W)$ and so the subscripts in
  $\Hom_{\CS(W)}(-,-)$ and $\End_{\CS(W)}(-)$ will be omitted. First we observe that
  $\HH_{\subb}(\End(\one))\cong \HH_{\subb}(R) \cong R \otimes \Lambda$ is a subalgebra of
  $\HH_{\subb}(\CS(W))$ by definition. Second, let $\sigma_w$ denote the class of the
  identity morphism of $\End(\Delta_w)$ in $\HH_{\subb}(\CS(W))$. Then
  $\End(\Delta_w)\simeq R$ is a bimodule over $\End(\one)\cong R$ where the
  right action of $R$ is standard and the left action of $R$ is twisted by $w$,
  and the same holds for the actions of $\HH_{\subb}(R)$ on
  $\HH_{\subb}(\End(\Delta_w))\simeq \HH_{\subb}(R)$. In other words,
$$
f(x,\theta)\cdot \sigma_w=\sigma_w\cdot f(w^{-1}(x),w^{-1}(\theta)).
$$
 It remains to prove that $w\mapsto \sigma_w$ is a group homomorphism. 
 For this, we construct an action of $\sigma_w$ on $\End(\one)$ by the following
 composition of homotopy equivalences (where all homs are taken in $\CS(W)$):
$$
\sigma_w\colon \End(\one)\simeq \Hom(\one,\Delta_w\Delta_w^{-1})
\simeq \Hom(\Delta_w,\Delta_w)\simeq \Hom(\one,\Delta_w^{-1}\Delta_w)\cong \End(\one).
$$
Here the middle isomorphism follows from the fact that $\Delta_w$ and
$\Delta_w^{-1}$ are biadjoint in $\CS(W)$. One can check that this action is
compatible with the shuffle multiplication in $\HH_0(\CS(W))$ and $\sigma_w$
acts as a permutation $(-1)^{\ell(w)}w$ on $\End(\one)$ (the signs that appear
here depend on the homological shift conventions for Rouquier complexes, see
\cite[Section 4.3 and Remark 4.26]{GW} for a discussion).
Since the representation of $W$ on $\End(\one)$ is faithful, we get that the
$\sigma_w$ generate a copy of $W$ inside $\HH_0(\CS(W))$.
 \end{proof}
 
\begin{remark}
The main result of \cite{MR3448186} states that $\HH_0(\SBim(W))\cong R\rtimes \C[W]$, so
our theorem is a natural generalization. The methods of proof in \cite{MR3448186},
however, were completely different and used cellularity of $\SBim(W)$.
\end{remark} 

\begin{remark}
One can check that the classes $\sigma_w=[\Delta_w]$ generate a copy of $W$ in a more direct way. 
Indeed, for simple reflections $s$ the objects $\Delta_{s}$ satisfy braid relations in $\CS(W)$,
and hence $[\Delta_{s}]$ satisfy them too. Furthermore,
$$
[\Delta_{s}]^2=[(\Delta_{s})^2]\simeq [\underline{B_s}(-1)\to B_s(1)\to R(2)]=[R]-[B_s]+[B_s]=[R]=1.
$$
The middle equation follows from \eqref{eq:cone}.
\end{remark}
 
Next we would like to understand the $\Ai$-structure on  $\HH_{\subb}(\SBim(W))$.

We can also describe the action of the Connes differential on Hochschild homology.

\begin{proposition}
The action of the Connes differential $\B$ on $\HH_{\subb}(\SBim(W))$, pulled back to
$R\otimes \Lambda \rtimes \C[W]$ via the isomorphism from
Theorem~\ref{thm:HHSBim}, is given by
\begin{equation}
\label{eq:connes for SBim}
\B=\sum_i \theta_i \frac{\partial}{\partial x_i},
\end{equation}
 where ${x_i,\theta_i}$ are dual bases in $V,V^*$. 
\end{proposition}

\begin{proof}
Recall the isomorphisms
\[ \HH_{\subb}(\SBim(W))\cong \HH_{\subb}(\CS(W))\cong \oplus_{w\in W} \HH_{\subb}(\End(T_w))\cong
\oplus_{w\in W} \HH_{\subb}(R).\] Observe that the second isomorphism is obtained by
certain retractions which commute with the action of the Connes differential
$\B$, so it is sufficient to know the action of $\B$ on
$\HH_{\subb}(\End(T_w))=\HH_{\subb}(R)$. The differential $\B$ on $\HH_{\subb}(R)$ is given by
\eqref{eq: de rham}, and it is obviously $W$-invariant. 
\end{proof}

\subsection{Towards understanding formality}

Recall that we have a chosen realization of the Coxeter group $W$, denoted $V$,
and $R=\Sym(V^\ast(-2))$ is the ring of polynomial functions on $V$, bigraded by
declaring linear functions to have degree $(2,0)$. Let us choose once and for
all pair of dual bases for $V$ and $V^\ast$.  After this choice, we will write
$R=\C[x_1,\ldots,x_n]$, where $x_i$ range over the dual basis in $V^\ast(-2)$.
We will need the following auxilliary objects as well.

Let $E$ denote the bigraded algebra freely generated by elements
$\xi_i,\xi_i^\ast$ ($1\leq i\leq n$) of degree $\deg(\xi_i)=(2,-1)$,
$\deg(\xi_i^\ast)=(-2,1)$, modulo
\[
[\xi_i,\xi_j]=0,\qquad[\xi_i^\ast,\xi_j^\ast],\qquad [\xi_i^\ast,\xi_j] =\begin{cases}1 & \text{if $i=j$}\\ 0 & \text{else} \end{cases}
\]

Let $M=E/I$, where $I\subset E$ is the left ideal generated by the $\xi_i^\ast$
($1\leq i\leq n$).  

We think of $\xi_i$ and $\xi_i^\ast$ as ranging over bases of $V^\ast$ and $V$,
respectively; in basis independent language one would write
$M=\Sym(V^\ast[1](-2))$ (which is the exterior algebra on $V^\ast$ with some
regrading), and $E=\End_\C(M)$.

\begin{definition}
We define the {\em cube complex} $K\in \CS(W)$ as the  twisted complex
$K=\one\otimes_\C M$ with twist $\alpha_K=\sum_i x_i\otimes \xi_i^\ast$.  Note
that this is simply the Koszul complex (an $n$-dimensional cube-like complex)
associated to the action of $x_1,\ldots,x_n$ on $\one$.
\end{definition}

\begin{example}
For $W=S_2$ and $V=\C$ (the reflection representation) we get
$$
K= \left( \one\xrightarrow{x_1-x_2} \one\right).
$$
\end{example}
\begin{example}
For $W=S_2$ and $V=\C^2$ (the permutation representation) we get
$$
K=\left( \one \xrightarrow{(x_1,x_2)}\C^2\otimes \one\xrightarrow{(x2,-x1)}\one  \right).
$$
\end{example}

Observe that  $\Tr(K)$ is the  twisted complex $\Tr(\one)\otimes_\C M$ with
twist $\alpha_K$. The action of $W$ on $\Tr(\one)$ and on $M$ induces the action
of $W$ on $\Tr(K)$, as the twist is $W$-invariant. The operators
$x_i,\theta_i\in \End_{\hTr}(\Tr(\one))$ also act on $\Tr(K)$, although it is
easy to see that the action of $x_i$ on $\Tr(K)$ is null-homotopic. Finally,
multiplication by $\xi_i^*$ also acts on $\Tr(K)$.

\begin{theorem}
\label{thm: K formal}
The endomorphism dg algebra $\End_{\hTr(\CS(W))}(\Tr(K))$ is formal and thus
quasi-isomorphic to its cohomology $\C[\theta_i,\xi_i^\ast] \rtimes \C[W]$. 
\end{theorem}

\begin{proof}
It is easy to see that $\End_{\hTr(\CS(W))}(\Tr(K))$ is quasi-isomorphic to the total complex of the bicomplex 
\begin{equation}
\label{eq: end K complex}
\left(\End_{\hTr(\CS(W))}(\Tr(\one))\otimes \End_{\C}(M), D\right)
\end{equation}
where the differential $D$ consists of the internal differential $d_{\CCC}$ of
$\End_{\hTr(\CS(W))}(\Tr(\one))$, which we identity with $\CCC_*(\CS(W))$, and the ``cube differential''
$d_1$ corresponding to the commutation with the twist $\alpha_K$.  

We can compute the homology of \eqref{eq: end K complex} using the following
spectral sequence. First we compute the homology of $d_{\CCC}$ and obtain
\[H^{\supb}(\End_{\hTr(\CS(W))}(\Tr(\one))\simeq \HH_{\subb}(\CS(W))\cong R\otimes \wedge(\theta_i)\rtimes
\C[W],\] 
so the $E_1$ page of the spectral sequence has the form
$$
E_1=H^{\supb}((\End_{\hTr(\CS(W))}(K),d_{\CCC})=R\otimes \C[\theta_i] \rtimes \C[W]\otimes  \End_{\C}(M).
$$
Next, we compute the homology of the differential $d_1$. We have
$$
H^{\supb}(R\otimes  \End_{\C}(M),[\alpha_K,-])\simeq H^{\supb}(R\otimes  E,[\alpha_K,-])\simeq  \wedge(\xi^*_i).
$$
Since the twist commutes with $W$ and all $\theta_i$, 
we get the $E_2$ page
$$
E_2=H^{\supb}(E_1,d_1)=\C[\theta_i,\xi_i^\ast] \rtimes \C[W].
$$
Observe that here $\theta_i$ have bidegree $(2,-1)$,
while $\xi^*_i$ have degree $(-2,1)$ and $\C[W]$ has
bidegree $(0,0)$. Therefore the bidegrees of all
homology generators in $E_2$ page are concentrated on the line $(2s,-s)$. 

On the other hand, any higher differentials $d_r$ have (total) 
degree $(0,1)$, so they cannot preserve this line and must
vanish. We conclude that the spectral sequence collapses at the $E_2$ page.

Finally, we need to prove that the Massey products $\mu_r$ vanish for $r\ge 3$.
Note that $\mu_r$ shift the homological degree by $(3-r)$ and do not change the
$q$-degree. Since all generators in the homology of $\End_{\hTr(\CS(W))}(K)$ are
concentrated on the line $(2s,-s)$, and $\mu_r$ move them out of this line for
$r\ge 3$, we conclude that $\mu_r=0$.
\end{proof}

\begin{corollary}
There is a functor $\Hom_{\hTr(\CS(W))}(\Tr(K),-)$ from the derived horizontal
trace $\hTr(\CS(W))$ to the category of $\Ai$-modules over dg algebra
$\End_{\hTr}(\Tr(K))$.
\end{corollary}

\begin{remark}
Although the dg algebra $\End_{\hTr(\CS(W))}(\Tr(K))$ is formal, it could have
nontrivial $\Ai$-modules. 
\end{remark}

\begin{example}
Similarly to the above, it is easy to see that
\[\Hom_{\hTr(\CS(W))}(\Tr(K),\Tr(\one))\simeq \C[\theta_i]\rtimes \C[W]\] as a
module over $\End_{\hTr}(\Tr(K))$.  However, we do not know if it has nontrivial
$\Ai$-products 
$$
\mu_r:\End_{\hTr(\CS(W))}(\Tr(K))^{\otimes r-1}\otimes \Hom_{\hTr(\CS(W))}(\Tr(K),\Tr(\one))\to \Hom_{\hTr(\CS(W))}(K,\Tr(\one)).
$$
\end{example}

\begin{conjecture}
The cyclic bar complex $\CCC_*(\CS(W))\simeq \End_{\hTr(\CS(W))}(\Tr(\one))$ is
formal as a dg algebra. In particular, all maps $\mu_d$ on the Hochschild homology 
$\HH_{\subb}(\CS(W))$ vanish for $d\ge 3$.

\end{conjecture}

\begin{remark}
The previous version of this paper on arXiv contained an incomplete proof of
this conjecture. The key step in the argument was using that the all deformation
retract data were $R$-bilinear. This is indeed the case except possibly for the
crucial homotopy retracting $\Hom(\Delta_{v},\Delta_{w})\simeq 0$  implied by
Proposition \ref{prop:Soergel-semi-orth}. One can choose a homotopy which is
$R$-linear on the left (or on the right), but it not known to us if we can
choose an $R$-bilinear homotopy. Therefore, it is unclear if one can choose an $R$-bilinear
retraction of the whole cyclic bar complex of $\CS(W)$ to its homology.
\end{remark}

We note that the Connes differential $\B$ induces an interesting endomorphism of
$\Tr(K)$. Recall that $\B=\sum \theta_i \frac{\partial}{\partial x_i}$, and we
can write the differential on $K$ as $d=\sum x_i\xi_i^{\ast}$, so $\B(d)=\sum
\theta_i \xi_i^{\ast}$. Here $\theta_i$ and $\xi_i^*$ correspond to two dual
bases in $V^{\ast}$ and $V$ respectively, so $\B(d)$ is simply the canonical
tensor in $V^{\ast}\otimes V$. It has bidegree $(0,0)$ and is invariant under
the action of $W$. 

\begin{example}
For $W=S_2$ and $V=\C$ we get the following endomorphism of $K$:
\begin{equation*}
\begin{tikzcd}
K=\big( \Tr(\one) \arrow{r}{x_1-x_2} \arrow{dr}{\theta_1-\theta_2}& \Tr(\one) \big) \\
K=\big( \Tr(\one) \arrow{r}{x_1-x_2} & \Tr(\one)\big).
\end{tikzcd}
\end{equation*}
\end{example}

\subsection{The dg monoidal trace of the Soergel category}
\label{sec:compactgen}
In this section, we let $W$ be an arbitrary finite Coxeter group and again write
$\CS(W)=\Ch^b(\SBim(W))$ for the dg monoidal category of bounded chain
complexes of Soergel bimodules for $W$ with a given
realisation. 

Any object in $\hTr(\CS(W))$ is homotopy equivalent to a twisted complex built
out of finite direct sums of $\hTr(B_w)$. The same applies to any object in
$\pretr{\hTr(\CS(W))}$. Indeed, any object $X$ in $\CS(W)$ is a complex built
out of finite direct sums of $B_w$, which we can write as an iterated cone.
Since $\hTr$ is a dg functor, we can write $\Tr(X)$ as an iterated cone built
out of $\hTr(B_w)$ which is a twisted complex. Note that $\hTr$ indeed sends a
complex to a twisted complex in general, see Example \ref{ex: full twist
derived}.

\begin{proposition}
Let $X$ and $Y$ be two objects in $\CS(W)$. Then $\Tr(Y)\Tr(\CS(W))\Tr(X)$ is homotopy
equivalent to a bounded complex of free finitely generated $R$-modules.  
\end{proposition}

In particular, the homs in $\hTr(\CS(W))$ are finite-dimensional in each bidegree.

\begin{proof}

Here we use the semiorthogonal decomposition from
Proposition~\ref{prop:Soergel-semi-orth} again, and we write $\CS=\CS(W)$ As in
Section~\ref{ss:semiortho}, the hom complexes in $\Tr(\CS)$ retract onto
semi-orthogonal hom complexes, i.e. for objects $X$ and $Y$ we have a retraction

\begin{align*}
  \Tr(Y)\Tr(\CS\Tr(X) & = \BB(\CS \otimes_{\CS\boxtimes \CS^{\op}}  \XS_{12}(Y,X)  \twoheadrightarrow 
\BB^\Gamma(\CS,\IS_\BS) \otimes_{\CS\boxtimes \CS^{\op}} \XS_{12}(Y,X),
\end{align*}
and further onto a complex which has chain groups $\bigoplus_{w_0}
(\Delta_{w_{0}}Y)\CS(X\Delta_{w_0})$ in degree $r=0$ and for $r>0$:
$$
\bigoplus_r \bigoplus_{w_0,\ldots,w_r} 
\susp{r} \Big( \Delta_{w_0}\CS\Delta_{w_1}\otimes \cdots \otimes \Delta_{w_{r-1}}\CS
\Delta_{w_r} \otimes (\Delta_{w_{r}}Y)\CS(X\Delta_{w_0}) \Big)
$$
with $w_0\leq w_1\leq\ldots \leq w_r$ in Bruhat order. This still leaves us
with an infinite complex since $w_i$ could repeat arbitrarily. However, such
infinite repeats form a copy of the two-sided bar complex of $\End(\Delta_w)=R$, which can
be retracted to a finite complex. Since $W$ was assumed to be finite, which
implies that chains in Bruhat order are finite, performing all such retractions we arrive at
a homotopy equivalent bounded complex. Freeness now follows since
morphism spaces in $\CS$ are free over $R$, see \cite{MR2329762} or \cite[Theorem
3.6]{MR3245013}.
\end{proof}

\section{The Soergel category in type A}
\label{sec:type A}
In this section we describe the derived horizontal trace of the Soergel category
in type A, using the computation of its derived vertical trace in
Section~\ref{sec:HH-Scat}. 

In Proposition~\ref{prop:verthor} we have seen that the (derived) vertical trace
$\CCC_*(\CS)$ of a monoidal (dg) category $\CS$ can be identified with the
endomorphisms of $\hTr(\one)$ in the (derived) horizontal trace $\hTr(\CS)$.
That this in fact determines the entire horizontal trace for Soergel bimodules
of type A is the upshot of the following theorem.

\begin{theorem}
  \label{th: annular simplification}
  The Karoubi completion of the triangulated hull of $\hTr(\SBim_n)$ is split-generated by
  $\Tr(\one)$. We have
  $$
  \Kar^{dg}(\pretr{\hTr(\SBim_n)})\simeq \Perf(\C[x_1,\ldots,x_n, \theta_1,\ldots, \theta_n] \rtimes \C[S_n]\mod).
  $$
  Here in the right hand side we have the category of perfect $\Ai$-modules over the $\Ai$-algebra
  $$
  \HH_{\subb}(\SBim_n)\cong R\otimes \Lambda \rtimes \C[S_n]=\C[x_1,\ldots,x_n, \theta_1,\ldots, \theta_n] \rtimes \C[S_n],
  $$
  where $R$ and $S_n$ are supported in cohomological degree zero and the variables $\theta_i$ have cohomological degree $-1$.
  \end{theorem}
  
   The proof of this theorem will occupy the rest of this section. Here we
   outline the strategy of the proof. First, we use the technology of Frobenius
   extensions to prove that $\Tr(B_{I})$ is a direct summand in a direct sum of
   several copies of $\Tr(\one)$, where $B_I$ is the Soergel bimodule
   corresponding to the longest element in a parabolic subgroup $W_I\subset W$,
   for any subset $I\subset S$ of simple transpositions. Next, we use an
   explicit ``annular simplification'' algorithm to present $\Tr(B)$ as a
   direct summand in a direct sum of $\Tr(B_I)$ for any Soergel bimodule $B$.
   Finally, we show that any complex of Soergel bimodules is mapped by the trace
   functor to a summand in a finite twisted complex built out of $\Tr(\one)$,
   thus completing the proof.

\subsection{Frobenius extensions and horizontal trace}
In this self-contained section we show that Soergel bimodules $B_I$ associated
to longest elements of finite parabolic subgroups $W_I$ of a Coxeter group $W$
always have traces $\Trzz(W_I)$ isomorphic to summands of $\Trzz(\one)$, after
Karoubi completion. We prove this in a slightly more general setting, using the
language of Frobenius extensions.

\begin{definition}
  A \emph{Frobenius extension} is an extension of commutative rings $\iota\colon
  A\hookrightarrow B$, such that $B$ is free and finitely generated as an
  $A$-module, equipped with a non-degenerate $A$-linear map $\partial\colon B\to
  A$, called the \emph{trace}. Here, \emph{non-degeneracy} asserts the existence
  of $A$-linear \emph{dual bases} $\{x_\alpha\}$ and $\{y_\alpha\}$ for $B$ such
  that $\partial(x_\alpha y_\beta)=\delta_{\alpha, \beta}$ (the Kronecker
  delta). 
 \end{definition}

\begin{example}
  \label{exa:FrobExt1}
  Let $(W,S)$ be a Coxeter system of finite rank as in \S \ref{sec:HH-Scat} and
  let $R$ again denote the base ring associated with a reflection faithful
  balanced realization over $\C$. For every finite parabolic subgroup $W_I$, let
  $R^I$ denote the subring of $W_I$-invariants in $R$.  Let $w_I\in W_I$ be the
  longest element.  Then $\iota\colon  R^I \to R$ is a graded Frobenius
  extension of rank $|W_I|$ with trace
  $\partial=\partial_I=\partial_{s_1}\cdots\partial_{s_r}$ where $w_I=s_1\cdots
  s_r$ is a reduced expression and
  \[\partial_s(f)=\frac{f-s(f)}{\alpha_s}.\]
  See e.g. \cite[Section 3]{Williamson-thesis}.
\end{example}

\begin{remark}
The subsets $I$ corresponding to finite parabolic subgroups $W_I$ are called
{\em finitary}. Note that we do not need $W$ to be finite in Example
\ref{exa:FrobExt1}. 
\end{remark}

\begin{example}
In type A, we have $R=\C[X_1,\dots,\lon_n]$ and $W=S_n$ and $R^{S_n}\hookrightarrow
R$ is a graded Frobenius extension of rank $n!$. An $R^{S_n}$-linear basis of
$R$ is given by the monomials $X_1^{a_1}X_2^{a_2}\cdots X_{n-1}^{a_{n-1}}$ where
$0\leq a_i\leq n-i$. Then we have 

\[\partial(X_1^{a_1}X_2^{a_2}\cdots X_{n-1}^{a_{n-1}})=\begin{cases} 1 &\text{ if } a_i=n-i \text{ for } 1\leq i \leq n,\\
0 & \text{ otherwise}\end{cases}.\] The basis dual to the monomial basis has
 elements $\prod_{k=1}^{n-1} (-1)^{b_k}e_{b_k}(X_{n+1-k},\dots, X_{n})$ where
 $b_k=k-a_{n-k}$ for $1\leq k\leq n-1$.
\end{example}

Given a Frobenius extension $\iota \colon A \to B$ with trace $\partial$, 
  we have the following maps of $B,B$-bimodules, which are best encoded diagrammatically: 
  \begin{itemize}
    \item the ``multiplication'' $\trival{blue}{.2}{.2}\colon B\otimes_A B\otimes_A B\to B\otimes_A A \otimes_A B \cong B\otimes_A B$, 
    \item the ``inclusion'' $\Sdot{blue}{.2}{-.2}\colon B \to B \otimes_A B$
    defined by $1\mapsto \sum_{\alpha} x_{\alpha} \otimes y_{\alpha}$, and 
    \item the ``trace'' $\Sdot{blue}{.2}{.2}\colon B\otimes_A B \to B$ given by $x\otimes y \mapsto xy$, 
  \end{itemize}
  which exhibit $B\otimes_A B$ as a Frobenius extension of $B$. 
  This is an instance of Jones' basic construction \cite{MR0696688}. Further we have:
 \begin{itemize}
    \item $\trival{blue}{.2}{-.2}\colon B\otimes_A B \to B\otimes_A B\otimes_A B$ given by $1\otimes 1 \mapsto 1\otimes 1 \otimes 1$.
\end{itemize}
These satisfy the relations:
$$
\begin{tikzpicture}[anchorbase,scale=.2]
  \draw[blue, line width=.5mm] 
  (-1,-1) to (0,0) to (1,-1) 
  (0,0) to (0,1)
  (-1,2) to (0,1) to (1,2) ;
\end{tikzpicture}
\;=\;
\begin{tikzpicture}[anchorbase,scale=.2]
  \draw[blue, line width=.5mm] 
  (-1,-1) to (-1,2) 
  (-1,0) to (1,1)
  (1,-1) to (1,2) ;
\end{tikzpicture}
\;=\;
\begin{tikzpicture}[anchorbase,scale=.2]
  \draw[blue, line width=.5mm] 
  (-1,-1) to (-1,2) 
  (-1,1) to (1,0)
  (1,-1) to (1,2) ;
\end{tikzpicture}
\,
,\quad
\sidedot{blue}{.2}{.2}
\;=\;
\sidedot{blue}{-.2}{.2}
\;=\;
\sidedot{blue}{.2}{-.2}
\;=\;
\sidedot{blue}{-.2}{-.2}
\;=\;
\begin{tikzpicture}[anchorbase,scale=.2]
  \draw[blue, line width=.5mm] 
  (0,-1) to (0,2);
\end{tikzpicture}
\,
,\quad
\begin{tikzpicture}[anchorbase,scale=.2]
  \draw[blue, line width=.5mm] (0,0) to (0,-1);
  \fill[blue] (0,0) circle (2.5mm);
  \draw[blue, line width=.5mm] (0,1) to (0,2);
  \fill[blue] (0,1) circle (2.5mm);
\end{tikzpicture}
\;=\; \sum_{\alpha} x_\alpha
\, 
\begin{tikzpicture}[anchorbase,scale=.2]
  \draw[blue, line width=.5mm] (0,2) to (0,-1);
\end{tikzpicture}
\,
y_\alpha
,\quad
\begin{tikzpicture}[anchorbase,scale=.2]
  \draw[blue, line width=.5mm] 
  (0,-1) to  (0,-.5) to (-1,0) to (-1,1) to (0,1.5) to (0,2)
  (0,-.5) to (1,0) to (1,1) to (0,1.5);
  \node at (0,.5) {\tiny$p$};
\end{tikzpicture}
\;=\;
\partial(p) \,
\begin{tikzpicture}[anchorbase,scale=.2]
  \draw[blue, line width=.5mm] 
  (0,-1) to (0,2);
\end{tikzpicture}
$$
We also define the following shorthand notation:
$$
\begin{tikzpicture}[anchorbase,scale=.2]
  \draw[blue, line width=.5mm] 
  (-1,-1) to [out=90,in=180] (0,0) to [out=0,in=90] (1,-1);
\end{tikzpicture}
\;:=\;
\begin{tikzpicture}[anchorbase,scale=.2]
  \draw[blue, line width=.5mm] 
  (-1,-1) to (0,0) to (1,-1) 
  (0,0) to (0,1);
  \fill[blue] (0,1) circle (2.5mm);
\end{tikzpicture}
\,
,\quad 
\begin{tikzpicture}[anchorbase,scale=.2]
  \draw[blue, line width=.5mm] 
  (-1,2) to [out=270,in=180] (0,1) to [out=0,in=270] (1,2);
\end{tikzpicture}
\;:=\;
\begin{tikzpicture}[anchorbase,scale=.2]
  \draw[blue, line width=.5mm] 
  (-1,2) to (0,1) to (1,2) 
  (0,0) to (0,1);
  \fill[blue] (0,0) circle (2.5mm);
\end{tikzpicture}\,
,\quad
\begin{tikzpicture}[anchorbase,scale=.2]
  \draw[blue, line width=.5mm] 
  (-1,-1) to (1,2) (1,-1) to (-1,2) ;
\end{tikzpicture}
\;:=\;
\begin{tikzpicture}[anchorbase,scale=.2]
  \draw[blue, line width=.5mm] 
  (-1,-1) to (0,0) to (1,-1) 
  (0,0) to (0,1)
  (-1,2) to (0,1) to (1,2) ;
\end{tikzpicture}
$$
These morphisms satisfy the expected string-straightening and vertex rotation relations.

\begin{definition}
In the following, we will write $B_I:=R\otimes_{R^I}\otimes R$ and call this a
\emph{generalised Bott--Samelson bimodule}. See also \cite{MR3502025}.
\end{definition}

\begin{lemma}
  \label{lem:hTrFE}
  With the same assumptions as in Example~\ref{exa:FrobExt1}, we have an isomorphism
  \[\Trzz(B_I) \cong \oplus_{|W_I|} \left(\Trzz(\one), \frac{[\Id_{B_I}]}{|W_I|}\right) \]
  in $\Kar\Trzz(\SBim(W))$, where we have used the identification
  \[\End_{\Trzz(\SBim(W))}(\Trzz(\one))\cong \HH_0(\SBim(W)) \] to describe the
   idempotents appearing on the right-hand side. 
\end{lemma}
\begin{proof}
The Frobenius extension provides an isomorphism $\phi\colon B_I\star B_I \cong \oplus_{|W_I|} B_I$, 
which admits a convenient diagrammatic description, c.f. \cite[Section 4]{MR3502025}. 
To make it explicit, let $\{x_\alpha\}$ and $\{y_\alpha\}$ denote dual bases of $R$ as a free $R^I$-module, 
where $\alpha \in W_I$. Then $\phi=\oplus_{\alpha \in W_I} \phi_\alpha$ with

$$\phi_\alpha\colon B_I\star B_I\to B_I, \quad \phi_\alpha( (r_1\star r_2)\star(r_3\star r_4)) 
:= \partial(x_\alpha r_2 r_3)  (r_1 \star r_4),
\quad \phi_\alpha = \begin{tikzpicture}[anchorbase,scale=.2]
  \draw[blue, line width=.5mm] 
  (-1,-1) to (-1,1) to (0,1.5) to (0,2)
  (1,-1) to (1,1) to (0,1.5);
  \node at (0,0) {\tiny$x_\alpha$};
\end{tikzpicture}$$

 and $\phi^{-1}=\sum_{\beta \in W_I} \phi^{-1}_\beta$ with

$$ \phi^{-1}_\beta \colon B_I \to B_I \star B_I, \quad \phi^{-1}_\beta(r_1 \star r_2) 
:= \sum_{\alpha\in W_I} (r_1 \star 1)\star(y_\beta \star r_2),
\quad \phi_\alpha = \begin{tikzpicture}[anchorbase,scale=.2]
  \draw[blue, line width=.5mm] 
  (-1,2) to (-1,0) to (0,-.5) to (0,-1)
  (1,2) to (1,0) to (0,-.5);
  \node at (0,1) {\tiny$y_\beta$};
\end{tikzpicture}
$$ 

Now in $\Kar\Trzz(\SBim(W)))$ we have:

\begin{align*}
\Trzz(B_I) 
&
\;=\; 
\begin{tikzpicture}[anchorbase,scale=.2]
  \draw[blue, line width=.5mm] 
  (-1,-2) to (-1,2);
\ann{3}
\end{tikzpicture}
\;=\;
\sum_{\alpha}
\begin{tikzpicture}[anchorbase,scale=.2]
  \draw[blue, line width=.5mm] 
  (-2,-2) to (-2,2)
  (-3,0) to (3,0);
  \node at (2,-1) {\tiny$x_\alpha$};
  \node at (2,1) {\tiny$y_\alpha$};
\ann{3}
\end{tikzpicture}
\;=\;
\sum_{\alpha,\beta}
\frac{1}{|W_I|}
\begin{tikzpicture}[anchorbase,scale=.2]
  \draw[blue, line width=.5mm] 
  (-2,-2) to (-2,2)
  (-3,0) to (-2,0) 
  (-2,.8) to (1,.8) to (2,0) to (3,0)
  (-2,-.8) to (1,-.8) to (2,0) ;
  \node at (2,-1.5) {\tiny$x_\alpha$};
  \node at (2,1.5) {\tiny$y_\alpha$};
  \node at (0,0) {\tiny$x_\beta y_\beta$};
\ann{3}
\end{tikzpicture}\\
&
\;=\;
\sum_{\alpha,\beta}
\frac{1}{|W_I|}
\begin{tikzpicture}[anchorbase,scale=.2]
  \draw[blue, line width=.5mm] 
  (-2,-2) to (-2,-.8) 
  (-2,.8) to (-2,2)
  (-1,-.8) to (-1,.8)
  (-3,.8) to (3,.8) 
  (-3,-.8) to (3,-.8);
  \node at (2,-1.5) {\tiny$x_\alpha$};
  \node at (2,1.5) {\tiny$y_\alpha$};
  \node at (0,0) {\tiny$ y_\beta$};
  \node at (-2,0) {\tiny$ x_\beta$};
\ann{3}
\end{tikzpicture}
\;=\;
\sum_{\alpha}
\frac{1}{|W_I|}
\begin{tikzpicture}[anchorbase,scale=.2]
  \draw[blue, line width=.5mm] 
  (-2,-2) to (-2,-.8) 
  (-2,.8) to (-2,2)
  (-3,.8) to (3,.8) 
  (-3,-.8) to (3,-.8);
  \node at (2,-1.5) {\tiny$x_\alpha$};
  \node at (2,1.5) {\tiny$y_\alpha$};
\ann{3}
\end{tikzpicture}
\end{align*}

Conversely we have:
$$
\frac{1}{|W_I|}
\begin{tikzpicture}[anchorbase,scale=.2]
  \draw[blue, line width=.5mm] 
  (-2,-1.2) to  (-2,1.2) 
  (-3,-1.2) to (3,-1.2) 
  (-3,1.2) to (3,1.2) ;
  \node at (0,0) {\tiny$x_\alpha y_\beta$};
\ann{3}
\end{tikzpicture}
\;=\;
\frac{1}{|W_I|}
\begin{tikzpicture}[anchorbase,scale=.2]
  \draw[blue, line width=.5mm] 
  (-2,-1.2) to  (-2,1.2) 
  (-3,0) to (-2,0) 
  (-2,-1.2) to (1,-1.2) to (2,0) to (3,0) 
  (-2,1.2) to (1,1.2) to (2,0) ;
  \node at (0,0) {\tiny$x_\alpha y_\beta$};
\ann{3}
\end{tikzpicture}
\;=\;
\frac{\delta_{\alpha, \beta} }{|W_I|}
\begin{tikzpicture}[anchorbase,scale=.2]
  \draw[blue, line width=.5mm] 
  (-3,0) to (3,0);
\ann{3}
\end{tikzpicture}
$$
This implies that $\frac{1}{|W_I|}
\begin{tikzpicture}[anchorbase,scale=.2]
  \draw[blue, line width=.5mm] 
  (-2,-1) to (-2,0.5) 
  (-3,0.5) to (3,0.5);
  \node at (2,-.25) {\tiny$x_\alpha$};
\annh{3}
\end{tikzpicture}$ 
and
$
\begin{tikzpicture}[anchorbase,scale=.2]
  \draw[blue, line width=.5mm] 
  (-2,1) to (-2,-.5) 
  (-3,-.5) to (3,-.5);
  \node at (2,.25) {\tiny$y_\alpha$};
\annh{3}
\end{tikzpicture}$
are the components of inverse isomorphisms 
\[\Trzz(B_I) \leftrightarrow \oplus_{|W_I|}  \left(\Trzz(\one), \frac{[\Id_{B_I}]}{|W_I|}\right).\] 
\end{proof}

\begin{remark}
An analogous argument in the derived setting shows that 
\[\hTr(B_I) \cong \oplus_{|W_I|} \left(\hTr(\one), \frac{[\Id_{B_I}]}{|W_I|}\right) \]
  in $\Kar\hTr(\SBim(W)))$. The main difference is that isotopies of diagrammatic 
  morphisms through the seam of the annulus 
  are now only possible only up to homotopy. One thus arrives at homotopy idempotents, 
  which split in the Karoubi envelope, as described in Section~\ref{sec:Karoubi}.
\end{remark}

We can describe the idempotent $\frac{[\Id_{B_I}]}{|W_I|}$ more explicitly. 
Indeed, it is well known that $B_I$ can be presented as a twisted complex consisting of 
$\Delta_w[\ell(w)]$ for $w\in W_I$, so by \eqref{eq:cone}
we get
\[
[\Id_{B_I}]=\sum_{w\in W_I}(-1)^{\ell(w)}[\Id_{\Delta_w}]=\sum_{w\in W_I} w
\]
where we identified $[\Id_{\Delta_w}]$ with $(-1)^{\ell(w)} w$ in $\HH_0(\SBim_W)=\HH_0(\CK^b(\SBim_W))$ 
using Theorem \ref{thm:HHSBim}. We get the following

\begin{corollary}
We have \[\hTr(B_I) \cong \oplus_{|W_I|} \left(\hTr(\one), \frac{1}{|W_I|}\sum_{w\in W_I}w\right). \]
\end{corollary}

\begin{example}
Let $W=S_n$ and $W_I=S_{k_1}\times\cdots \times S_{k_r}$ for $k_1+\ldots+k_r=n$. Then 
\[
\hTr(B_I)\simeq \bigoplus_{\lambda} (\hTr(\one),\ee_{\lambda})^{\oplus x_{\lambda,k_1,\dots,k_r}}
\]
where the sum is over all partitions $\lambda$ of $n$, $\ee_{\lambda}$ is an idempotent 
in $\C[S_n]$ corresponding to the irreducible representation $V_{\lambda}$ and the graded multiplicities are given by 
\[
  x_{\lambda,k_1,\dots,k_r}:= [k_1]!\cdots [k_r]!c_{k_1,\ldots,k_r}^{\lambda}  
\]
where $c_{k_1,\ldots,k_r}^{\lambda}$ denotes the multiplicity of the Schur function $s_{\lambda}$
in the product $h_{k_1}\cdots h_{k_r}$.
\end{example}

\subsection{Explicit annular simplification in type A}

Throughout this section, we work with type $A$ Soergel bimodules corresponding to the action of 
$W=S_n$ on $V=\C^n$.
We abbreviate the notation
$B_i:=B_{s_i}$ for the Bott-Samelson bimodules associated to the simple
reflections $s_i\in S_n$. Moreover, if $|i-j|=1$, we also consider the Soergel
bimodule $B_{iji}=B_{jij}:= R\otimes_{R^{\langle s_i,s_j\rangle }} R$. Further,
to declutter many expressions in this subsection, we will omit the $\star$
indicating the composition of bimodules.
 
\newcommand{\stari}{}

\begin{lemma}
\label{lem: last one}
Any Bott-Samelson bimodule for $S_n$ is isomorphic to a direct summand of a
direct sum of Bott-Samelson bimodules, in each of which $B_{n-1}$ appears at
most once.
\end{lemma}

\begin{proof}
This is a standard argument which uses $B_i\stari B_i \cong B_i(1)\oplus
B_i(-1)$, $B_i\stari B_j \cong B_j\stari B_i$ if $|i-j|>1$, and $B_i \stari B_j
\stari B_i \cong B_{iji} \oplus B_i$ if $|i-j|=1$. 
\end{proof}

For each subset $I=\{i_1<\cdots<i_k\}\subset S=\{1,\cdots, n-1\}$, we define 
the \emph{Coxeter-Bott-Samelson bimodule} $\underline{B}_I:=B_{i_1}\stari B_{i_2}\stari \cdots\stari B_{i_k}$.

\begin{lemma}
\label{lem:resolve by cox}
For any Soergel bimodule $B$ the trace $\Tr(B)$ is isomorphic to a direct summand in the 
direct sum of traces of Coxeter-Bott-Samelson bimodules as defined above. 
\end{lemma}

\begin{proof}
The proof is by induction on $n\geq 1$.  The base case is trivial.  
By Lemma \ref{lem: last one} we can present $B$ as a direct summand in the direct 
sum of objects of the form $X$ or $X\stari B_{n-1}\stari Y$, where
$X,Y\in \SBim_{n-1}$. 
The first case is taken care of by induction. In the second case we
can use Lemma \ref{lem:trace relation} to replace $\Tr(X\stari B_{n-1}\stari Y)$ by $\Tr(Y\stari
X\stari B_{n-1})$. Now by Lemma \ref{lem: last one} we can either write $Y\stari
X\subset_{\oplus}X'$ or $Y\stari X\subset_{\oplus}X'\stari B_{n-2}\stari Y'$ for some $X',Y'\in
\SBim_{n-2}$. Thus we either have $Y\stari X\stari B_{n-1}\subset_{\oplus} X'\stari B_{n-1}$ or
$$
\Tr(Y\stari X\stari B_{n-1})\subset_{\oplus} \Tr(X'\stari B_{n-2}\stari Y'\stari B_{n-1})\cong
\Tr(Y'\stari X'\stari B_{n-2}\stari B_{n-1}),
$$
which is taken care of by induction.
\end{proof}

For the following let $\lon_n$ denote the indecomposable Soergel bimodule corresponding to the
longest element in $S_n$.

\begin{lemma}
\label{lem: Xn relation}
In $\SBim_{n+1}$ we have:
$$
\lon_n\stari B_n\stari \lon_n\cong [n-1]!\lon_{n+1}\oplus [n]!\lon_{n}
$$
where we use quantum numbers to indicate multiple direct summands with grading shifts.
\end{lemma}

\begin{example}
For $n=2$ we have $X_2=B_1$, so $B_1\stari B_2\stari B_1\cong B_{121}\oplus B_{1}.$
\end{example}

\begin{proof}
Let $R=\C[X_1,\dots,X_{n+1}]$
\begin{align*}
    \lon_n\stari B_n\stari \lon_n 
&= R \otimes_{R^{S_n\times S_1}} R \otimes_{R^{S_1\times \cdots \times S_1 \times S_2}} R \otimes_{R^{S_n\times S_1}} R \\
& \cong [n-1]! R \otimes_{R^{S_n\times S_1}} R^{S_{n-1}\times S_1\times S_1} \otimes_{R^{S_{n-1} \times S_2}}
    R^{S_{n-1}\times S_1\times S_1} \otimes_{R^{S_n\times S_1}} R  \\
& \cong [n-1]! R \otimes_{R^{S_{n+1}}} R \oplus [n-1]![n] R \otimes_{R^{S_{n}\times S_1}} R  \\
& = [n-1]!\lon_{n+1}\oplus [n]!\lon_{n}
\end{align*}
Here we have used the well-known ``square-switch'' isomorphism of
singular Bott-Samelson bimodules to proceed to the third line.
This can, for example, be deduced from \cite[Lemma 11.2]{MR3234803}---
the corresponding statement for matrix factorizations---by taking homology with respect to the positive differential 
and forgetting the negative differential, in the terminology of \emph{loc.~cit.}.
\end{proof}

In the previous proof we have used the fact that $\C[X_1,\dots, X_{n-1}]$ is a free module of (graded) rank $[n-1]!$ 
of $\C[X_1,\dots, X_{n-1}]^{S_{n-1}}$ (in fact, a Frobenius extension).

\begin{lemma}
\label{lem: cox to parabolic}
In $\Kar\hTr(\SBim_n)$  we have that $\Tr(B_1\stari \cdots\stari B_{n-1})$ is isomorphic to
a summand in a direct sum of traces of Soergel bimodules corresponding to the longest
elements in parabolic subgroups in $S_n$.
\end{lemma}

\begin{proof} 
We prove a more general statement: $\Tr(\lon_{k} B_{k}\cdots B_{n-1})$  is
isomorphic to a summand in a direct sum of Soergel bimodules corresponding to
the longest elements in parabolic subgroups.  Let us use induction in $k$, starting from $k=n$
where we have just $\lon_n$. For the induction step $k+1\mapsto k\geq 1$ we have
$$
\lon_k\stari \lon_k\stari B_{k}\cdots B_{n-1}\cong [k]!\lon_k\stari B_{k}\stari \cdots \stari  B_{n-1}
$$
and
$$
\Tr(\lon_k\stari \lon_k\stari B_{k}\cdots B_{n-1})\cong\Tr(\lon_kB_{k}\cdots B_{n-1}\stari \lon_k)
\cong\Tr(\lon_k\stari B_k\stari \lon_k\stari B_{k+1}\stari \cdots \stari B_{n-1}),
$$
and by Lemma \ref{lem: Xn relation} we have
$$
\lon_k\stari B_k\stari \lon_k\stari B_{k+1}\stari \cdots \stari B_{n-1}
\cong[k-1]!\lon_{k+1}\stari B_{k+1}\stari \cdots \stari B_{n-1}
\oplus [k]!\lon_{k}\stari B_{k+1}\stari \cdots \stari B_{n-1}.
$$
The claimed statement now follows from the case $k=2$.
\end{proof}

Let us give an example computation to highlight the differences compared to the
annular Khovanov--Rozansky invariant of Queffelec--Rose \cite{MR3729501}

\begin{example}[The full twist on two strands]
\label{ex: full twist derived}
  Consider the category of Soergel bimodules of type $A_1$. 
It is well-known that the Rouquier complex of the full twist on two strands can be expressed as
\[ F(\sigma^2) \simeq \left( \underline{B(-1)}\xrightarrow{x^l_1-x^r_1} B(1)
\xrightarrow{\text{unzip}} R(2) \right) \] To compute the derived horizontal
trace class of this complex we will use the following tools. First of all, just
as in the underived horizontal trace, we have homotopy equivalences
$\hTr(B)\simeq \wedge^2(-1)\oplus \wedge^2(1)$ and $\hTr(R)\cong \wedge^2 \oplus
S^2$, where $\wedge^2$ and $S^2$ indicate the isotypic components of $\hTr(R)$
under the natural $S_2$ action. Furthermore we now observe $x^l_1-x^r_1 = d(\bbar
x_1 \bbar \Id_B) = d(w(x_1;B))$ as well as $\textrm{unzip} \circ w(x_1;B)  =
\theta_1\circ \textrm{unzip} $. Then the derived annular simplification proceeds as follows:
\[
\begin{tikzpicture}[baseline=0em]
\tikzstyle{every node}=[font=\scriptsize]
\node (a) at (0,0) {$\Tr(F(\sigma^2))$};
\node at (1.5,0) {$\simeq $};
\node (b) at (3,0) {$\Big(\underline{\Tr(B(-1))}$};
\node (c) at (7,0) {$\Tr(B(1))$};
\node (d) at (11,0) {$\Tr(R(2)) \Big)$};
\node at (1.5,-2) {$\simeq $};
\node (ab) at (3,-2) {$\Big(\underline{\Tr(B(-1))}$};
\node (ac) at (7,-2) {$\Tr(B(1))$};
\node (ad) at (11,-2) {$\Tr(R(2)) \Big)$};
\node at (1.5,-4) {$\simeq $};
\node (bb) at (3,-4) {$\Big(\underline{\wedge^2(-2)\oplus\wedge^2}$};
\node (bc) at (7,-4) {$\wedge^2 \oplus \wedge^2(2)$};
\node (bd) at (11,-4) {$S^2(2)\oplus \wedge^2(2)\Big)$};
\node at (1.5,-6) {$\simeq $};
\node (cb) at (3.3,-6) {$ \underline{\wedge^2(-2)} \;\; \oplus\;\; \Big(\underline{\wedge^2}$};
\node (cc) at (7,-6) {$\wedge^2$};
\node (cd) at (11,-6) {$S^2(2)\Big)$};
\draw[frontline,->,>=stealth,shorten >=1pt,auto,below,node distance=1.8cm,thick]
(ab) to[bend right=10] node {$-\text{unzip} \circ w(x_1;B)$} (ad)
(bb) to[bend right=10] node {$\sqmatrix{x\theta & \theta \\ x \theta & \theta}$} (bd)
(cb) to[bend right=10] node {$\theta$} (cd);
\path[->,>=stealth,shorten >=1pt,auto,node distance=1.8cm,
  thick]
(b) edge node  {$x^l_1-x^r_1$} (c)
(c) edge node {$\text{unzip}$} (d)
(ac) edge node {$\text{unzip}$} (ad)
(bc) edge node {$\sqmatrix{x & 0 \\ x & \Id} $} (bd)
(cc) edge node {$x$} (cd);
\path[->,>=stealth,shorten >=1pt,auto,node distance=1.8cm,
  thick, dashed]
(b) edge node {$\Id$} (ab)
(c) edge node {$\Id$} (ac)
(d) edge node {$\Id$} (ad)
(b) edge node {$w(x_1;B)$} (ac);
\end{tikzpicture}.
\]
Where we abbreviate $x:=(x_1-x_2)/2$ and $\theta:=(\theta_2-\theta_1)/2$. Since
we work over $\C$, we may also rescale to get $x=x_1-x_2$ and
$\theta=\theta_1-\theta_2$. (Note that the second line indeed shows a twisted
complex, and the dashed arrows encode a degree zero closed invertible morphism
between twisted complexes.) In the underived horizontal trace the long arrow
would be zero, and the complex would split into three direct summands.
\end{example}

\subsection{Proof of Theorem \ref{th: annular simplification}}
\label{sec:proof annular simp}

In this section we prove Theorem \ref{th: annular simplification}. First we
specialize some general facts about the $\Ai$-category $\Perf(A)$ proved in section
\ref{ss:perf definitions} to the case of $A:=R\otimes\Lambda\rtimes \C[S_n]$ with maps $\mu_d$ induced by the retraction from the cyclic bar complex.

Recall that the generators $\theta_i$ of $\Lambda$ have cohomological degree
$-1$, so $A$ is indeed supported in nonpositive cohomological degrees. Since $A$
has zero differential, we have $\Kar^{dg}\langle A\rangle =\Kar\langle
A\rangle$ where $\langle A\rangle$ is the category of free $A$ modules of finite rank. Because we work
over a field of characteristic zero, $\Kar\langle A\rangle$ is a semisimple
category with finitely many indecomposable objects labeled by irreducible
representations of $S_n$. The category $\Perf(A)=\pretr{\Kar\langle A\rangle}$
consists of twisted complexes built out of these objects. By Theorem \ref{th: perf complete}
the category $\Perf(A)$ is homotopy idempotent complete.

\begin{proof}[Proof of Theorem \ref{th: annular simplification}]
The quasi-equivalence between $\hTr(\SBim_n)$ and $\Perf(A)$ is constructed in two steps. 

First, recall that the endomorphism algebra of $\Tr(\one)$ can be identified with the cyclic bar complex
$\CCC(\SBim_n)$. We have therefore a functor 
\begin{equation}
\label{eq: htr to cc}
\Hom_{\hTr}(\Tr(\one),-)\colon \hTr(\SBim_n)\to {\rmod}\CCC(\SBim_n)
\end{equation}
By Lemma \ref{lem:resolve by cox} we can resolve any object in $\hTr(\SBim_n)$
by the traces of Coxeter-Bott-Samelson bimodules. By Lemma \ref{lem: cox to
parabolic} the trace of every such bimodule is equivalent to a summand in the
sum of traces of indecomposable Soergel bimodules corresponding to the longest
elements of parabolic subgroups. Finally, by Lemma \ref{lem:hTrFE} the trace of
any such bimodule is equivalent to a summand in the direct sum of several copies
of the trace of the identity bimodule. 

This means that any object in the essential image of \eqref{eq: htr to cc} is
homotopy equivalent to a direct summand in a free $\CCC(\SBim_n)$-module, and
\eqref{eq: htr to cc} defines a quasi-fully faithful quasi-functor\footnote{Note
that the category of $\CCC(\SBim_n)$-modules homotopy equivalent to direct
summands of free modules is quasi-equivalent to
$\Kar^{dg}\langle\CCC(\SBim_n)\rangle$ }
$$
\Hom_{\hTr}(\Tr(\one),-)\colon \hTr(\SBim_n)\to \Kar^{dg}\langle\CCC(\SBim_n)\rangle
$$
The corresponding functor
\begin{align}
\label{eq: htr to perf}
\Kar^{dg}(\pretr{\hTr(\SBim_n)}) & \to \Kar^{dg}(\pretr{\Kar^{dg}\langle\CCC(\SBim_n)\rangle)})
\\
& =
\Kar^{dg}(\Perf(\CCC(\SBim_n)))
\end{align}
obtained by successively applying $\pretrf$ and $\Kar^{dg}$ on both sides is then a quasi-equivalence. 

For the second step, observe that by Theorem \ref{thm:HHSBim} the dg algebra
$\CCC(\SBim_n)$ deformation retracts onto $A$. 
As in Lemma \ref{lem:perf quasi equivalent}, we get quasi-equivalences $\langle
\CCC(\SBim_n) \rangle\simeq \langle A\rangle$ and
$$
\Kar^{dg}(\Perf(\CCC(\SBim_n)))\simeq \Kar^{dg}(\Perf(A)).
$$
Since $\Perf(A)$ is Karoubian, $\Kar^{dg}(\Perf(A))\simeq \Perf(A)$. 

By combining \eqref{eq: htr to perf} with all these quasi-equivalences, we conclude that 
$$
\Kar^{dg}(\pretr{\hTr(\SBim_n)})\simeq \Perf(A).\vspace{-.6cm}
$$ 
\end{proof}

\begin{corollary}
We have the following equivalence of categories:
$$
\Kar(\pretr{\Trzz(\SBim_n)})\simeq \Ch^b(\Kar(\langle A\rangle))\simeq \Perf(R\rtimes \C[S_n]).
$$
In particular, after taking the pretriangulated hull and Karoubi completion the functor
$\hTr(\SBim_n)\to \Trzz(\SBim_n)$ can be identified with the
forgetful functor $\epsilon\colon \Perf(A)\to \Ch^b(A)$ defined in section
\ref{ss:perf definitions}.
\end{corollary}

\begin{proof}
The algebra $A$ is supported in nonpositive cohomological degrees, with
$R\rtimes \C[S_n]$ in cohomological degree zero. Therefore 
$$
\Ch^b(\Kar(\langle A\rangle))\simeq \Ch^b(\Kar(R\rtimes \C[S_n]))=\Perf(R\rtimes \C[S_n]).\vspace{-.6cm}
$$
\end{proof}

\subsection{A derived annular Khovanov-Rozansky invariant}
\label{sec:derived-AKhR}

In \cite{MR3729501,GW} the (underived) traces of web categories
were related to annular Khovanov-Rozansky invariants, and to the Khovanov-Rozansky
homology of links in $\mathbb{R}^3$. In this subsection we review this construction using the
category $\SBim_n$ and its trace.

Given a braid word $\underline{\b}$ on $n$ strands, let $F(\underline{\b})$ again denote the
Rouquier complex for $\underline{\b}$ as defined in Section \ref{subsec: SBim def}, which
we now consider as an object in $K^b(\SBim_n)$, the bounded homotopy category of
$\SBim_n$. We have already mentioned that these complexes satisfy braid relations up to
(canonical) homotopy equivalence. In fact, braid cobordisms induce natural chain maps (up to
homotopy) between Rouquier complexes \cite{MR2721032}.

Consider the underived trace functor $\Trzz\colon \SBim_n \to \Trzz(\SBim_n)$,
with the target considered as embedded in the Karoubi completion
$\Kar(\Trzz(\SBim_n))$. Recall that the latter is equivalent to the category of
graded projective $R\rtimes \C[S_n]$-modules. In the following definition we use the
functor $\Trzz$, extended to the bounded homotopy categories of the source and target.

\begin{definition}
  The \emph{annular Khovanov--Rozansky invariant} of a braid word
  $\underline{\b}$ on $n$ strands is defined as
  \[\AKhR(\hat{\underline{\b}}):=\Hom_{\Trzz}(\Trzz(\one), \Trzz(F(\underline{\b})))\in
  \Perf(R\rtimes \C[S_n])\] 
\end{definition}

By virtue of factoring through the underived trace, $\AKhR$ is a categorical
invariant of braid conjugacy classes (a.k.a. annular links with a coherent
orientation) which is natural under annular link cobordisms (preserving the
coherent orientation). More precisely, $\AKhR$ is defined on $n$-strand braid
words for each $n\geq 0$ separately, but these invariants fit together to give a
\emph{monoidal} annular link invariant, see \cite{GW} for details. Relatives of
this notion of annular Khovanov--Rozansky invariant have previously been studied
in \cite{MR3729501,QRS}. Another ahistorical aspect of our presentation here is that
the interest in \emph{annular} Khovanov--Rozansky invariants rose well after the
construction of the triply-graded Khovanov--Rozansky homology of links in
$\mathbb{R}^3$ \cite{MR2421131, MR2339573}, which categorifies the HOMFLY-PT
polynomial. We will comment on their relationship in the next section.

We are now ready to define a derived annular Khovanov--Rozansky invariant.

\begin{definition} The \emph{derived annular Khovanov--Rozansky invariant} is defined on
closures of $n$-strand braid words $\underline{\b}$ as 
  \[\AKhR_{\mathrm{dg}}(\hat{\underline{\b}}):=\Hom_{\Tr}(\hTr(\one), \hTr(F(\underline{\b})))\in
  \Perf(R\otimes \Lambda\rtimes \C[S_n])\] 
  where $F(\underline{\b})$ is the Rouquier complex
  of $\underline{\b}$, $\hTr$ is the universal dg monoidal trace, and we use that
  $\Hom_{\Tr}(\hTr(\one),-)$ realises the equivalence from Theorem~\ref{th: annular
  simplification}. As for the underived annular Khovanov--Rozansky invariant, this
  construction is functorial under braid-like annular link cobordisms between braid
  closures up to homotopy.
\end{definition}

One important feature of $\AKhR_{\mathrm{dg}}(\hat{\underline{\b}})$ is that the $2\pi$
rotation of the annular link $\hat{\underline{\b}}$ typically induces an interesting
endomorphism of $\AKhR_{\mathrm{dg}}(\hat{\underline{\b}})$, while it always induces the identity map on 
$\KhR(\hat{\underline{\b}})$ and its annular version.

Recall that by \eqref{eq:connes for SBim} the Connes differential $\B=\sum
\theta_i\frac{\partial }{\partial x_i}$ defines a derivation on the algebra
$\HH_{\subb}(\SBim_n)\cong R\otimes \Lambda\rtimes \C[S_n]$.  This allows us to define an
interesting endofunctor on the category of twisted $R\otimes \Lambda\rtimes
\C[S_n]$-modules.    Given a  twisted complex $(A,d)$ we have:
$$
0=\B(d^2)=\B(d)d+d\B(d),
$$ 
so that $\B(d)$ is always a closed (degree 1) endomorphism of $(A,d)$.

\begin{conjecture}
The action of the rotator $w_X$ (defined in Remark \ref{rem:rotator}) on twisted
complexes $X$ built out of summands of $\Tr(\one)$ is homotopic to $\Id+\B+\text{higher order terms}$,
where $\B$ is the action of Connes differential defined above.
\end{conjecture}

\begin{example} Let $f$ be an arbitrary polynomial in $R=\End_{\SBim_n}(\one)$,
let $X=\Cone(f)$. Since the rotator $w_{\one}:\Tr(\one)\to \Tr(\one)$ is trivial,  
by \eqref{eq:traciator two term} the rotator $w_X$ is given by the morphism
$$
\begin{tikzcd}
  \big( \Tr(\one) \arrow{r}{f} \arrow{d}{\Id} \arrow{dr}{w(f;\one)} & \Tr(\one)\big)  \arrow{d}{\Id}\\
  \big( \Tr(\one) \arrow{r}{f}   & \Tr(\one) \big) \\
\end{tikzcd}
$$
It is easy to see that $w(f;\one)=\bbar f\bbar \Id$ is homotopic to $\B(f)$.

\end{example}

Another important feature is that the derived annular Khovanov--Rozansky
invariant of an annular link has an action of the derived center
$\ZC^{\text{dg}}(\SBim_n)$. The second author has shown with Ben Elias that the
Rouquier complex of the full twist braid $\FT_n$ (together with suitable
half-braiding data) is an object of the derived Drinfeld center \cite{EHcenter}.  The corresponding
endofunctor of the derived trace sends $\AKhR_{\mathrm{dg}}(\hat{\underline{\b}})$ to
$\AKhR_{\mathrm{dg}}(\widehat{FT_n \underline{\b}})$, i.e. it cuts the thickened annulus
containing the annular link $\hat{\underline{\b}}$ and re-glues it after a $2\pi$ twist to
create $\widehat{FT_n \underline{\b}}$. 

This operation suggests that $\AKhR_{\mathrm{dg}}$ should be considered as an
invariant of (coherently oriented) links in $S^1\times D^2$, which can be
computed by choosing an $I$-bundle structure on $S^1\times D^2$, but which comes
with the data necessary to change $I$-bundle structure.

\subsection{Triply graded homology and the Hochschild cohomology of Soergel bimodules}
In this section we focus on the Hochschild cohomology of individual Soergel
bimodules $M$ in $\SBim(W)$ (we will soon specialise to $\SBim_n$), not on the
Hochschild (co)homology of the category $\SBim_n$. 

We define
$\HH^i(M):=\Ext^i(R,M)$ where the $\Ext$ groups are computed in the category of
$R,R$-bimodules. In particular, $\HH^0(M)=\Hom(R,M)$. 

Given a complex $C=[\ldots \to M_{j}\to M_{j+1}\to \ldots]$ of Soergel
bimodules, we define a complex of graded vector spaces
$$
\HH^i(C):=[\ldots \to \HH^i(M_{j})\to \HH^i(M_{j+1})\to \ldots].
$$

\begin{definition}
Given a braid word $\underline{\b}$ on $n$ strands, the triply-graded {\em
Khovanov-Rozansky homology} of the braid closure $\hat{\beta}$ is defined as
$$
\KhR(\hat{\beta}):= H^{\supb}(\oplus_{i}\HH^i(F(\underline{\b}))),
$$
where $F(\underline{\b})$ is the associated Rouquier complex of Soergel
bimodules in $\SBim_n$ defined in Section \ref{subsec: SBim def}. The vector
space $\KhR(\hat{\beta})$ is triply-graded by $q$-degree, cohomological degree
in the complex $F(\underline{\b})$, and Hochschild degree $i$. 
\end{definition}
As defined, $\KhR$ is a braid conjugacy invariant, i.e. an invariant of links in
a thickened annulus that are obtained from braid closures. However, after an
overall grading shift, $\KhR$ becomes invariant under the second Markov move,
and thus an invariant of links in $\mathbb{R}^3$. Since braid cobordisms induce
natural chain maps (up to homotopy) between Rouquier complexes, one also has
induced morphisms between the Khovanov-Rozansky homologies. If we also take into
account the morphisms induced by braid conjugation on the level of Hochschild
homology, this can be summarised by saying that $\KhR$ is functorial under braidlike
annular link cobordisms between braid closures.

Next we explain how the triply-graded Khovanov--Rozansky homology can be
recovered from the annular Khovanov--Rozansky invariant (derived or underived).
For an alternative but related approach see \cite{QRS}. 

It is well-known that for arbitrary Soergel
bimodules $M$ and $N$ there is a natural isomorphism
\begin{equation}
\label{eq: hh as trace}
\HH^i(M\otimes N)\cong \HH^i(N\otimes M)
\end{equation}
of $q$-graded vector spaces. In other words, $\HH^i$ is a trace-like linear
functor and hence factors through $\Trzz$. Moreover, in type A Rasmussen proved
\cite[Proposition 4.6]{MR3447099} that $\HH^i(M)$ is free over $R$. (More
generally this is known for all Weyl groups, see Webster--Williamson
\cite{MR2852120}).

\begin{lemma}
\label{lem: hh from hooks}
For a Soergel bimodule $M$ in $\SBim_n$ we have
$$
\HH^i(M)\cong \Hom_{\Trzz}(\ee_{i,1^{n-i}}\Trzz(\one),\Trzz(M))\oplus \Hom_{\Trzz}(\ee_{i-1,1^{n-i+1}}\Trzz(\one),\Trzz(M)).
$$
where $\ee_{i,1^{n-i}}$ is the idempotent in $\C[S_n]$ corresponding to the hook
partition $(i,1^{n-i})$. In particular,
$$
\Hom(\one,M)\cong\HH^0(M)\cong\Hom_{\Trzz}(\ee_{-}\Trzz(\one),\Trzz(M))
$$
where $\ee_{-}=\ee_{1^n}$ is the antisymmetrizer in $\C[S_n]$. The isomorphisms
above are natural in $M$.
\end{lemma}

\begin{proof}
By \eqref{eq: hh as trace} the functor $\HH^i$ is a trace-like functor from
$\SBim_n$ to graded vector spaces, so as discussed in
Section~\ref{section:trace-like} it factors through $\hTr(\SBim_n)$ and defines
a functor $\HH^i\colon \hTr(\SBim_n)\to \Vect$. By Theorem \ref{th: annular
simplification} we have that $\Tr(M)$ is isomorphic to the direct sum of direct
summands of $\Tr(\one)$, so it is sufficient to check it for $M=\Tr(\one)$. Now 
$$
\HH^i(R)\cong R\otimes \wedge^i(\C^n) \cong R\otimes (V_{i,1^{n-i}}\oplus V_{i-1,1^{n-i+1}}).
$$
while 
$$
\Hom_{\Trzz}(\Trzz(\one),\Trzz(\one))\ee_{\lambda}
\cong (R\rtimes \C[S_n]) \ee_{\lambda} \cong R\otimes V_{\lambda}
$$
for any irreducible representation $V_{\lambda}$ of $S_n$.
\end{proof}

\begin{theorem}
\label{thm: hh from hooks complex}
For any complex $C$ of Soergel bimodules one has the following isomorphism of
complexes of graded vector spaces:
$$
\HH^i(C)\cong \Hom_{\Trzz}(\ee_{i,1^{n-i}}\Trzz(\one),\Trzz(M))\oplus \Hom_{\Trzz}(\ee_{i-1,1^{n-i+1}}\Trzz(\one),\Trzz(C)).
$$
\end{theorem}

\begin{proof}
The functor $\Trzz$ sends $C$ to
$$
\Trzz(C)=[\ldots \to \Trzz(M_{j})\to \Trzz(M_{j+1})\to \ldots]
$$
with no higher differentials. Now the statement follows from Lemma \ref{lem: hh from hooks}.
\end{proof}

 Theorem \ref{thm: hh from hooks complex}  implies that Khovanov-Rozansky homology
can be computed by first applying the functor $\Trzz$ to $F(\underline{\b})$, then
running the annular simplification from Theorem \ref{th: annular
simplification}.

\begin{proposition}
For a braid word $\underline{\b}$ on $n$ strands we have:
\[
  \KhR(\hat{\beta}) \cong \oplus_{i} \Hom_{\Trzz}((\ee_{i,1^{n-i}}+\ee_{i-1,1^{n-i+1}})\Trzz(\one),\Trzz(F(\underline{\b})))
  \]
\end{proposition}

\begin{remark}
In \cite{GW} the first and third authors defined the {\em evaluation} functor $\Trzz(\SBim_n)\to \Vect$ which sends 
$\Trzz(\one)$ to $R$ as a module over $R\rtimes \C[S_n]$. This is equivalent to the above since
\[
\Hom_{\Trzz}(\ee_{-}\Trzz(\one),\Trzz(\one))=(R\rtimes \C[S_n])\ee_{-} =R.
\]
Similarly, higher Khovanov-Rozansky homology $\HH^i$ can be computed by
evaluating $\hTr(\one)$ to $\HH^i(R)=R\otimes \wedge^i(\C^n)$ as a module over
$R\rtimes \C[S_n]$.
\end{remark}


\end{document}